%% file: diss.tex
\documentclass[12pt, a4paper, titlepage]{book}

\usepackage{amsmath, amsthm, amsfonts, amssymb, mathbbol}
\usepackage{ifthen, enumerate, url}
\usepackage{pgf}
\usepackage{color}
\usepackage{appendix}
\usepackage[pdfborder={0 0 0},draft=false,breaklinks=true]{hyperref}

\newtheorem{thm}{Theorem}
\newtheorem{cor}[thm]{Corollary}
\newtheorem{prop}[thm]{Proposition}
\newtheorem{lemma}[thm]{Lemma}
\newtheorem{exa}[thm]{Example}

\newtheorem{ass}[thm]{Assumptions}
\newtheorem{defin}[thm]{Definition}
\newtheorem{rem}[thm]{Remark}
\newtheorem{rems}[thm]{Remarks}

\theoremstyle{definition}

\DeclareMathOperator{\mv}{\mathsf{v}}
\DeclareMathOperator{\me}{\mathsf{e}}
\DeclareMathOperator{\Li}{\mathcal{L}}
\DeclareMathOperator{\E}{\mathsf{E}}
\DeclareMathOperator{\V}{\mathsf{V}}
\DeclareMathOperator{\G}{\mathsf{G}}
\DeclareMathOperator{\pir}{\pi^{-1,r}_{i,0}}
\DeclareMathOperator{\pjr}{\pi^{-1,r}_{j,0}}
\DeclareMathOperator{\re}{{\rm Re}}
\DeclareMathOperator{\range}{{\rm Range}}
\DeclareMathOperator{\Ker}{{\rm Ker}}
\DeclareMathOperator{\sign}{{\rm sign}}

\DeclareMathOperator{\dom}{{\rm Dom}}
\DeclareMathOperator{\supp}{{\rm supp}}
\DeclareMathOperator{\im}{{\rm Im}}
\DeclareMathOperator{\id}{{\rm id}}
\DeclareMathOperator{\etasg}{(e^{ta})_{t \geq 0}}
\DeclareMathOperator{\etbsg}{(e^{tb})_{t \geq 0}}

\DeclareMathOperator{\Sfin}{S_{\rm fin}}
\DeclareMathOperator{\Vfin}{\mathsf{V}_{\rm fin}}
\DeclareMathOperator{\Gfin}{\mathsf{G}_{\rm fin}}
\DeclareMathOperator{\Efin}{\mathsf{E}_{\rm fin}}
\DeclareMathOperator{\dv}{d_{\mathsf{V}}}
\DeclareMathOperator{\de}{d_{\mathsf{E}}}
\DeclareMathOperator{\It}{\tilde{\mathcal I}}
\DeclareMathOperator{\K}{\begin{pmatrix} P_Y &0 \\ 0 & P_Y \end{pmatrix} }
\begin{document}
%\maketitle
\setcounter{tocdepth}{2}

\input{titlepage}

\newpage
\vspace*{12cm}
{\large
\begin{tabbing}
Dekan:\hspace{3cm} \=Prof.\,Frank Stehling\\[1cm]
\\
Erstgutachter:\>Prof.\,Wolfgang Arendt\\[0.4cm]
Zweitgutachter:\>Prof.\,Frank Steiner\\[0.4cm]
Drittgutachter:\>Prof.\,Joachim von Below\\[0.4cm]
\\
\\
Tag der Promotion:\>13. Juni 2008
\end{tabbing}
}
\pagestyle{plain}
\pagenumbering{roman}
\tableofcontents

\newpage
\addcontentsline{toc}{chapter}{Introduction}
\input{introduction}

\chapter{Sesquilinear forms on product Hilbert spaces}\label{hesse}
This Chapter is devoted to the study of sesquilinear forms that are defined on the direct sum oh Hilbert spaces.
We start proving general properties of these forms (Sections~\ref{sec:findimarg}-\ref{sec:gencoercivity}) 
and then we move to the investigation of properties of the associated evolution equation (Sections~\ref{sec:operator}-\ref{sec:symmetriesmatrix}).
Two illustrative applications are also discussed (Sections~\ref{sec:wave}-\ref{sec:dynamic}). Finally, we devote a part of the Chapter (Sections~\ref{sec:nondiag} and~\ref{sec:histremI}) to general considerations.
\input{sec_findimarg}

\input{sec_systemsI}
\input{sec_gencoercivity}
\input{sec_operator}

\input{sec_evoleq}

\input{sec_symmetriesmatrix}
\input{sec_wave}
\input{sec_dynamic}
\input{sec_nondiag}
\newpage
\input{sec_histremI}
\chapter{Network equations}\label{network}

In this chapter we will describe and investigate the theory of parabolic equations on network-shaped structures. To be more precise, we want to consider differential equations taking place on the edges of a network. 

Our main goal is to prove relations between graph theoretic properties and features of the solutions of strongly coupled diffusion equations on networks.
To this end, we start discussing basic graph theoretic issues in the first section.

\input{sec_infintro}
\input{sec_graphoper}

\input{sec_operatordomain}
\input{sec_systemsII}
\input{sec_infinitesymmetries}
\input{sec_treelayer}
\input{sec_general}
\newpage
\input{sec_histremII}

%\chapter{Applications}
%\section{General boundary coupling}
%\section{Nonautonomous and nonlinear equations}
%\section{Relation to quantum mechanics}
%\section{Historical remarks}

\begin{appendix}
\chapter{Appendix}
\input{sec_sesquilinear}
%\newpage
%\input{sec_intermezzo}
\end{appendix}
\newpage
\addcontentsline{toc}{chapter}{Bibliography}
\bibliographystyle{plain} 
\bibliography{literatur}

\newpage
\addcontentsline{toc}{chapter}{Zusammenfassung in deutscher Sprache}
\input{sec_zusammen}

\newpage
\vspace*{15cm}
\noindent
{\bf Erkl\"arung:}\\ \vspace{0.1cm}

\noindent
Die vorliegende Arbeit habe ich selbst\"andig angefertigt und keine anderen als die angegebenen Quellen und Hilfmittel benutzt sowie die w\"ortlich oder inhaltlich
\"ubernommene Stellen als solche erkenntlich gemacht.\vspace{1cm}

%\hfill
%(Stefano Cardanobile)
\newpage
\input{sec_leben}
\newpage
\input{sec_publikationen}
\end{document}

%% file: titlepage.tex
\thispagestyle{empty}
\vspace{2cm}
\begin{center}
\Huge 
\sc 
Diffusion systems\\ 
and\\
heat equations on networks\\
\vspace{4cm}
\Large
Dissertation\\
\vspace{1cm}
\rm
\large
der Fakult\"at f\"ur Mathematik und Wirtschaftswissenschaften\\
der Universit\"at Ulm\\
zur Erlangung des Grades eines\\ 
Doktors der Naturwissenschaften\\
\vspace{5cm}
vorgelegt von\\
\sc
Stefano Cardanobile\\
\rm
aus Bari\\
\vfill
2008
\end{center}

%% file: introduction.tex
\chapter*{Introduction}
Sesquilinear forms have played an important role in the analysis of elliptic equations for a long time.
A possible starting point can be seen in Bernhard Riemann's doctoral thesis~\cite[Art.~18]{Rie51}, where the celebrated Dirichlet's Principle is formulated.
In short, the latter states that given an open, bounded domain $\Omega$, the solution of the Laplace equation $\Delta \psi(x)=0$ with homogeneous boundary conditions
$\psi_{|\partial \Omega}=0$ 
is the function for which the Dirichlet integral
\begin{equation}\label{intro0}\tag{DI}
\int_\Omega |\nabla \psi(x)|^2 dx
\end{equation}
attains a minimum in the set of the once differentiable functions that vanish on the boundary of $\Omega$.

The Dirichlet integral of $\psi$ can be seen as the value in the point $\psi$ of the quadratic form $a: H^1_0(\Omega) \to \mathbb C$ defined
by
$$
a(\psi):=\int_\Omega \big( \nabla \psi(x) \mid \overline{\nabla \psi(x)} \big) dx,\quad  \psi \in H^1_0(\Omega).
$$
An important question is whether it is possible to connect properties of the Laplace equation and of its quadratic form with properties of the heat equation
\begin{equation}\tag{HE}
\left\{\begin{array}{rcll}
\frac{\partial }{\partial t}\psi(t,x)&=& \Delta \psi(t,x), & t\geq 0,\\
\psi(t,x) &=&0,& x \in \partial \Omega,\\
\psi(0,x)&=& f(x)& x \in \Omega.
\end{array}\right.
\end{equation}

After the second world war, Tosio Kato (\cite{Kat55}, \cite[Chap.~6]{Kat66}) recognised that for this purpose it is convenient
to consider the \emph{sesquilinear form} $a: H^1_0(\Omega) \times H^1_0(\Omega) \to \mathbb C$ defined by
$$
a(\psi,\psi'):=\int_\Omega \big( \nabla \psi(x) \mid \overline{\nabla \psi'(x)} \big) dx, \quad \psi,\psi' \in H^1_0(\Omega).
$$
He also formulated the general concept of sesquilinear forms and constructed an autonomous theory for such mappings.

We shortly explain the basic idea, in a slightly different manner from that introduced by Tosio Kato. 
Having fixed a complex Hilbert space $V$, 
a sesquilinear form $(a,V)$ is a mapping $a: V \times V \to \mathbb C$ which is linear in the first component and antilinear in the second one.
If the mapping $a:V \times V \to \mathbb C$ is continuous and a second Hilbert space $H$ is given, such that $V$ is densely embedded in $H$,
it is possible to canonically associate an operator $(A,D(A))$ to the form $(a,V)$ by a result due to Peter Lax and Arthur Milgram, see~\cite{LaxMil54}.
Under suitable assumptions on the form $(a,V)$ there exists a family of bounded operators $(S_t)_{t \geq 0}$
such that $(S_t f)_{t \geq 0} $ is the solution of the abstract Cauchy problem
\begin{equation}\tag{ACP}
\left\{\begin{array}{rcll}
\frac{d }{d t}\psi(t)&=& A\psi(t), & t\geq 0,\\
\psi(0)&=& f,& f \in H,
\end{array}\right.
\end{equation}
where $(A,D(A))$ is the operator associated with the form.
For this reason, we write $\etasg:=(e^{tA})_{t \geq 0}:=(S_t)_{t \geq 0}$. 
Once this family has been obtained,  
it is possible to show that some properties of the family of operators $\etasg$ can be characterised in terms of properties of the sesquilinear form $(a,V)$.

After the fundamental work of Tosio Kato, the use of sesquilinear forms for the analysis of partial differential equations has developed autonomously. 
Among other important results, a striking theorem was proved by El-Maati Ouhabaz in \cite{Ouh96}. There, a closed convex set $C \subset H$ is considered and the property 
$$
\psi \in C \Longrightarrow e^{ta}\psi \in C, \qquad  t \geq 0,
$$
is studied. We will refer to this property as to the \emph{invariance} of $C$ under $\etasg$.
El-Maati Ouhabaz proved that under suitable assumptions on the form $(a,V)$ 
the invariance of $C$ is equivalent to $P_CV \subset V$ and $a( \psi, \psi-P_C\psi) \leq 0$ for all $\psi \in V$.
This result will be the key tool of our investigations.

\medskip
We start our work observing that a sesquilinear form $(a,V)$ on the direct sum $$V:= \bigoplus_{i \in I} V_i$$ of Hilbert spaces $V_i$ admits a matrix representation
obtained by setting
$$
a_{ij}(\psi,\psi'):= a(\mathbb 1_j \otimes \psi , \mathbb 1_i \otimes \psi' ).
$$
Following this idea, it is possible to develop a matrix theory for such sesquilinear forms and to characterise or, at least, 
to give sufficient conditions for properties of $(a,V)$ in terms of properties of the individual mappings.

In the case that $V_i=V_j$ for all $i,j \in I$, this matrix theory can be used to investigate \emph{symmetry properties} of the semigroup $\etasg$
associated with the form $(a,V)$ on a certain Hilbert space $H$, i.e., to investigate invariance of subspaces of a certain form under the action of $\etasg$.

This is due to the fact that the invariance criterion proved by El-Maati Ouhabaz can be simplified in the case of closed subspaces.
Denote $P_\mathcal Y$ the orthogonal projection of $H$ onto a closed subspace $\mathcal Y$. 
Then, $e^{ta}\mathcal Y \subset \mathcal Y$ for all $t \geq 0$ if and only if
\begin{equation}\label{intro1}
P_\mathcal Y V \subset V
\end{equation}
and
\begin{equation}\label{intro2}
a(\psi,\psi')=0, \qquad \psi \in \mathcal Y, \psi' \in \mathcal Y^\perp.
\end{equation}

On the one hand, it turns out that the condition~\eqref{intro1} is always satisfied in the case of subspaces representing symmetries, 
if $V$ is a direct sum of Hilbert spaces.
Thus, invariance is equivalent to the orthogonality condition~\eqref{intro2}, and this can be characterised in terms of properties of the mappings $a_{ij}$.

On the other hand, if the form domain contains coupling terms, i.e., if $V \subset \bigoplus_{i \in I} V_i$ but $V$ is not an ideal, it is more difficult to develop a general matrix theory for forms. However, it is possible to systematically investigate special classes of spaces and couplings.

A possible example is the case of the Laplace operator on a network. Here, the form domain is defined by
$$
V:=\{\psi \in \bigoplus_{i \in I} H^1(0,1): \psi(0)\oplus \psi(1)  \in Y \}
$$
where $Y$ is a suitable subspace of $\ell^2(I) \bigoplus \ell^2(I)$, and the action of the form is defined by the one-dimensional Dirichlet integral.

Again, it is possible to systematically study symmetry properties of the associated semigroup $\etasg$. 
For the heat equation, in particular, the orthogonality condition~\eqref{intro2} is always satisfied.
Invariance is thus equivalent to the admissibility condition~\eqref{intro1}. 
A natural question is whether it is possible to characterise admissibility by graph theoretic properties.
This task is not trivial, but it is also possible to extend some results to more general cases than network equations.

\medskip
We start by observing that a sesquilinear form $(a,V)$ on a Hilbert space $V$ of the form $V:=\bigoplus_{i \in I} V_{i}$ 
can be thought of as a matrix of sesquilinear mappings $(a_{ij})_{i,j \in I}$.
Thus, some sufficient conditions for properties of the form $(a,V)$ can be proved applying linear algebraic (or elementary functional analytic) arguments
to suitable complex-valued matrices constructed from the form $(a,V)$. This is done in Section~\ref{sec:findimarg}.

However, these arguments fail in many cases even in the easy task of characterising continuity or coercivity of infinite form matrices.
In Section~\ref{sec:systemsI} we consider the case of an infinite strongly coupled system and we characterise for such systems both continuity and coercivity.

In Section~\ref{sec:gencoercivity} we move back to the investigation of coercivity in the most general case, and we show that it suffices to analyse the finite restrictions of the form $(a,V)$.
In the case of two-dimensional matrices it is even possible to characterise coercivity in terms of properties of the single mappings.

We start in Section~\ref{sec:operator} to address the issue of evolution equations, and we try to identify the domain of the operator associated with the form $(a,V)$.
We continue our investigations of evolution equations in Section~\ref{sec:evoleq}, 
where we characterise well-posedness and contractive properties of the semigroup $\etasg$ in terms of properties of the single mappings $a_{ij}$.
As an example, we prove in Theorem~\ref{positivity} that the semigroup $\etasg$ generated by a matrix of sesquilinear mappings is positive if and only if
the semigroups generated by the sesquilinear forms on the diagonal are positive and the mappings off-diagonal are negative mappings.

We completely devote Section~\ref{sec:symmetriesmatrix} to symmetry properties, i.e., to invariance of subspaces of a certain form. 
Theorem~\ref{symmetries} gives a complete characterisation of these invariance properties.

Section~\ref{sec:wave} and Section~\ref{sec:dynamic} contain two illustrative applications of the techniques developed in the first part of Chapter~\ref{hesse}. 
In Section~\ref{sec:wave} we consider a strongly damped wave equation, and in Section~\ref{sec:dynamic} a heat equation with dynamic boundary conditions.
In both cases, we first obtain well-posedness for a large choice of parameters by means of Proposition~\ref{perturbinter}, 
and then we investigate qualitative properties of the solutions.

Finally, in Section~\ref{sec:nondiag} we introduce the topic of non-diagonal domain. 
Network equations represent maybe the most important case of non-diagonal domains and this will be the object of Chapter~\ref{network}. 

The results of Chapter~\ref{hesse} are motivated and discussed in Section~\ref{sec:histremI}. 
In particular, we shortly discuss the advantages of a matrix formalism for sesquilinear forms compared to the matrix formalism for operators.
Further, we shortly describe the history of the different issues that we have addressed in Chapter~\ref{hesse}.

We begin Chapter~\ref{network} by explaining informally the basic ideas of defining sesquilinear forms on networks. 
Section~\ref{sec:infintro} can be seen as a quick course in integration by parts on networks.

Since we want to consider infinite networks in a $L^2$-setting, we show in Section~\ref{sec:graphoper} that all definitions that are usual for finite networks
also make sense for infinite networks. In particular, we investigate operator theoretic properties of the incidence matrices. 
Once we have done this, we can define in~\eqref{formdomain}-\eqref{networkform} the form domain and the action of the form $(a,V)$ that will be the object of the following sections.

In Section~\ref{sec:systemsII} we investigate symmetry properties for systems of diffusion equations on networks.
This section can be seen as a continuation of Section~\ref{sec:systemsI}.
We are able in Theorem~\ref{charadmiss} to completely characterise the admissibility of projections connected to symmetries.

Both Section~\ref{sec:infinitesymmetries} and Section~\ref{sec:treelayer} are an application of the results in the previous section.
In the first one, we address the issue of irreducibility for the heat equation on infinite networks, and we characterise those networks for which the heat equation is irreducible. In Section~\ref{sec:treelayer} we investigate symmetry properties of special types of networks. 

In Section~\ref{sec:general} we leave the framework of networks, and we show that systems of diffusion equations on $H^1(0,1)$ 
can be seen as a heat equation on a network, if the boundary conditions satisfy some properties.

The concept of symmetry is explained from the physical point of view in Section~\ref{sec:histremII}.
We distinguish between \emph{space} and \emph{gauge} symmetries and we show that the symmetries that we have investigated in our work are {gauge symmetries} in a genuine physical sense.

Finally, we discuss in the rest of Section~\ref{sec:histremII} the history of the different issues that we have discussed in the previous sections.

\nopagebreak
\vfill
\noindent
{\bf Acknowledgements:}\\\vspace{-0.3cm}

\nopagebreak
I am extremely grateful to my supervisor Prof.\,Wolfgang Arendt. 
I own him my mathematical knowledge, and I admire his exemplary way of thinking of mathematics and science.

\nopagebreak
\smallskip
My colleagues Dr.\,Delio Mugnolo and Robin Nittka were exciting partner for mathematical discussions and good friends: 
I am deeply indebted to them and their patience.

\nopagebreak
\smallskip
I warmly thank my diploma supervisors Prof.\,Rainer Nagel and Prof.\,Ulf Schlotterbeck for their interest to my new research and their constructive criticisms.

\nopagebreak
\smallskip
I also express my gratitude to

\nopagebreak
\smallskip
Prof.\,Joachim von~Below, PD\,Jens Bolte, Prof.\,Klaus Engel,  Jun.\,Prof.\,Balint Farkas, PD\, Markus Haase, Markus Kunze, 
Heiner Markert, Prof.\,G\"unther Palm, Dr.\,Olaf Post, Prof.\,Frank Steiner

\nopagebreak
\smallskip
for their interest in my research and for the many interesting discussions.

\nopagebreak
\smallskip
The Institute for Applied Analysis and the Graduate School ``Mathematical Analysis of Evolution, Information and Complexity" in Ulm are stimulating environments and excellent places for a researcher.

\smallskip
I dedicate this work to Ulla K\"onig. I am grateful for her strong encouragement and support in the last months of my work.

%% file: sec_findimarg.tex
\pagenumbering{arabic}
\section{Finite dimensional arguments}\label{sec:findimarg}
The goal of this introductory section is twofold. 
First, we introduce the matrix representation of a form defined on the direct sum of Hilbert spaces $V:=\bigoplus_{i\in I} V_i$.
Then, we show that it is possible to deduce properties of such a form by properties of its matrix representation.
We recall the definition of the direct sum of Hilbert spaces.
\begin{defin}\label{formdomdef}
Let $I=\{1,\ldots,m\}$ or $I=\mathbb N$ be an index set and $V_i$ be Hilbert spaces for all $i \in I$. We consider on the space
\begin{equation}\label{formdomaindiag}
 V:=\bigoplus_{i\in I}V_i=\{(\psi_i)_{i \in I}: \psi_i \in V_i \mbox{ and }\sum_{i \in I} \|\psi_i\|^2_{V_i}<\infty  \},
\end{equation}
the scalar product $(\cdot \mid \cdot ) : V \times V \to \mathbb C$ by
\begin{equation}\label{scalar}
(\psi \mid \psi'):=\sum_{i \in I} (\psi_i \mid \psi'_i)_{V_i}.
\end{equation}
The inner product space $\big(V, (\cdot \mid \cdot)\big)$ is complete and thus it is Hilbert space.
\end{defin}
In the case of a finite set $I$, the direct sum coincides with the set theoretical Cartesian product. 
Therefore, we call the Hilbert space $V$ a \emph{product Hilbert space.}

Before we formally introduce the concept of a sesquilinear form on $V$ induced by a family of mappings, 
we recall how it is possible to define a matrix associated with a sesquilinear form on $\mathbb R^m \times \mathbb R^m$.
Later in this section, we show that the same idea also applies for sesquilinear forms on infinite-dimensional Hilbert spaces.

Consider a linear mapping $L$ on the Hilbert space $\mathbb C^m$, on which we have fixed an orthonormal basis $B=\{b_1,\ldots,b_m\}$. 
Further, we consider the matrix representation $A=(a_{ij})_{i,j=1,\ldots,m}$ of this linear mapping $L$. Since
$(Lb_j \mid b_i)= (a_{\cdot j} \mid b_i)=a_{ij},$
it is possible to determine the entries of the matrix representation (with respect to $B$) of $L$  in terms in terms of the scalar product.

Moreover, the sesquilinear form $a:\mathbb C^m\times \mathbb C^m \to \mathbb C$  defined by 
$$ a(u,v):=(Lu \mid v)$$
can be expressed in terms of the one-dimensional sesquilinear mappings $$a_{ij}:\langle b_j \rangle \times \langle b_i \rangle \to \mathbb C, \qquad a_{ij}(\lambda_j b_j, \mu_i b_i):=\lambda_j\bar{\mu_i}a_{ij}.$$ To see this, fix two vectors $u=\sum_{i=1}^m \lambda_i b_i$ and $v=\sum_{i=j}^m \mu_j b_j$. 
Computing now
\begin{eqnarray*}
a(u,v)& = & a\left(\sum_{i=1}^m \lambda_i b_i,\sum_{j=1}^m \mu_j b_j\right) \\
&=&\left(L \sum_{i=1}^m \lambda_i b_i \mid \sum_{j=1}^m \mu_j b_j\right)\\
&=&\sum_{i,j=1}^m a_{ij}( \lambda_j b_j , \mu_i b_i),
\end{eqnarray*}
shows how it is possible to define a sesquilinear form on  the space $H$ by means of a family of sesquilinear forms $a_{ij}$. 
This formal computation also holds if the (possibly countably many) vectors $b_j$ are elements of an infinite dimensional Hilbert space.

We prove two simple results. First, we define a sesquilinear form on the space $V$ by means of a family of sesquilinear mapping. 
Second, we show that each form on such a product space always induces a family of sesquilinear mappings $a_{ij}:V_j \times V_i \to \mathbb C$. 
Observe that the case of operator matrices is much more involved, see Section \ref{sec:histremI} for historical considerations on this topic.

We start by discussing the first implication, i.e., we show under which conditions a family of continuous sesquilinear mappings induces a continuous sesquilinear form on $V$. 
Recall that the mapping $a_{ij}:V_j\times V_i \to \mathbb C$ is said to be \emph{continuous} if for some $M_{ij } \geq 0$ the estimate 
$$
|a_{ij}(\psi,\psi')| \leq M_{ij}\|\psi\|_{V_j}\|\psi'\|_{V_i}, \qquad  \psi \in {V_j}, \psi'\in {V_i}
$$ 
holds. Fix now a family of continuous mappings $a_{ij}: V_j \times V_i \to \mathbb C$ and define an infinite matrix $M:=(M_{ij})_{i,j\in I}$.
We associate to the matrix $M$ an operator (which we also denote by $M$) defining it on its maximal domain
$$
D(M_{max}):=\{ x \in \ell^2(I): (M_{ij}x_j)_{j \in I}\in \ell^1(I) \mbox{ for all $i$ and} \sum_{i \in I} \big| \sum_{j \in I} M_{ij}x_j|^2 < \infty\}.
$$
Since the mappings $a_{ij}$ are continuous, we formally compute
$$ 
|\sum_{i,j \in I} a_{ij}(\psi_j,\psi_i)| \leq \sum_{i,j \in I} M_{ij}\|\psi_j\|_{V_j}\|\psi_i\|_{V_i} = M(\|\psi_i\|_{V_i})_{i\in I},
$$
and this makes sense if and only if $(\|\psi_i\|_{V_i})_{i\in I} \in D(M_{max})$.
In fact, we have proved the following result.
\begin{lemma}\label{boundedform}
Consider a family of continuous sesquilinear mappings 
$a_{ij}: V_j \times V_i \to \mathbb C$, and assume $a_{ij}$ to be continuous with continuity constants $M_{ij}$ for all $i,j \in I$, i.e.,
$$
|a_{ij}(\psi,\psi')| \leq M_{ij}\|\psi\|_{V_j}\|\psi'\|_{V_i}, \qquad  \psi \in {V_j}, \psi'\in {V_i}.
$$ 
If the linear operator $M$ on $\ell^2(I)$ defined by
$M:=(M_{ij})_{i,j \in I}$ is bounded,
then the sesquilinear form $ a: V \times V \to \mathbb C$ defined by
\begin{equation}\label{form}
a(\psi,\psi'):=\sum_{i,j\in I} a_{ij}(\psi_j,\psi'_i), \qquad \psi,\psi' \in V
\end{equation}
is continuous with constant $\|M\|_{\ell^2(I)}$.
\end{lemma}

Conversely, if $a$ is a continuous sesquilinear form on $V$, then it is possible to give a matrix representation of the form $a$.
\begin{lemma}\label{matrixform}
Let $a:V \times V \to \mathbb C$ be a continuous sesquilinear form on a product Hilbert space $V$. Then there exist uniquely determined continuous sesquilinear mappings $a_{ij}:V_j\times V_i \to \mathbb C$ such that $a(\psi)=\sum_{i,j \in I} a_{ij}(\psi_j,\psi_i)$.
\end{lemma}
\begin{proof}
We start defining the projection $\pi_i \in \Li (V,V_i)$ by
$$
\pi_i(\psi):=\psi_i.
$$
Observe now that $\pi_i$ does possess a (not uniquely determined) isometric right inverse $\pir \in \Li(V_i,V)$ defined by
$$
(\pir(\psi))_j:=
\begin{cases}
\psi & j=i,\\
0& \mbox{otherwise,}
\end{cases}
\qquad \psi \in V_i.
$$
One sees that the mapping $\pir$ satisfies the identity
\begin{equation}\label{piridentity}
\sum_{i\in I}\pir (\pi_i(\psi))=\psi, \qquad \psi \in V.
\end{equation}
Define now mappings $a_{ij}:V_j \times V_i \to \mathbb C$ by
\begin{equation}\label{factorforms}
a_{ij}(\psi, \psi'):=a(\pjr (\psi), \pir (\psi')),
\end{equation}
whose sesquilinearity is clear because of the analogous property of $(a,V)$. 

We prove that the mappings $a_{ij}$ are continuous. Denote $M$ the continuity constant of $(a,V)$ and compute for $\psi \in V_j, \psi' \in V_i$
\begin{eqnarray*}
|a_{ij}(\psi, \psi')|&=&|a(\pjr (\psi), \pir (\psi'))|\\
 &\leq& M \|\pjr (\psi)\|_V \|\pir (\psi')\|_V\\
&=&M \|\psi\|_{V_j} \|\psi'\|_{V_i}.
\end{eqnarray*}
We now show that the series $\sum_{i,j\in I}a_{ij}(\psi_j,\psi_i)=a(\psi)$ for all $\psi \in V$. 
By the definition of $a_{ij}$ and the sesquilinearity of $(a,V)$ we compute
\begin{eqnarray*}
\sum_{i,j\in I}a_{ij}(\psi_j,\psi_i)&=&\sum_{i,j\in I} a(\pjr (\pi_j(\psi)), \pir (\pi_i(\psi')))\\
&=&a(\sum_{j\in I}\pjr (\pi_j(\psi)),\sum_{i\in I}\pir (\pi_i(\psi)))\\
&=& a(\psi).
\end{eqnarray*}
Here all series converge since $(a,V)$ is continuous.
To see that the mappings $a_{ij}$ are unique, assume that there exists a second family 
$b_{ij}: V_j \times V_i \to \mathbb C$ satisfying $a(\psi,\psi')=\sum_{i,j \in I} b_{ij}(\psi_j,\psi'_i)$.
Computing
\begin{eqnarray*}
a(\pjr(\psi), \pir(\psi')) & = & \sum_{i,j\in I} b_{ij}(\pjr(\psi)_j, \pir{\psi'}_i)\\
&=& b_{ij}(\psi,\psi')
\end{eqnarray*}
shows that $a_{ij}=b_{ij}$.
\end{proof}

At a first glance, both above results seem quite satisfying. 
In fact, Lemma~\ref{matrixform} shows that the analogy between matrices of complex numbers and matrices of sesquilinear mappings also holds in general. 
Therefore,  we sometimes write with an abuse of notation $a=:(a_{ij})_{i,j\in I}$ for the sesquilinear form defined in \eqref{form}.  
However, the continuity bound given in Lemma \ref{boundedform} is far away from being optimal, even in the case of optimal bounds $M_{ij}$ and finite dimensional Hilbert spaces $V_1=V_2=\mathbb C$.  
\begin{exa}\label{regularoperators}
Consider the sesquilinear mappings $a_{ij}: \mathbb C \times \mathbb C \to \mathbb C$ defined by 
$$
a_{11}(z_1,z_2):=a_{12}(z_1,z_2):=a_{22}(z_1,z_2):= z_1 \overline{z_2}, 
$$
and
$$
a_{21}(z_1,z_2):=-z_1\overline{z_2}.
$$
Then the bound given in Lemma \ref{boundedform} on the continuity constant is the norm of the matrix 
$$
M=\begin{pmatrix}
1&1\\1&1
\end{pmatrix},
$$
whereas the optimal continuity norm of $(a,V)$ is given by the norm of the matrix
$$
M'=\begin{pmatrix}
1&1\\-1&1
\end{pmatrix}.
$$
It turns out that $\sqrt{2}=\|M'\|<\|M\|=2.$

In Section~\ref{sec:histremI} we shortly discuss the theoretical background of this phenomenon.
\end{exa}
Despite such unsatisfactory phenomena, it is possible to use both arguments given in the Lemmas \ref{boundedform} and \ref{matrixform} and obtain a characterisation of continuous sesquilinear forms in the case of a finite index set $I$.
\begin{cor}\label{finiteform}
Consider a family of sesquilinear mappings $\{a_{ij}: i,j \in I\}$ on a finite index set $I$. Then the matrix of sesquilinear forms $a=(a_{ij})_{i,j \in I}$ is continuous if and only if all mappings $a_{ij}$ are continuous.
\end{cor}
We now want to give some finite dimensional arguments for the {coercivity} of a form. 
Recall that a sesquilinear form $a:V \times V \to \mathbb C$ is said to be \emph{coercive}, 
if there exists an $\alpha>0$ such that $\re a(\psi) \geq \alpha \|\psi\|_V^2$ for all $\psi \in V.$ 
Analogously, we say that a bounded linear operator $L$ on $\ell^2(I)$ is coercive if $\re (Lv \mid v)\geq \alpha \|v\|^2_{\ell^2(I)}$, 
and this is equivalent to the positive definiteness of the symmetric operator $L+L^*$. 
Consider the matrix representation $A=(a_{ij})_{i,j\in I}$ of the linear operator $L$ with respect to an orthonormal basis $(b_i)_{i\in I}$. 
For the vectors of the orthonormal basis $B$ compute
$$
\re a_{ii}=\re (Lb_i \mid b_i) \geq \alpha \|b_i\|_{\ell^2(I)}=\alpha.
$$
As a consequence, the elements on the main diagonal are bounded from below by the coercivity constant $\alpha$.
The same argument also works in the for a form matrix $(a,V)$.

\begin{lemma}\label{coercivformnec}
Consider a sesquilinear form $a:V\times V \to \mathbb C$. If $a=(a_{ij})_{i,j\in I}$ is coercive with constant $\alpha$, then the sesquilinear forms $a_{ii}: V_i \times V_i \to \mathbb C$ are coercive with constant $\alpha$.
\end{lemma}
\begin{proof}
Let $\psi \in V_i$ and compute
$$
a_{ii}(\psi_i)=a(\pir(\psi_i)) \geq \alpha \|\pir(\psi_i)\|^2_V=\alpha \|\psi_i\|^2_{V_i}.
$$
This is the claim.
\end{proof}
Using the same argument as in Lemma \ref{boundedform}, we can also give a sufficient condition for the coercivity of a form $a=(a_{ij})_{i,j\in I}$.
\begin{prop}\label{coercivform}
Consider a sesquilinear form $a=(a_{ij})_{i,\in I}$ of continuous mappings $a_{ij}$ and denote the continuity constants by $M_{ij}$. 
Assume the diagonal forms $a_{ii}$ to be coercive with constants $\alpha_{ii}$. Define the matrix $A$ by
$$
A_{ij}:=\begin{cases} \alpha_{ii} & i=j\\ -M_{ij}&i\neq j. \end{cases}
$$
If the matrix $A$ is a bounded operator, and it is coercive  with constant $\alpha$, 
then the form $(a,V)$ is coercive with the same constant.
\end{prop}
\begin{proof}
The result is a consequence of the following computation
\begin{eqnarray*}
\re a(\psi)&=&\re \sum_{i\in I} a_{ii}(\psi_i) + \re \sum_{i\neq j \in I} a_{ij}(\psi_j,\psi_i)\\
&\geq& \sum_{i\in I} \alpha_{ii} \|\psi_i\|^2_{V_i} - \sum_{i\neq j \in I} M_{ij} \|\psi_j\|_{V_j}\|\psi_i\|_{V_i}\\
&=& A(\|\psi_i\|_{V_i})_{i\in I} \geq \alpha \|\psi\|^2_V.
\end{eqnarray*}
\end{proof}
As in the above result about boundedness, the coercivity estimate is not optimal. 
In the next section we will discuss a standard example for which it is possible to characterise both continuity and coercivity.

%% file: sec_systemsI.tex
\section{Systems of diffusion equations}\label{sec:systemsI}

We have seen that it is not a trivial task to characterise continuity and coercivity of form matrices, even in the case of a finite index set $I$. 
Assuming that all Hilbert spaces $V_i$ are equal and that the mappings $a_{ij}$ have the same structure, i.e., 
considering the case of a strongly coupled system allows us to give better results. 

To be more precise, we assume $V_i=H^1_0(\Omega)$ or $V_i=H^1(\Omega)$ for a bounded domain $\Omega \subset \mathbb R^d$. 
Moreover, $C$ will denote some bounded operator--valued function $C: \Omega \to \mathcal L(\ell^2(I))$. 
In this setting, we define the form $a: V \times V \to \mathbb C$ by
\begin{equation}\label{systemsform}
a(\psi,\psi'):= \int_{\Omega} (C(x) \nabla \psi(x) \mid \nabla \psi'(x) ) dx.
\end{equation}
We want to discuss two different topics in this section.
First, we address the issue of continuity for systems defined on $H^1(\Omega)$.
Second, we investigate coercivity properties for systems defined on $H^1_0(\Omega)$.
We want to focus our attention on infinite systems and therefore, we restrict ourselves to the case $I = \mathbb N$.

In order to fix the ideas, let us summarise the standing assumptions for this section.
\begin{ass}
In this section we always assume the following.
\begin{itemize}
\item The set $\Omega$ is a domain in $\mathbb R^d$.
\item The function $C: \Omega \to \mathcal L(\ell^2(\mathbb N)), x\mapsto (c_{ij}(x))_{i,j \in \mathbb N}$ is componentwise measurable. 
\item If there is no danger of confusion, we denote by $C$ also the range of the function $C$, i.e., the operator family $C:=\{C(x): x\in \Omega\}$.
\item The Hilbert space $H$ is defined as in \eqref{formdomaindiag} by
\begin{equation}\label{hilbertspace}
H:=\bigoplus_{i\in \mathbb N} L^2(\Omega).
\end{equation}
\end{itemize}
\end{ass}
Since the function $C$ is an operator for all $x$ and since $\psi(x) \in \ell^2(\mathbb N)$ for almost every $x$, 
it is possible to define for almost every $x\in \Omega$ and for all $n \in \mathbb N$
\begin{equation}\label{multoper}
\big(M_C \psi\big)_n (x) := \sum_{m\in\mathbb N} c_{nm}(x) \psi_m(x).
\end{equation}
The function $C$ is componentwise measurable. 
Thus, also $M_C \psi : x \mapsto \big(\sum_{m\in\mathbb N} c_{nm}(x) \psi_m(x))_{n\in \mathbb N}$ is componentwise measurable.
As a consequence, we can define $M_C$ on its maximal domain
$$
D(M_{C_{max}}):=\{ \psi \in \ell^2(E): M_C \psi \in \ell^2(E)\}.
$$
In the following, we will need estimates for the gradient of $\psi$. 
We thus introduce a multiplication operator $M_C^d$ on $H^d$. This is defined on 
$$
D(M^d_{C_{max}}):=\bigoplus_{i=1,\ldots,d}D(M_{C_{max}})
$$ 
by
\begin{equation}\label{multoperd}
\big(M^d_C \psi\big) := (M_C \psi^1, \ldots, M_C \psi^d), \quad \psi=(\psi^i)_{i=1,\ldots,d} \in D(M^d_{C_{max}}).
\end{equation}
In the next Lemma, we show that $H$ admits different representations.
\begin{lemma}\label{Hrepresent}
Assume that $X_i$ are $\sigma$-finite measure spaces. Then
\begin{equation}
\bigoplus_{i\in \mathbb N} L^2(X_i) \simeq L^2\big(\bigoplus_{i\in \mathbb N} X_i \big).
\end{equation}
If $X_i=X$ for all $i \in \mathbb N$ both spaces are also isometrically isomorphic to  $L^2(X,\ell^2(\mathbb N)).$
In particular,
\begin{equation}\label{sysnorm}
\|\psi\|^2_{L^2(\bigoplus_{i \in \mathbb N}  X)}= \int_{X} \|\psi(x)\|^2_{\ell^2} dx.
\end{equation}
\end{lemma}
\begin{proof}
For the sake of the completeness, we recall the construction of the measure space $\bigoplus_{i\in \mathbb N} X_i$.
The direct sum of sets is defined by
$$
X:=\bigoplus_{i\in \mathbb N}X_i:=\{(i,x): i \in \mathbb N, x \in X_i\}.
$$
We define a $\sigma$-algebra $\Sigma$ on $X$. Observe that each subset $A\subset X$ has the form 
$$
A=\bigoplus_{i \in  \mathbb N} A_i, \quad A_i:=\{x \in X_i: \exists y \in A, y=(i,x) \}.
$$
We define a $\sigma$-algebra $\Sigma$ by requiring that $A \in \Sigma$ if and only $A_i$ is Lebesgue measurable for each $i \in I.$
It is a simple exercise to check that $\Sigma$ is indeed a $\sigma$-algebra on $X$.
For all $A \in \Sigma$ define
$$
\lambda(A):=\sum_{i \in \mathbb N} \lambda_i(A_i).
$$
Since all terms in the above sum are positive, the series converges if and only if it converges absolutely. 
Thus, convergence is independent of the order of summation and the expression $\sum_{i \in I}\lambda_i(A_i)$ makes sense. 
Then, $\lambda(A)=0$ if and only if $\lambda_i(A_i)=0$ for all $i \in \mathbb N$. 
Since the Cartesian product of countable sets is countable, we see that the triple $(X,\Sigma,\lambda)$ is a $\sigma$-finite measure space. 
Define a mapping $\phi: L^2(X) \to H$ by 
$$
(\phi(\psi))_i(x):=\psi_i(x) \qquad x \in X_i.
$$
The mapping $\phi$ is bijective. We prove that it is isometric. 
Fix to this aim $\psi \in L^2(X)$ and observe that
the sets $X_i$ are a measurable partition of $X$. Thus, by the definition of the norm in $L^2(X)$ and
the monotone convergence theorem, we obtain
\begin{eqnarray*}
\|\psi\|^2_{L^2(X)} & =& \int_{X} |\psi(x)|^2 dx =  \int_X \sum_{i \in \mathbb N}| \mathbb 1_{X_i} \psi(x)|^2 dx\\
& = & \sum_{i\in I} \int_{X} | \mathbb 1_{X_i} \psi(x)|^2 dx =  \sum_{i \in \mathbb N}\int_{X_i} |\psi_i(x)|^2 dx.
\end{eqnarray*}
The latter is the definition of the norm in $H$.

\smallskip
Assume now that $X_i = \Omega$ for all $i \in \mathbb N$ and define a mapping $\phi: H \to L^2(\Omega, \ell^2(\mathbb N))$ by
$$
(\phi(\psi))_i(x):=\psi_i(x) \qquad x \in X_i.
$$
Since $\psi_i \in L^2(\Omega)$ for all $i\in \mathbb N$, the function $f:x \mapsto \sum_{i \in \mathbb N}|\psi_i(x)|^2$ is defined almost everywhere and is a measurable function as series of positive measurable functions. Thus, the mapping $\phi$ is well-defined and it is clearly bijective.  
The formula~\eqref{sysnorm} follows again from the monotone convergence theorem, and so the claim is proven.
\end{proof}
We are now in the position of attacking the main problem.
As first result, we characterise the continuity of the form $(a,V)$ in the case $V_i= H^1(\Omega)$ for all $i \in \mathbb N$. 
As usual, we denote $V$ the Hilbert space defined in \eqref{formdomaindiag}. 
In order to achieve the goal of characterising continuity, we need to obtain sharp estimates for the $L^2$ norm of the gradient of functions in $V$,
and for this optimal bounds for the operator $M^d_C$ are needed. 
We start showing that the norm boundedness of the range of $C$ is equivalent to the boundedness of the operator $M_C$.

\begin{lemma}\label{multoperbound}
Recall the definition~\eqref{hilbertspace} of the Hilbert space $H$. The following assertions are equivalent for an arbitrary $M \in \mathbb R$.
\begin{enumerate}[a)]
\item The operator $M_C$ is bounded by $M$, i.e., $\|M_C\|_{\mathcal L(H)}\leq M$.
\item The operator family $C$ is uniformly essentially bounded, i.e., $\|C(x)\|_{\mathcal L(\ell^2)} \leq M$ {for almost all } $x \in \Omega.$
\end{enumerate}
\end{lemma}
\begin{proof}
Assume the family $C$ to be uniformly essentially bounded. Use the definition of the norm in $H$ to see
$$
\|M_C \psi\|^2_H =  \sum_{n\in \mathbb N} \int_\Omega \left|(C(x)\psi(x))_n\right|^2 dx.
$$
Since the series converges absolutely, this yields
$$
\|M_C \psi\|^2_H \leq  M \int_{\Omega}  \|\psi(x)\|^2_{\ell^2} dx = M\|\psi\|^2_H,
$$
and this implies a).

We now prove the converse implication. Assume $M_C$ to be bounded by $M$ and fix an arbitrary normed sequence $y:=(y_i)_{i \in \mathbb N} \in \ell^2$. Define now functions in $L^2(\Omega)$ by
$$
\psi_i:=\chi_{\omega} y_i,
$$
where $\lambda(\omega) < \infty.$

Since $y \in \ell^2$ the vector $\psi:=(\psi_i)_{i\in \mathbb N}$ is an element of $H$, whose norm is given by
$$
\|\psi\|_H= \sum_{i\in N} \|\psi_i\|^2_{L^2}=\lambda(\omega)\sum_{i\in N} |y_i|^2 = \lambda(\omega).
$$
As a consequence, the boundedness of $M_C$ yields
\begin{equation}\label{hilf1}
\|M_C \psi\| \leq M \lambda(\omega).
\end{equation}
Observe that
$$
M_C\psi(x)=\begin{cases} C(x)y & x\in \omega, \\ 0 & \mbox{otherwise}. \end{cases}
$$
Now are we in the position of completing the proof. Fix an arbitrary $x_0 \in \Omega$ and set $\omega_j=B_{\frac{1}{j}}(x_0)$. 
Let us denote $\beta(d)$ the volume of the unit ball in $\mathbb R^d$. 
As a consequence $\lambda(\omega_j)= \frac{\beta(d)}{j^d}$. By the Lebesgue-Besicovitch Theorem (see e.g. \cite[Theo.~1.7.1]{EvaGar92}) we obtain for almost all $x_0 \in \Omega$
$$
\lim_{j\to\infty} \frac{j^d}{\beta(d)} \int_{\omega_j} (M_C \psi)_i(x) dx 
=\lim_{j\to\infty} \frac{j^d}{\beta(d)} \int_{\omega_j} (C(x)y)_i dx
=(C(x_0)y)_i.
$$
for all $i\in \mathbb N.$ 
Since~\eqref{hilf1} holds, by the monotone convergence theorem we obtain
$$
\|C(x_0)y\| = \lim_{j\to\infty} \frac{j^d}{\beta(d)} \|M_C \psi\|.
$$
Now, passing to the limit yields
$$
\|C(x_0)y\| = \lim_{j\to\infty} \frac{j^d}{\beta(d)} \|M_C \psi\| \leq \lim_{j\to\infty} \frac{j^d}{\beta(d)} \lambda(\omega_j)M =M.
$$
what we had to prove.
\end{proof}

\begin{rem}\label{vecmultoper}
The operator $M_C^d$ defined in \eqref{multoperd} is a bounded operator on $H^d$ if and only if $M_C$ is a bounded operator on $H$.
\end{rem}
We are now going to define the sesquilinear form whose continuity we want to characterise. We introduce to this aim some notations. 
For a function $\psi \in V$ define the gradient of $\psi$ by
$$\nabla\psi:=(\nabla \psi_n)_{n\in \mathbb N}, \qquad \psi \in V.$$
We observe that, by virtue of the definition of $V$, for all $\psi \in V$ 
$$(\| \nabla\psi_n \|_{L^2(\Omega)})_{n\in \mathbb N} \in \ell^2(\mathbb N),$$ and so $\nabla \psi \in \mathcal L(V,\ H^d)$. 
Recall that for functions $\psi,\psi'$ the scalar product on $H^d$ is defined by
$$
(\psi \mid \psi')_{H^d}=\sum_{k=1}^d (\psi^k \mid \psi'^k)_H= \sum_{k=1}^d \sum_{i=1}^\infty(\psi^k_i \mid \psi'^k_i)_{L^2(\Omega)}
$$ 
On finite sequences of $V$ it is possible to define a sesquilinear form $a_C$ by
\begin{equation}\label{SEI}
a_C(\psi,\psi'):=(M^d_C \nabla \psi \mid \nabla \psi')_{H^d}.
\end{equation}
In the next Proposition we characterise the continuity of the form $a_C$ on the whole space $V$.
\begin{prop}\label{contsystems}
Assume the set $\Omega$ to have finite measure. 
Then the form $a_C$ defined in \eqref{SEI} is continuous as a mapping $a_C:V \times V \to \mathbb C$
if and only if the operator family $C$ is essentially uniformly bounded.
\end{prop}

\begin{proof}
Assume the form $a_C$ to be continuous with constant $M$. In particular, for each $k=1,\ldots,d$ the form $a_{C,k} : V \times V \to \mathbb C$ defined by
\begin{equation}\label{SEIpartial}
a_{C,k}(\psi,\psi'):=\left(M_C \frac{\partial}{\partial x_k}\psi \mid \frac{\partial}{\partial x_k} \psi'\right)_{H}
\end{equation}
is continuous. Thus,
$$
\left|
\left(M_C \frac{\partial}{\partial x_k}\psi \mid \frac{\partial}{\partial x_k} \psi'\right)_{H}
\right| \leq M_1, \quad \mbox{for all } \psi,\psi' \in B_1^V(0).
$$
The operator $\frac{\partial }{\partial x_k} \in \mathcal L(H^1(\Omega),L^2(\Omega))$ is a bounded and surjective linear operator.
So it is an open mapping and this implies
$$
\left|
\left(M_C \psi \mid  \psi'\right)_{H}
\right| \leq M_2, \quad \mbox{for all } \psi,\psi' \in B_1^H(0).
$$
The formula
$$
\|M_C\|=\sup_{\psi,\psi'\in B_1^H(0)} (M_C \psi \mid \psi'),
$$
together with the Lemma \ref{multoperbound} shows the claim.

Conversely, assume the operator family $C$ to be bounded. 
By Lemma \ref{multoperbound} and Remark \ref{vecmultoper} the operator $M_C^d$ is bounded. Thus, the estimate
$$
|a_C(\psi,\psi')| = |(M^d_C J\psi \mid J \psi')| \leq \|M^d_C\| |(J\psi \mid J \psi')| \leq \|M^d_C\|\|\psi\|_V\|\psi'\|_V
$$
yields the continuity of $a_C$, thus completing the proof.
\end{proof}

\medskip
We now turn our attention to coercivity properties of the form defined in~\eqref{systemsform} 
in the case $V_i=H^1_0(\Omega)$ for all $i\in \mathbb N$ and $V$ as in \eqref{formdomaindiag}. 
Analogously as in the previous part of the section we want to prove the equivalence between the coercivity of $(a,V)$ and the that of essential uniform coercivity of $C$. In order to prove this result the equivalence of the norm on $V$ and the scalar product defined by $[u \mid v]:= (\nabla u \mid \nabla v)_{L^2}$ is needed. 
Domains having this property are said to be of \emph{Poincar\'e type}.
In particular, bounded domains and domains lying in a strip are of {Poincar\'e type}, see \cite[Prop.~3.4.1]{Are06}.
\begin{thm}\label{coercsys}
Assume the operator--valued function $C$ to be of class $C(\overline{\Omega})$, uniformly bounded and accretive. 
The form $(a,V)$ is coercive if and only if the operator family $C$ is uniformly coercive.
\end{thm}
As before, we split the proof into two steps. First, we prove the equivalence between the uniform coercivity of $C$ and the coercivity of $M_C$. Second, we transfer the information about $M_C$ to the form $(a,V)$.
\begin{lemma}\label{multopercoerc}
If $\Omega$ is of Poincar\'e type and $C$ is essentially uniformly bounded, the following assertions are equivalent.
\begin{enumerate}[a)]
\item The family $C$ is essentially uniformly coercive with constant $\alpha$, i.e.,
$$
\re (C(x) y \mid y)_{\ell^2} \geq \alpha \|y\|^2_{\ell^2} \quad \mbox{for almost all } x \in \Omega;
$$
\item The operator $M_C$ is coercive with constant $\alpha$, i.e.,
$$
\re (M_C \psi \mid \psi)_H \geq \alpha \|\psi\|^2_H, \quad \mbox{for all } \psi \in H.
$$
\end{enumerate}
\end{lemma}
\begin{proof}
Assume the family $C$ to be essentially uniformly coercive. We thus compute
\begin{eqnarray*}
\re (M_C \psi \mid \psi)_H &=& \re \sum_{n\in \mathbb N} \int_\Omega \sum_{m\in \mathbb N} c_{nm}(x)\psi_m(x) \overline{\psi_n(x)} dx\\
&=&\re \int_\Omega  \sum_{n,m\in \mathbb N}  c_{nm}(x)\psi_m(x) \overline{\psi_n(x)} dx\\
&=&\re \int_\Omega  (C(x)\psi(x) \mid \psi(x)) dx\\
&\geq& \int_\Omega \alpha \|\psi(x)\|^2_{\ell^2} dx = \alpha \|\psi\|^2_H ,
\end{eqnarray*}
where in the last step we used the formula~\eqref{sysnorm}.

The converse implication is a consequence of the Lebesgue-Besicovitch Theorem, and can be obtained exactly in the same way as in the analogous implication of Lemma \ref{multoperbound}.
\end{proof}
Exactly as before, $M_C$ is coercive if and only if $M_C^d$ is coercive. 
Observe that since $\Omega$ is of Poincar\'e type, the inequality $\| \nabla \psi_i \|^2_{L^2(\Omega)} \geq M_1 \|\psi_i\|^2_{H^1_0(\Omega)}$ 
holds for all $\psi_i \in H^1_0(\Omega)$.
Thus,
$$
\sum_{i \in \mathbb N} \|\nabla \psi_i\|^2_{L^2(\Omega)}= \|\nabla \psi\|^2_H \geq M_1 \|\psi\|^2_V.
$$
We can now move to the proof of the main result.
\begin{proof}[Proof of Theorem \ref{coercsys}]
First, we assume the family $C$ to be essentially uniformly coercive. Because of Lemma~\ref{multopercoerc} and of the Poincar\'e property
\begin{eqnarray*}
\re a(\psi) = (M^d_C \nabla \psi \mid \nabla \psi) \geq \alpha \|\nabla \psi\|^2_{H} \geq \alpha M_1 \|\psi\|^2_V.
\end{eqnarray*}
Thus, the form $a$ is coercive with constant $\alpha M_1$.

Conversely, assume that $C$ is not uniformly coercive. Thus, there exists a sequence $(x_k)_{k \in N} \subset \Omega$ and
a sequence $(y^k)_{k\in \mathbb N} \subset \ell^2(\mathbb N), \|y^k\|_{\ell^2(\mathbb N)}=1$ such that 
$$
\re (C(x_k) y^k \mid y^k) \leq \frac{1}{k}\|y^k\|^2_{\ell^2(\mathbb N)}.
$$
In particular, since $C$ is of class $C(\overline{\Omega})$, it is also uniformly continuous and there exists an uniform radius $r > 0 $ such that 
\begin{equation}\label{unifradius}
\re (C(x) y^k \mid y^k)\leq \frac{2}{k}\|y^k\|^2_{\ell^2(\mathbb N)}, \qquad x \in B_r(x_k) \cap {\Omega}.
\end{equation}
Consider now real valued functions $\psi^k \in C^\infty_c(\Omega)$ such that $\supp \psi^k \subset B_{r}(x_k)$.
Define $\widehat{\psi^k} \in V$ by
$$
(\widehat{\psi^k})_i:= y^k_i \psi^k,\qquad i,k \in \mathbb N.
$$
Observe that $\|\widehat{\psi^k}\|^2_V=\|{\psi^k}\|^2_{H^1(\Omega)}$ since $\|y^k\|_{\ell^2(\mathbb N)}=1$. Furthermore
$$
\frac{\partial \widehat{\psi^k}(x)}{\partial x_\ell} = \frac{\partial {\psi^k}(x)}{\partial x_\ell} y^k, \qquad \ell=1, \ldots,d,\, k \in \mathbb N,\, x \in \Omega,
$$
and so by~\eqref{unifradius}
\begin{eqnarray*}
(C(x) \frac{\partial \widehat{\psi^k}(x)}{\partial x_\ell} \mid \frac{\partial \widehat{\psi^k}(x)}{\partial x_\ell}) 
&=&(C(x) \frac{\partial {\psi^k}(x)}{\partial x_\ell} y^k \mid \frac{\partial {\psi^k}(x)}{\partial x_\ell} y^k)\\
&\leq& \frac{2}{k} |\frac{\partial {\psi^k}(x)}{\partial x_\ell}|^2 \|y^k\|^2_{\ell^2(\mathbb N)}\\
&=&\frac{2}{k} |\frac{\partial {\psi^k}(x)}{\partial x_\ell}|^2
\end{eqnarray*}
for all $ \ell=1,\ldots, d, k \in \mathbb N$ and all $x \in B_r(x_k) \cap {\Omega}.$
The estimate
\begin{eqnarray*}
\re a(\widehat{\psi^k})&=& \re \sum_{\ell=1}^d \int_{\Omega} (C(x) \frac{ \partial \widehat{\psi^k}(x)}{\partial x_\ell} 
\mid  \frac{\partial \widehat{\psi^k}(x)}{\partial x_\ell}) dx\\
&=& \re \sum_{\ell=1}^d\int_{B_r(x_k) \cap \Omega} (C(x) \frac{\partial \widehat{\psi^k}(x)}{\partial x_\ell} 
\mid \frac{\partial \widehat{\psi^k}(x)}{\partial x_\ell}) \\
& \leq & \frac{2}{k} \sum_{\ell=1}^d \int_{B_r(x_k) \cap \Omega} |\frac{\partial{\psi^k}(x)}{\partial x_\ell}|^2 dx \\
& \leq & \frac{2}{k} \|\psi^k\|^2_{H^1(\Omega)} =  \frac{2}{k} \|\widehat{\psi^k}\|^2_{V}
%&\leq & \re \frac{2}{k} \|\widehat{\psi^k}\|^2_V.
%&=& \sum_{\ell=1}^d \int_{B_r(x_k)\cap \Omega} (C(x) \widehat{\varphi^k}(x) \mid  \widehat{\varphi^k}(x)) dx \\
%&=& d \int_{B_r(x_k)\cap \Omega} (C(x)  \widehat{\varphi^k}(x) \mid  \widehat{\varphi^k}(x)) dx\\
%&\leq& \frac{2d}{k} \int_{B_r(x_k)\cap \Omega} \| \widehat{\varphi^k}(x)\|^2_{\ell^2(\mathbb N)} dx 
%=  \frac{2d}{k} \int_{B_r(x_k)\cap \Omega}  |\varphi^k(x)|^2 dx= \frac{2d}{k}.\\
\end{eqnarray*}
%By the Poincar\'e inequality, $\|\widehat{\psi^k}\|_V \geq C \|\widehat{\varphi^k}\|_H=C$, and so the form $a$ is not coercive.
shows that the form $(a,V)$ is not coercive.
\end{proof}

In the special case of systems we have characterised continuity and coercivity properties of the form $a$.
In the next section we characterise coercivity in another special case.

%% file: sec_gencoercivity.tex
\section{Coercivity}\label{sec:gencoercivity}
The characterisation of the coercivity of a form matrix is in general a complex task. 
We present two results in this section: the first gives a characterisation in terms of the finite restrictions of the form. In the second one, we discuss the coercivity in the case of a two-dimensional matrix.

We start by fixing the standing assumptions for this section.
\begin{ass}
During this section we always assume the following.
\begin{itemize}
\item The form $(a,V)$ is continuous.
\item Fix $I'\subset I$. 
We denote $a_{|I'}$ the the restriction of the form to $V_{|I'}:=\prod_{i\in I'} V_i$ and we call it the \emph{restriction of the form $(a,V)$ induced by $I'$}.
\end{itemize}
\end{ass}
\begin{thm}\label{coercrest}
The following assertions are equivalent.
\begin{enumerate}[a)]
\item The form $(a,V)$ is coercive with constant $\alpha$ 
\item Each finite $I' \subset I$  induces a coercive restriction with constant $\alpha$.
\end{enumerate}
\end{thm}
\begin{proof}
The necessity of the condition b) can be seen by defining the projection $\pi_{I'} \in \mathcal L(V,V_{|I'})$
$$
(\pi_{I'}(\psi))_i:=\psi_i, \qquad \mbox{for all } i\in I',
$$
and arguing as in the proof of Lemma \ref{coercivformnec}. 
The sufficiency of the condition b) in the case of a finite index set $I$ is tautological. So, we only have to prove the sufficiency in the case $I=\mathbb N$.

To this aim, let $\psi \in V, \|\psi\|_V =1$. 
We have to show $\re a(\psi) \geq \alpha$. 
First observe that each $\psi \in V$ can be uniquely decomposed as 
$$
\psi=\psi' + \psi'',\qquad  \psi' \in \bigoplus_{i \in I'} V_i, \psi'' \in \bigoplus_{i \in I''} V_i.
$$
Here $I''= \mathbb N\setminus I'$. By sesquilinearity we obtain
\begin{equation}\label{gencoerchilf1}
a(\psi)= a(\psi') + a (\psi'') + a(\psi',\psi'') + a(\psi'',\psi').
\end{equation}
Since $a(\psi')=a_{I'}(\pi_{I'}(\psi))$, we obtain
\begin{equation}\label{gencoerchilf2}
\re a(\psi') \geq \alpha \|\psi'\|^2_V.
\end{equation}
Further, by the continuity of $(a,V)$, we also see that
\begin{equation}\label{gencoerchilf3}
|a (\psi'') + a(\psi',\psi'') + a(\psi'',\psi')| \leq M \|\psi''\|^2_V + 2 M \|\psi'\|_V\|\psi''\|_V.
\end{equation}
Thus, using $\|\psi'\| \leq 1$ and estimating as in the proof of Proposition \ref{coercivform},
\begin{equation}\label{gencoerchilf4}
-(\re a (\psi'') + \re a(\psi',\psi'') + \re a(\psi'',\psi')) \geq -M\| \psi'' \|_V (\|\psi''\|_V + 2).
\end{equation}
Observe now that $\|\psi'\|^2_V + \|\psi''\|^2_V=1$, i.e.,
\begin{equation}\label{gencoerchilf5}
\|\psi'\|_V = \sqrt{1-\|\psi''\|^2_V}.
\end{equation} 
Combining equations \eqref{gencoerchilf1}-\eqref{gencoerchilf2}-\eqref{gencoerchilf3}-\eqref{gencoerchilf4}-\eqref{gencoerchilf5} we obtain the estimate
\begin{equation}\label{gencoerchilf6}
\re a(\psi) \geq \alpha \|\psi'\|^2_V - M(1-\|\psi'\|_V)-2 M \sqrt{1- \|\psi'\|^2}.
\end{equation}
We want to estimate the right hand-side of~\eqref{gencoerchilf6} from below.
Define the auxiliary function
$$
f(t):= \alpha t^2 - M (1-t) - 2M \sqrt{1-t^2}.
$$
Thus, we can reformulate \eqref{gencoerchilf6} as $\re a(\psi') \geq f(\|\psi'\|_V)$.
Since $f$ is continuous and $f(1) = \alpha$, we obtain
$$
\limsup_{n \to \infty} \re a(\psi')\geq \lim_{t \to 1} f(t) \geq \alpha,
$$
where $I':=\{1,\ldots,n\}$. Since $\lim_{n \to \infty} \psi' =\psi$ by definition and $(a,V)$ is a continuous form
$$
\re a(\psi)=\lim_{n\to \infty} \re a(\psi')=  \limsup_{n \to \infty} \re a(\psi') \geq \alpha,
$$
and the proof is complete.
\end{proof}
Next, we show that for two by two matrices it is possible to characterise the coercivity in terms of property of the single mappings $a_{ij}$. 
To this aim, for a finite index set $I$ and a vector $w:=(w_i)_{i \in I}, w_i > 0$ we define a weighted equivalent norm $\|\cdot\|_w$ on $V$ by
$$
\|\psi\|^2_w:= \sum_{i\in I} w_i\|\psi\|^2_{V_i}.
$$
\begin{prop}\label{propcoerc2x2}
Consider as index set $I:=\{1,2\}$ and coercive sesquilinear forms $a_{ii}: V_i \times V_i \to \mathbb C$ for $i \in I$, with optimal constants $\alpha_1,\alpha_2$, respectively. Then the following assertions are equivalent.
\begin{enumerate}
\item The form $(a,V)$ is coercive.
\item There exist  $\alpha >0$ such that for all $\psi=\begin{pmatrix} \psi_1 \\ \psi_2 \end{pmatrix} \in V$ and for $w=\begin{pmatrix}\alpha_1 \\ \alpha_2 \end{pmatrix}$
\begin{equation}\label{coerc2x2}
\re a_{12}(\psi_2,\psi_1)+ \re a_{21}(\psi_1,\psi_2) \geq  \alpha \|\psi\|^2_w - 2\sqrt{\alpha_1\alpha_2} \|\psi_1\|_{V_1}\|\psi_2\|_{V_2}.
\end{equation}
\end{enumerate}
\end{prop}
\begin{proof}
If there is no danger of confusion, we do not specify in which space we are computing the norm. For a real parameter $\lambda >0$
\begin{eqnarray*}
\re a(\psi)&=&\re a(\lambda \psi, \frac{1}{\lambda} \psi)\\
&=&\re \lambda^2 a_{11}(\psi_1,)+ \re \frac{1}{\lambda^2}a_{22}(\psi_2)\\
&&+ \re a_{12}(\psi_2,\psi_1)+ \re a_{21}(\psi_1,\psi_2)\\
&\geq& \alpha_1 \lambda^2\|\psi_1\|^2_{} + \alpha_2\frac{1}{\lambda^2} \|\psi_2\|^2_{}+ \re a_{12}(\psi_2,\psi_1)+ \re a_{21}(\psi_1,\psi_2)\\
&=:&f(\lambda)
\end{eqnarray*}
The auxiliary function $f : (0,\infty) \to \mathbb R$ has derivative
$$
f'(\lambda)=2\left( \alpha_1 \|\psi_1\|^2_{}  \lambda- \alpha_2\|\psi_2\|^2\lambda^{-3}\right),
$$
which vanishes only at
$$
\lambda_0=\sqrt[4]{\frac{ \alpha_2\|\psi_2\|^2_{} }{ \alpha_1 \|\psi_1\|^2_{}}}
=\sqrt{\frac{ \sqrt{\alpha_2}\|\psi_2\|_{} }{\sqrt{\alpha_1} \|\psi_1\|_{}}}.
$$
Inserting $\lambda_0$ into $f$ yields:
$$
f(\lambda_0)=2 \sqrt{\alpha_1\alpha_2} \|\psi_1\| \|\psi_2\|+\re a_{12}(\psi_2,\psi_1)+ \re a_{21}(\psi_1,\psi_2).
$$
The strict convexity of $f$ implies that $\lambda_0$ is the global minimum of the function $f$. 
Since the constants are optimal and $\psi_1,\psi_2$ can be chosen independently of each other, the coercivity is equivalent to
$$
f(\lambda) \geq \alpha \|\psi\|^2_V
$$
for some $\alpha > 0$. Since $\alpha_1, \alpha_2 > 0$, the norm in $V$ is of course equivalent to the weighted norm $\|\cdot\|_w$, where $w:=(\alpha_1,\alpha_2)$.
We deduce that the coercivity with respect to $\|\cdot \|$ is equivalent to the coercivity with respect to $\|\cdot\|_w$, i.e.,  to the existence of $\alpha >0$ such that
$$
2 \sqrt{\alpha_1\alpha_2} \|\psi_1\| \|\psi_2\|+ \re a_{12}(\psi_2,\psi_1)+ \re a_{21}(\psi_1,\psi_2) 
\geq \alpha \left( \alpha_1\|\psi_1\|^2+ \alpha_2 \|\psi_2\|^2\right),
$$
what we had to prove.
\end{proof}
We illustrate a possible use of this criterion.
\begin{exa}
On a domain $\Omega$ of class $C^\infty$ consider the operator matrix
$$
H:=\begin{pmatrix}
\Delta & M_V \\
M_V & \Delta
\end{pmatrix}
$$
on $D(H):=D(\Delta_{Dirichlet}) \times D(\Delta_{Dirichlet})$. 
Here $M_V$ denotes a bounded multiplication operator which is defined by
$$
M_V f (x):= V(x)f(x).
$$
The form matrix associated with the operator is defined on the Hilbert space $V:=H^1_0(\Omega) \times H^1_0(\Omega)$ by setting
$$
a_{11}(\psi,\psi'):=a_{22}(\psi,\psi'):= \int_\Omega \nabla \psi(x) \cdot \overline{\nabla \psi'(x)} dx, \qquad \psi,\psi' \in H^1_0(\Omega),
$$
$$
a_{12}(\psi,\psi'):=-\int_{\Omega} V(x) \psi(x) \overline{\psi'(x)} dx, \qquad \psi,\psi' \in H^1_0(\Omega),
$$
and
$$
a_{21}(\psi',\psi):= -\int_{\Omega} V(x) \psi'(x) \overline{\psi(x)} dx, \qquad \psi,\psi' \in H^1_0(\Omega).
$$
Assume now that $V(x) \in i \mathbb R$ for all $x \in \mathbb R$. Compute for $\psi=\begin{pmatrix}\psi_1\\ \psi_2\end{pmatrix} \in V$
\begin{eqnarray*}
&&\re (a_{12}(\psi_2,\psi_1)+ a_{21}(\psi_1,\psi_2))\\ 
&=& - \re \int_{\Omega} V(x) \left(\psi_1(x) \overline{\psi_2(x)}+\psi_2(x) \overline{\psi_1(x)} \right) dx\\
&=&0.
\end{eqnarray*}
Since the right hand-side of \eqref{coerc2x2} is larger than 0 if $\psi \neq 0 = \psi'$, the condition b) in Proposition \ref{propcoerc2x2} is not satisfied. 
As a consequence, the form $(a,V)$ is not coercive.
\end{exa}

%% file: sec_operator.tex
\section{Identification of the operator domain}\label{sec:operator}
Our main interest are evolution equations, and we will address this topic in details in Section~\ref{sec:evoleq}.
In that context, we will relate properties of the form $(a,V)$ to properties of solutions of evolution equations governed by the associated operator $(A,D(A))$, as described in Appendix~\ref{sec:sesquilinear}. 
Therefore, it is important to determine the domain of the operator associated with the sesquilinear form. This the object of this section.

In order to define the operator associated with the form $(a,V)$ we have to specify a Hilbert space $H$ such that $V \hookrightarrow H$ densely.
To obtain such a space in our context, we assume that the Hilbert spaces $V_i$ are uniformly and densely embedded into Hilbert spaces $H_i$, 
i.e.,  we choose a family $H_i$ such that $V_i \hookrightarrow H_i$ for all $i\in I$, and the norms of the canonical injections are bounded from above.
In order to fix the ideas, we state the standing assumptions.
\begin{ass}\label{injections}
During the rest of this chapter we always assume the following.
\begin{itemize}
\item The Hilbert spaces $V_i$ are continuously and densely embedded into the Hilbert spaces $H_i$ for all $i \in I$.
\item The Hilbert space $H$ is defined as
$$
H:=\bigoplus_{i \in I} H_i.
$$
\item The canonical injections are denoted $\phi_j: V_j \to H_j $.
\item The canonical, set theoretical injection $\phi: \prod_{j \in I} V_j \to \prod_{j \in I} H_j$ is defined componentwise by
$$
\pi_j(\phi(\psi)):=\phi_j(\pi_j(\psi)).
$$
\end{itemize}
\end{ass}
In the case of a finite set $I$ there are no difficulties connected with the boundedness of the injection $\phi$.
In the general case, however, the mapping $\phi$ is an unbounded operator from $V$ to $H$, see Example~\ref{multper}.
We recall some embedding results, proving them for the sake of completeness.
\begin{lemma}\label{contembed}
The following equivalences hold.
\begin{enumerate}
\item The injection $\phi$ is continuous if and only if the injections $\phi_j$ are uniformly continuous, i.e.,  if and only if there exists $M \geq 0$ such that
\begin{equation}\label{injcont}
\|\psi\|_{H_i} \leq M \|\psi\|_{V_i}, \qquad \mbox{for all } i \in I, \psi \in V_i.
\end{equation}
\item Assume $I$ to be infinite. The injection $\phi$ is a compact operator if and only if the injections $\phi_j$ are compact operators 
and  
$$
\lim_{j \to \infty} \|\phi_j\|_{\mathcal L(V_j,H_j)}=0.
$$
\end{enumerate}
\end{lemma}
\begin{proof}
a) If the estimate \eqref{injcont} holds, then 
$$
\|\psi\|^2_H = \sum_{i \in I} \|\psi_i\|^2_{H_i} \leq M \sum_{i \in I} \|\psi_i\|^2_{V_i} = M\|\psi\|^2_V
$$
yields the continuity of the injection $\phi$.

Assume now the injection $\phi$ to be continuous with constant $M$ and fix arbitrary $i\in  I$ and $\psi_i \in V_i$.
Convince yourself that the identity
$$
\phi_i(\psi_i)=\phi(\pir(\psi_i))
$$
holds and observe that $\pi^{-1,r}_{i,0}: H_i \to H$ is also isometric. Thus,
$$
\|\phi_i(\psi_i)\|_{H_i}= \|\phi(\pir(\psi_i))\|_H \leq M \| \pir(\psi_i) \|_V = M \|\psi_i\|_{V_i}.
$$

\smallskip
b) First assume that $\lim_{j \to \infty}\|\phi_j\|_{\mathcal L(V_j,H_j)} =0$ and approximate $V$ by $V^{(m)}:=\bigoplus_{j \leq m} V_j$ and $H^{(m)}$ analogously. 
Observe that the canonical injections $\phi^{(m)} \in \mathcal L(V^{(m)}, H^{(m)})$ are compact operators. 
Using this fact, we define a compact operator $\phi_m \in \mathcal K(V,H)$ by 
$$
\phi_m\psi :=\pi^{-1,r}_{\{1,\ldots,m\},0}(\phi^{(m)}(\pi_{\{1,\ldots,m\}}\psi)).
$$
This operator satisfies $\pi_i(\phi_m \psi)= \psi_i$, $i\leq m$, and $\pi_i(\phi_m \psi)= 0$, $i >m$.
Computing
$$
\|(\phi-\phi_m) \psi\|^2_H \leq \sum_{j \in I} \|\phi_j\|_{\mathcal L(V_j,H_j)}^2 \|\psi_j\|^2_{V_j} 
\leq \sup_{j \geq m}\|\phi_j\|^2_{\mathcal L(V_j,H_j)} \|\psi\|^2_V.
$$
shows that $\phi_m$ converges to $\phi$ in the operator norm.
Since $\mathcal K(V,H)$ is a closed ideal the first implication is proved.

\smallskip
Assume now that $V \hookrightarrow H$ compactly and fix an arbitrary infinite subset $I' \subset I$. 
Consider $\psi_j \in V_j$ such that 
$$\|\psi_j\|_{V_j}=1 \quad \mbox{and} \quad \|\psi_j\|_{H_j} \geq \frac{\|\phi_j\|_{\mathcal L(V_j,H_j)}}{2}.$$
For the sake of the simplicity of the notation denote for the rest of the proof $\psi^j:=\pjr(\psi_j)$.

Since $\pjr$ is isometric, $\psi^j$ lies in a bounded subset of $\bigoplus_{j \in I'} V_j$
and since $\bigoplus_{j \in I'} V_j \hookrightarrow \bigoplus_{j \in I'} H_j$ compactly, 
there exists a subsequence $j_k$ such that $\psi^{j_k}$ converges in $H$.

The sequence $(\psi^j)_{j \in I}$ converges pointwise to 0, i.e.,  $\lim_{j \to \infty} \pi_{k}(\psi^j)=0$ for all $k \in I'$.
In fact, each component is eventually 0. Thus, $\psi^{j_k}$ converges in $H$ to 0. 
Since $\pjr$ is isometric $\|\phi_{j_k}\|_{\mathcal L(V_{j_k},H_{j_k})}$ converges to 0, too. 

Since $I'$ is infinite, it is possible to identify it with the subsequence
$$
(\|\phi_{j_k}\|_{\mathcal L(V_{j_k},H_{j_k})})_{j_k\in I}, I'=\{j_k: k \in I\}
$$
of $(\|\phi_j\|_{\mathcal L(V_j,H_j)})_{j\in I}$.

We have proved that there exists a subsubsequence 
$$\|\phi_{j_{k_\ell}}\|_{\mathcal L(V_{j_{k_\ell}},H_{j_{k_\ell}})}, \qquad \ell=1,2,\ldots$$ converging to 0.
This means that $(\|\phi_j\|_{\mathcal L(V_j,H_j)})_{j\in I}$ converges to 0.
\end{proof}
It has to be stressed that the Hilbert space $H$ in which the form domain $V$ is embedded determines the operator associated with the sesquilinear form:
changing the state space $H$ and considering the same form $(a,V)$ leads to a different operator on $H$.
By this method,  sesquilinear forms can be used in order to investigate also second-order operators that are not in divergence form, 
see Section~\ref{sec:histremI} for references.
The following is an illustrative example of the problems that can arise with such techniques in our context.

\begin{exa}[Multiplicative perturbations of the Laplacian]\label{multper}
For a bounded domain $\Omega$ consider $V_i=H^1(\Omega)$ for all $i \in \mathbb N$. 
Let $M \in C^1(\overline{\Omega})$ be a strictly positive, bounded function whose strictly positive minimum is given by $m:=\min_{x\in \Omega} M(x)$. 
Denote $(\Delta_{Neumann}, D(\Delta_{Neumann}))$ the Laplace operator with Neumann boundary condition on $\Omega$ and define the operators $(A_i,D(A_i))$ by
$$
A_if(x):=  \frac{M(x)}{i}\Delta_Nf (x), \quad  D(A_i):=D(\Delta_{Neumann}).
$$
Set $H_i:=L^2(\Omega,\frac{i}{ M(x)}d\lambda )$. Then the operators $A_i$ are associated 
with the sesquilinear forms $a_i: H^1(\Omega) \times H^1(\Omega) \to \mathbb C$ 
$$
a_i(\psi,\psi'):= \int_\Omega \nabla \psi(x) \overline{\nabla \psi'(x)} dx.
$$
Although $V_i \hookrightarrow H_i$ for all $i \in \mathbb N$,  $V$ is not continuously embedded into $H$.
As a consequence $V$ is not suitable as a form domain for  the form $(a,V)$ in the renormed Hilbert space $H$.
In fact, this is due to the degeneracy of the diagonal multiplication operator $\frac{M}{i}\otimes Id$.
\end{exa}

We want to identify the domain of the operator associated with the form $(a,V)$.
In the following we give a possible characterisation of the domain.
\begin{prop}\label{operdomain}
Consider a continuous, densely defined form $(a,V)$ such that $V \hookrightarrow H$. 
The domain of the associated operator $(A,D(A))$ satisfies
$$
D(A)=\left\{(\psi_i)_{i \in I} \in V: \exists f \in H, \forall \psi' \in V , i \in I:
\sum_{j \in I} a_{ij}(\psi_j,\psi'_i)= (f_i,\psi'_i)_{H_i} \right\}
$$
\end{prop}
\begin{proof}
Denote by $X$ the space in the claim. By the definition, the domain of the operator $D(A)$ is given by
$$
D(A):=\left\{ \psi \in V: \exists f \in H, \forall \psi' \in V : a(\psi,\psi')=(f,\psi')_H  \right\}.
$$
In particular, if $\psi \in D(A)$ there exists a $f$ as in the above definition for all $\psi':=\pir(\psi'_i)$ with $\psi'_i \in V_i.$ 
Plugging $\psi'$ into the expression of the form $(a,V)$ shows that $\psi \in X.$

\smallskip
Conversely, fix $f \in H.$
Fix an arbitrary $\psi'  \in V$. Using the embedding of $V \hookrightarrow H$, we see that the sum $\sum_{i\in I} (f_i, \psi'_i)=(f,\psi)$ converges. 
Further,
$$
(f,\psi')=\sum_{i\in I} (f_i,\psi'_i)
=\sum_{i \in I} \sum_{j\in I} a_{ij}(\psi_j,\psi'_i), 
= a(\psi,\psi')
\quad \mbox{ for all } \psi' \in V
$$
shows that $\psi \in D(A)$.
\end{proof}

%% file: sec_evoleq.tex
\section{Evolution equations}\label{sec:evoleq}
During this section we always assume that the Assumptions~\ref{injections} hold and that the canonical injections $\phi_j$ are continuous, uniformly on $j$. 
In the Hilbert space $H$, consider the abstract Cauchy problem
\begin{equation}\tag{ACP}\label{ACP}
\left\{\begin{array}{rcll}
\frac{d}{dt}{\psi}(t) &=& A\psi(t), & t\geq 0, \\
\psi(0)&=&f,& f\in H,
\end{array}\right.
\end{equation}
where the operator $(A,D(A))$ is the operator associated in $H$ with a sesquilinear form $(a,V)$. 
Then, it is possible to deduce properties of the solution of the equations from properties of the sesquilinear forms. 
In particular, we are interested in deducing properties of the solution of \eqref{ACP} by arguments of linear algebraic type 
applied on properties of the \emph{single} mappings $a_{ij}$.
All general results about sesquilinear forms that we will need in the following are stated without proof in Appendix~\ref{sec:sesquilinear}.

We start by discussing the well-posedness of the problem \eqref{ACP}.
By Theorem \ref{WPform} the operator $(A,D(A))$ associated with a form $(a,V)$ 
generates an analytic semigroup $\etasg$ on $H$ if the form $(a,V)$ is $H$-elliptic, 
i.e., if there exists $\omega \in \mathbb R$ such that
$$
\re a(\psi) + \omega \|\psi\|^2_H \geq \alpha \|\psi\|^2_V
$$
holds for all $\psi \in V$. In this case the system \eqref{ACP} is well-posed and the solutions $\psi(t)$ to the initial data $f$
are given by $e^{ta}f$ for all $t \geq 0$.

Our first goal is to give sufficient conditions for the ellipticity. 
We use the same idea as in Lemma \ref{boundedform} and Proposition \ref{coercivform}.
Assume that the form $(a,V)$ is continuous. Since the form $a_{ij}$ are also continuous for all $i,j$, there exist $\alpha_{ij}<0, \omega_{ij} \in \mathbb R$ 
such that the estimate
\begin{equation}\label{continsingle}
|a_{ij}(\psi,\psi')| \leq -\alpha_{ij}\|\psi\|_{V_j}\|\psi'\|_{V_i} + \omega_{ij} \|\psi\|_{H_j}\|\psi'\|_{H_i}
\end{equation}
holds for all $\psi \in V_j, \psi' \in V_i$ and all $i \neq j$. 
Assume now that the forms $a_{ii}$ are $H_i$-elliptic. Thus, there exist $\alpha_{ii} > 0, \omega_{ii} \in \mathbb R$ such that
\begin{equation}\label{ellipticsingle}
\re a_{ii}(\psi) + \omega_{ii} \|\psi\|_{H_i}^2 \geq \alpha_{ii}\|\psi\|^2_{V_i} 
\end{equation}
holds for all $\psi \in V_i$. Combining these facts, we obtain the following result.
\begin{prop}[Well-posedness of \eqref{ACP}]\label{ellipticityform}
Assume the following.
\begin{itemize}
\item The form $(a,V)$ is continuous and \eqref{continsingle} holds.
\item The forms $a_{ii}$ are elliptic and \eqref{ellipticsingle} holds.
\item The matrix $A=(\alpha_{ij})_{i,j \in I}$ defines a positive definite operator, i.e., there exists $\alpha >0$ such that
$\re (A v \mid v) \geq \alpha \|v\|^2_{\ell^2(I)}$ for all $v \in \ell^2(I)$.
\item The matrix $\Omega=(\omega_{ij})_{i,j \in I}$ defines an operator in $\ell^2(I).$
\end{itemize}
Then $(a,V)$ is elliptic and satisfies
$$
\re a(\psi) + \|\Omega\|_{\mathcal L(\ell^2(I))}\|\psi\|^2_H \geq \alpha \|\psi\|^2_V.
$$
\end{prop}
\begin{rem}
If the form $(a,V)$ is elliptic, computing $a(\psi_i)=a_{ii}(\psi)$ shows that also $a_{ii}$ is elliptic, with the same constants.
Thus, the ellipticity with uniform constants of $a_{ii}$ is a necessary condition for the form $(a,V)$ to be elliptic.
\end{rem}
\begin{proof}[Proof of Prop.~\ref{ellipticityform}]
The ellipticity is a consequence of the estimate
\begin{eqnarray*}
\re a(\psi)&=&\re \sum_{i \in I} a_{ii}(\psi_i)+ \re \sum_{j \not= i} a_{ij}(\psi_j,\psi_i)\\
&\geq &\alpha_{ij} \|\psi_j\|_{V_j} \|\psi_i\|_{V_i} - \omega_{ij} \|\psi_j\|_{H_j}\|\psi_i\|_{H_i}\\
&\geq& \alpha \|\psi\|^2_V - \|\Omega\| \|\psi\|^2_H.
\end{eqnarray*}
\end{proof}
 
We now discuss a second way to prove ellipticity. 
Assume the diagonal entries to be uniformly elliptic and the off-diagonal entries to be bounded as mappings $a_{ij}:H_j \times H_i \to \mathbb C$.
Then the resulting matrix of forms will be elliptic. In fact, it holds also if the off-diagonal entries are bounded in some interpolation space.
Recall that an interpolation space between $V$ and $H$ of order $\alpha \in [0,1)$ is any linear space $H^\alpha$ such that
$$
V \hookrightarrow H^\alpha \hookrightarrow H, \quad \mbox{and} \quad   \|\psi\|^\alpha_{V}\|\psi\|^{1-\alpha}_{H} \leq M\|\psi\|_{H^\alpha}, \quad \psi \in V.
$$
For finite index sets, Lemma ~\ref{delioperturb} yields a criterion for the ellipticity of the form $(a,V)$.

\begin{prop}\label{perturbinter}
For a finite index set $I$, assume that $H^\alpha=\bigoplus_{i \in I} H^\alpha_i$ is an interpolation space between $V$ and $H$ 
and that the mappings $a_{ij}:H^\alpha_j \times V_i \to \mathbb C$ are continuous. If $(a,V)$ is continuous and $a_{ii}$ are elliptic for all $i \in I$, the form $(a,V)$ is elliptic.
\end{prop}
\begin{proof}
Observe that for finite index sets $I$, $H^\alpha:=\bigoplus_{i\in I}H^\alpha_i$ is an interpolations space of order $\alpha$ between $V$ and $H$, 
whenever $H^\alpha_i$ are such between $V_i$ and $H_i$.

So the sesquilinear mapping $b=(a_{ij})_{i\neq j } : H^\alpha \times V \to \mathbb C$ is bounded. 
Thus, we can apply Lemma~\ref{delioperturb} and see that the form $(a,V)$ is elliptic.
\end{proof}
We now turn our attention to the cosine operator functions, see Appendix~\ref{sec:sesquilinear}. 
We give a sufficient condition for $(A,D(A))$ to generate a cosine operator function. 
Two applications of this result can be found in Section~\ref{sec:wave} and Section~\ref{sec:dynamic}.
\begin{prop}\label{cosine}
Consider a finite index set $I$ and an elliptic, continuous form $a=(a_{ij})_{i,j\in I}$.
The following assertions hold.
\begin{enumerate}[(1)]
\item Assume that there exists $M\geq 0$ such that for all $\psi\in V_j$ and $\psi'\in V_i$
\begin{enumerate}[(i)]
\item $|\im a_{ii} (\psi)|\leq M \| \psi \|_{V_i} \| \psi \|_{H_i}$ for all $i\in I$ and moreover 
\item There exists a set  $J \subset K:=\{(i,j) \in I \times I:  i\neq j \}$ such that
\begin{itemize}
\item for all $i\neq j$ either $(i,j) \in J$ or $(j,i) \in J$, and
\item either $|\im (a_{ij}(\psi,\psi')+{a_{ji}(\psi',\psi)})|\leq M\| \psi\|_{V_j}\| \psi'\|_{H_i}$, for all $(i,j) \in J,$ 
\item or $|\im (a_{ij}(\psi,\psi')+{a_{ji}(\psi',\psi)})|\leq M\| \psi'\|_{H_j}\| \psi \|_{V_i}$, for all $(i,j) \in K\setminus J.$
\end{itemize}
\end{enumerate}
Then the operator $A$ associated with $(a,V)$ generates a cosine operator function with associated phase space $V \times H$. In particular, $A$ generates an analytic semigroup of angle $\frac{\pi}{2}$ on $H$.
\end{enumerate}
\end{prop}
%\item Conversely, if $A$ generates a cosine operator function, then for all $i \in I$ also the operator $A_{ii}$ associated with $a_{ii}$ generates a cosine operator function.
\begin{proof}
We want to apply Proposition~\ref{crouzeix}.
Under the assumptions in (1), we have
\begin{eqnarray*}
| \im a (\psi)|&\leq& |\sum_{i\in I} \im a_{ii}(\psi_i,\psi_i)|
+|\sum_{i\not=j} \im a_{ij}(\psi_j,\psi_i)|\\
&= &|\sum_{i\in I} \im a_{ii}(\psi_i,\psi_i)|+|\sum_{(i,j) \in J} \im (a_{ij}(\psi_j,\psi_i)+a_{ji}(\psi_i,\psi_j))|\\
&\leq &\sum_{i\in I} |\im a_{ii}(\psi_i,\psi_i)|+\sum_{(i,j) \in J} |\im (a_{ij}(\psi_j,\psi_i)+a_{ji}(\psi_i,\psi_j))|\\
&\leq &\left\{\begin{array}{l}
M \sum_{i \in I} \| \psi_i\|_{V_i} \| \psi_i\|_{H_i} + M\sum_{(i,j)\in J} \| \psi_j\|_{V_j}\| \psi_i\|_{H_i}\\
M \sum_{i\in I}  \| \psi_i\|_{V_i} \| \psi_i\|_{H_i} + M\sum_{(i,j)\in J} \| \psi_i\|_{V_i}\| \psi_j\|_{H_j}
 \end{array}\right.\\
&\leq& \tilde{M}\| \psi\|_{V} \| \psi\|_{H}
\end{eqnarray*}
for some constant $\tilde{M}\geq 0$. So, the first claim follows from Proposition~\ref{crouzeix}.
The second claim is a consequence of~\cite[Thm.~3.14.17]{AreBatHie01}.
\end{proof}
%In order to prove (2) recall that $(A,D(A))$ generates a cosine operator function, 
%then the numerical range 
%$$
%W(a):=\{a(\psi): \psi \in V, \|\psi\|_H=1\}
%$$ 
%lies in a parabola for some equivalent scalar product on $H$. 
%Such a scalar product on $H$ induces an equivalent scalar product on $H_i$, too, and therefore also the numerical range $W(a_{ii})$ lies in the same parabola, for all $i\in I$. 
%Thus, $A_{ii}$ generates a cosine operator function by Crouzeix--Haase's result.
These criteria give sufficient conditions to prove well--posedness of the problem \eqref{ACP}.
Our next aim is to characterise properties of the solutions. 
First, we show under which conditions the semigroup generated by a form $(a,V)$ is real, positive or $L^\infty$-contractive. 
To this end we restrict ourselves to the case $H_i=L^2(X_i)$ and we introduce the set
\begin{eqnarray*}
C^\infty_i& := & B^\infty_{X_1}\times \ldots B^\infty_{X_{i-1}}\times L^2(X_i)\times B^\infty_{X_{i+1}}\times\ldots \\
&= &\{(\psi_j)_{j \in I}:\psi_j \in B^\infty_{X_j} \mbox{ for all } j\neq i, \mbox{ and }\psi_i \in L^2(X_i) \}.
\end{eqnarray*}
Here we denote by $B^\infty_X$ the unitary ball of the space $L^\infty(X)$.

\begin{thm}\label{positivity}
Consider a family of Hilbert spaces $H_i:=L^2(X_i), V_i \hookrightarrow H_i$ uniformly, such that $X_i$ is a $\sigma$-finite measure space for all $i \in I$.
Assume that $(a,V)$ is a continuous, elliptic and accretive form. The following assertions hold.
\begin{enumerate}
\item The semigroup $\etasg$ is real  if and only the two following conditions hold
\begin{itemize}
\item For all $i \in I$, $\psi \in V_i \Rightarrow \re \psi \in V_i.$
\item For all $i,j \in I$ the mapping $a_{ij}$ is real, i.e., $\im a_{ij}(\psi,\psi')=0 $ for all real functions  $\psi \in V_j,\psi' \in V_i.$
\end{itemize}
\item If the semigroup $\etasg$ is real, then it is positive if and only if the two following conditions hold
\begin{itemize}
\item For all $i \in I$ the semigroup $(e^{ta_{ii}})_{t \geq 0}$ is positive.
\item For all $i\neq j$ the sesquilinear mapping $-a_{ij}$ is positive, i.e.,
$$
a_{ij}(\psi,\psi') \leq 0, \qquad \mbox{for all } \psi \in V_j^+,\psi' \in V_i^+.
$$
\end{itemize}
\item If the semigroup $\etasg$ is contractive, then it is $L^\infty$-contractive if and only if the two following conditions hold
\begin{itemize}
\item For all $i \in I$, if $\psi \in V_i, (1\wedge |\psi|)\sign \psi \in V_i$
\item  The estimate
\begin{eqnarray*}
\sum_{j\neq i}|a_{ij}(\psi_j,(|\psi_i|-1)^+ \sign \psi)| \\
\leq \re a_{ii}((|\psi_i| \wedge 1)\sign \psi_i,(|\psi_i|-1)^+ \sign \psi)
\end{eqnarray*}
holds for all $ \psi \in V \cap C_i^\infty$.
\end{itemize}
\end{enumerate}
\end{thm}
We recall that by definition $\psi \in H$ is real, respectively, positive, if and only $\psi_i$ is almost everywhere real, respectively, positive, for all $i \in I.$
\begin{proof}
In order to apply the characterisations of real, positive and $L^\infty$-contractive semigroups arising from Theorem~\ref{ouh} (see also \cite{Ouh04}),
represent $H$ as $L^2(X)$ where $X$ is a $\sigma$-finite measure space, as discussed in Lemma~\ref{Hrepresent}.

\smallskip
a) Reality. First assume that both conditions in the claim hold. Observe that the projection on real cone of $L^2(X)$ is given by $P(\psi_i)_{i \in I}=(\re \psi_i)_{i \in I}$. We have to prove $\psi \in V \Rightarrow \re \psi \in V$, but this is clear since $V$ is the direct sum and the first condition holds. Further, computing
$$
\im a(\re \psi, \im \psi') = \sum_{i,j \in I} \im a_{ij}(\re \psi_j, \im \psi'_i)=0
$$
shows that also the second condition in Theorem \ref{ouh} holds.

\smallskip
b) Positivity. We first observe that the projection on the positive cone of the Hilbert space $H$ is given by $P(\psi_i)_{i \in I}=((\re \psi_i)^+)_{i \in I}$. Denote $P_i$ the projection of $H_i$ on its positive cone. Then, the domain condition in \ref{ouh} is equivalent to the condition $P_iV_i \subset V_i$ for all $i \in I.$

In order to prove the equivalence to the algebraic condition of Theorem~\ref{ouh}, we compute
\begin{eqnarray*}
\re a(P\psi, (I-P) \psi &=& - \re \sum_{i,j\in I} a_{ij}(\re \psi_j^+, \re \psi_i^-)\\
&=&-\left( \re\sum_{i \in I} a_{ii}(\re \psi_i^+, \re \psi_i^-) \right.\\
&& \left. + \re\sum_{i\neq j} a_{ij}(\re \psi_j^+, \re \psi_i^-) \right).
\end{eqnarray*}
The conditions stated in the theorem are sufficient. To prove that they are also necessary, fix arbitrary $\psi \in H_j, \varphi \in H_i$, recall that setting $a_{ij}(\psi , \varphi) = a(\pjr(\psi), \pir(\varphi))$ yields a representation for $(a,V)$ and distinguish the following two cases. If $i=j$, then
$$
0 \leq \re a(\re( \pir(\psi))^+, \re(\pir(\varphi))^- ) = a_{ii}(\re (\pir(\psi))^+ , \re(\pir(\varphi))^-)
$$
is equivalent to the algebraic condition for the positivity of the semigroups $(e^{ta_{ii}})_{t \geq 0}$, due to the surjectivity and positivity of $\pir.$ If, conversely $i\neq j$, then 
$$
0 \leq \re a(\re( \pjr(\psi))^+, \re(\pir(\varphi))^- ) = a_{ij}(\re (\pjr(\psi))^+ , \re(\pir(\varphi))^-)
$$
is equivalent to off-diagonal condition in the theorem, against by the surjectivity and positivity of $\pir, \pjr$.

\smallskip
c) Contractivity in $L^\infty(X)$. We use~\cite[Thm.~2.14]{Ouh04}. The semigroup is $L^\infty$-contractive if and only if 
\begin{equation}\label{ouhcon}
\begin{array}{c}
\psi \in V \Rightarrow (1\wedge | {\psi}|)\sign \psi \in V\\ 
\;\mbox{and}\;\\
\re a((1\wedge | \psi|)\sign \psi,(| \psi|-1)^+ \sign \psi)\geq 0.
\end{array}
\end{equation}
One sees that $\psi\in  V \Rightarrow (1\wedge | {\psi}|)\sign \psi \in V$ 
if and only if $ \psi \in V_i\Longrightarrow(1\wedge |\psi|)\sign \psi\in V_i$ for all $i\in I$. 
We have to prove the equivalence of the estimates in b) and \eqref{ouhcon}. Let first $\psi\in C^\infty_i$. Then
$$
(1\wedge |\psi|) \sign \psi = \pi^{-1,r}_{i, \psi} (1 \wedge \psi_i \sign \psi_i)
$$
and 
$$
(| \psi|-1)^+ \sign \psi= \pir((|\psi_i|-1)^+\sign \psi_i).
$$ 
Accordingly, 
\begin{eqnarray*}
0&\leq& \re a((1\wedge | \psi|)\sign\psi , (| \psi|-1)^+ \sign \psi)\\
&=&\sum_{j\not=i} \re a_{ij}(\psi_j,(|\psi_i|-1)^+\sign \psi_i)\\
&&+\re a_{ii}((1\wedge|\psi_i|),(|\psi_i|-1)^+\sign \psi_i)
\end{eqnarray*}
for all $\psi \in  {C}^\infty_i\cap { V}$ and all $i\in I$. Due to the sesquilinearity of $a_{ij}$, this also implies
$$0\leq \sum_{j\not=i}\re a_{ij}(\pm \psi_j,(|\psi_i|-1)^+\sign \psi_i)
+\re a_{ii}((1\wedge|\psi_i|),(|\psi_i|-1)^+\sign \psi_i)$$
for all $\psi \in { C}^\infty_i\cap {V}$, all $i\in I$, and all $\alpha\in{\mathbb C}$, $|\alpha|\leq 1$.
This yields the claimed criterion. The converse implication can be proven analogously.
\end{proof}

\begin{rems}
(1) Perturbing the form $(a,V)$ by a multiple of the scalar product, i.e., defining $a_\omega : V \times V \to \mathbb C$ by
$$
a_\omega(\psi,\psi') := a(\psi,\psi') - \omega(\psi \mid \psi')_H
$$
leads to a semigroup $(e^{ta_\omega})_{t \geq 0}= (e^{\omega t} e^{ta})_{t \geq 0}$.
 
Observe that  since $e^{\omega t} > 0$ for all $\omega, t \in \mathbb R$ the semigroup $e^{ta}$ is real, respectively, positive, if and only if
the semigroup $e^{ta_\omega}$ is real, respectively, positive.
Since every elliptic form is also accretive up to a translation, 
one sees that the accretivity condition in Theorem~\ref{positivity} can be dropped in part a) and b).

(2) In the case of systems, i.e., in the case $X_i=\Omega$ for all $i \in I$, the off-diagonal mappings are in fact sesquilinear forms.
So, the condition of positivity becomes quite restrictive.
For a discussion of the positivity of sesquilinear forms see Section~\ref{sec:histremI}.
\end{rems}

The adjoint of the operator $A$ is associated with the form $a^*: V\times V \to \mathbb C$ defined by
$a^*(\psi, \psi'):= \overline{ a(\psi',\psi)}$. 
If $(a^*,V)$ is contractive in $L^\infty$, then the semigroup $e^{ta}$ extrapolates to a family
of strongly continuous semigroups on all spaces $L^p(X)$. %riferimento bibliografico?
Thus, applying the criterion c) in Proposition~\ref{positivity} to the adjoint form $(a^\star,V)$,
we obtain a characterisation of form matrices generating an extrapolating semigroup.

\begin{cor}\label{ultra}
Consider a family of Hilbert spaces $H_i:=L^2(X_i), V_i \hookrightarrow H_i$ uniformly, such that $X_i$ is a $\sigma$-finite measure space for all $i \in I$.
Assume $a=(a_{ij})_{i,j \in I}: V \times V \to \mathbb C$ to be accretive. 
Then the semigroup $(e^{ta})_{t\geq 0}$ extrapolates to a family of contractive $C_0$-semigroups $(e^{ta_p})_{t \geq 0}$ on $L^p(X)$ for all $p\in[1,\infty)$
if and only if for all $i\in I$
\begin{itemize}
\item $\psi\in V_i\Longrightarrow(1\wedge |\psi|)\sign \psi \in V_i;$
\item For all $\psi \in  V \cap { C}^\infty_i,$
\begin{eqnarray*}
&&\sum_{j\not= i} |a_{ij}(\psi_j,(|\psi_i|-1)^+\sign \psi_i)|\\
&&\leq \re a_{ii}((1\wedge|\psi_i|)\sign \psi_i,(|\psi_i|-1)^+\sign \psi_i);
\end{eqnarray*}
\item  For all $\psi\in { V}\cap {C}^\infty_i $
\begin{eqnarray*}
&&\sum_{j\not= i} |a_{ji}((|\psi_i|-1)^+\sign\psi_i,\psi_j)|\\
&&\leq \re a_{ii}((|\psi_i|-1)^+\sign \psi_i,(1\wedge|\psi_i|)\sign \psi_i).
\end{eqnarray*}
\end{itemize}
\end{cor}

\begin{proof}
Let us first assume the semigroup $(e^{ta})_{t\geq 0}$ to extrapolate to a family of contractive $C_0$-semigroups on $L^p(X)$, $p\in[1,\infty)$, 
and hence in particular to be $L^\infty$-contractive.
%riferimento bibliografico? 
Since also the unit ball of $L^1(X)$ is left invariant, it follows by duality that the semigroup $(e^{ta^*})_{t\geq 0}$ is $L^\infty$-contractive. 
Here $a^*$ denotes the adjoint form of $(a,V)$, which by definition is given by
${a}^*(\psi,\psi')=\overline{{a}(\psi',\psi)}=\sum_{i,j=1}^m \overline{a_{ji}(\psi'_i,\psi_j)}$, $\psi,\psi'\in{V}$. 
Since $a^*$ is accretive if and only if $(a,V)$ is accretive, 
we can apply Theorem~\ref{positivity} to $(e^{ta})_{t\geq 0}$ and $(e^{ta^*})_{t\geq 0}$ and obtain the conditions in the claim.

Conversely, since both $(a,V)$ and $a^*$ are accretive, it follows from the first two conditions 
and Theorem~\ref{positivity} that $(e^{ta})_{t\geq 0}$ is $L^\infty$-contractive. 
Moreover, since also $a^*$ is accretive, it follows from the above conditions and Theorem~\ref{positivity} 
that $(e^{ta^*})_{t\geq 0}$ is $L^\infty$-contractive, too. 
Thus, by standard interpolation and duality arguments $(e^{ta})_{t\geq 0}$ extrapolates 
to a family of contractive semigroups on $L^p(X)$, $p\in [1,\infty)$, that are strongly continuous for all $p>1$. 
Finally, contractivity implies strongly continuity of the extrapolated semigroup also in $L^1(X)$, cf.~\cite[7.2.1]{Are04}.
\end{proof}
A second corollary of the above theorem is the characterisation of \emph{ultracontractive} semigroups. 
Assume the semigroup $(T(t))_{t\geq 0}$ to have consistent realisations on all $L^p(X)$ spaces, $p\in [1,\infty]$.
Then the semigroup is \emph{ultracontractive of dimension $d$} if there is a constant $c> 0$ such that for all $p,q\in [1,\infty], f\in L^p(X)$ 
the estimate
$$
\|T(t)f\|_{L^q(X)} \leq c t^{\frac{-d}{2}|p^{-1}-q^{-1}|} \|f\|_{L^p(X)}
$$  
holds.
\begin{cor}\label{ultra2}
Assume that the three conditions in Corollary~\ref{ultra} hold. Let $d > 2$ a real number.
Then $\etasg$ is ultracontractive of dimension $d$ if and only if $V_i$ is continuously embedded in $L^{\frac{2d}{d-2}}(X_i)$ for all $i \in I$.
\end{cor}
\begin{proof}
This is a consequence of Corollary~\ref{ultra} and \cite[Thm.~7.3.2]{Are04}.
\end{proof}

Again using criteria derived from Theorem \ref{ouh} it is possible to investigate domination properties of semigroups.
Recall that the semigroup $(e^{ta})_{t \geq 0}$ \emph{dominates} the semigroup $(e^{tb})_{t\geq 0}$ if $\etasg$ is a positive semigroup
and if for all $t\geq 0$ $e^{ta} \geq |e^{tb}|$.

\begin{prop}\label{domination}
Consider a family of Hilbert spaces $H_i:=L^2(X_i), V_i \hookrightarrow H_i$ uniformly, such that $X_i$ is a $\sigma$-finite measure space for all $i \in I$.
Assume $a=(a_{ij})_{i,j \in I}: V \times V \to \mathbb C$ to be continuous, elliptic and accretive and $(e^{ta})_{t\geq 0}$ to be positive. 
Consider another densely defined, continuous, ${H}$-elliptic sesquilinear form $b:=(b_{ij})_{i,j \in I}:{W}\times{W}\to {\mathbb C}$, 
${W}=\bigoplus_{i\in I} W_i$. 
Then $(e^{ta})_{t\geq 0}$ dominates $(e^{tb})_{t\geq 0}$ if and only if the following conditions hold.
\begin{itemize}
\item $W_i$ is an ideal of $V_i$ for all $i \in I$,
\item $\re b_{ii}(\psi,\psi') \geq a_{ii}(|\psi|,|\psi'|)$ for all $\psi,\psi'\in V_i$ such that $\psi,\overline{\psi'}\geq 0,\; i\in I$, and
\item $|\re b_{ij}(\psi,\psi')| \leq -a_{ij}(|\psi|,|\psi'|)$ for all $\psi \in V_j, \psi'\in V_i, i\in I.$
\end{itemize}
\end{prop}
\begin{proof}
First recall that from Theorem \ref{ouh} it is possible to derive the following equivalent conditions for domination:
\begin{itemize}
\item ${\mathcal W}$ is an ideal of ${V}$;
\item $\re b (\psi,\psi')\geq a (|\psi|,|\psi'|)$ for all $\psi,\psi'\in  W$ such that $\psi\overline{\psi'}\geq 0$.
\end{itemize}
Assume that $\etasg$ dominates $\etbsg$.

Since $W$ has diagonal form, $W$ is an ideal of $V$ if and only if each component $W_i$ is an ideal of $V_i$.

Let $i_0=j_0$ and $\psi,\psi' \in V_{i_0}$ such that $\psi\overline{\psi'}\geq 0$. So,
$$\pi^{-1,r}_{i_0,0}(\psi)\overline{\pi^{-1,r}_{i_0,0}(\psi')}\geq 0.$$
Computing 
\begin{eqnarray*}
\re b_{i_0i_0}(\psi,\psi')&=&\re b (\pi^{-1,r}_{j_0,0}(\psi),\pi^{-1,r}_{i_0,0}(\psi'))\\
&\geq& a(|\pi^{-1,r}_{j_0,0}(\psi)|,|\pi^{-1,r}_{i_0,0}(\psi')|)\\
&=& a_{i_0i_0}(|\psi|,|\psi'|)
\end{eqnarray*}
shows that the second condition is necessary. 

For $i_0\not=j_0$, let $\psi \in V_{j_0}$ and $\psi'\in V_{i_0}$, so that 
$\pi^{-1,r}_{j_0,0}(\psi)\overline{\pi^{-1,r}_{i_0,0}(\psi')}=0
=(-\pi^{-1,r}_{j_0,0}(\psi))\overline{\pi^{-1,r}_{i_0,0}(\psi')}$. 
Then, 
\begin{eqnarray*}
\pm \re b_{i_0j_0}(\psi,\psi')&=&\re b_{i_0j_0}(\pm \psi,\psi')\\
&=& \re b (\pm \pi^{-1,r}_{j_0,0}(\psi),\pi^{-1,r}_{i_0,0}(\psi'))\\
&\geq& a(|\pi^{-1,r}_{j_0,0}(\psi)|,|\pi^{-1,r}_{i_0,0}(\psi')|)\\
&=& a_{i_0j_0}(|\psi|,|\psi'|),
\end{eqnarray*}
thus proving that also the third condition is necessary. 

\smallskip
To check the converse implication let $\psi,\psi' \in W$ and compute 
$$
\re b(\psi,\psi') = \re \sum _{i,j\in I} b_{ij}(\psi_j,\psi'_i)\geq\sum _{i,j\in I} a_{ij}(|\psi_j|,|\psi'_i|)=a(|\psi|,|\psi'|).
$$
The proof is now complete since we have already shown that the conditions involving ideals are equivalent.
\end{proof}

%% file: sec_symmetriesmatrix.tex
\section{Symmetry properties}\label{sec:symmetriesmatrix}
In this section we are going to study \emph{symmetry} properties of the semigroup $\etasg$ in the case that $H_i=L^2(X)$ for all $i \in I$.
Recall that $H=\bigoplus_{i \in I} H_i$. Thus, in this case $H$ can be interpreted as $H=L^2(X,\ell^2(I))$, as discussed in Lemma~\ref{Hrepresent}, and, 
in this way, functions in $H$ are vector-valued functions.
The word \emph{symmetry} is used to denote the invariance of a particular class of subspaces of the state space $H$,
whose connection to the physical use of the word symmetry is discussed in Section~\ref{sec:histremII}.

Fix a Hilbert basis $(e_i)_{i \in I}$ of the space $\ell^2(I)$.
Since all functions $\psi \in H$ are vector-valued it is possible to decompose such functions using the Hilbert basis of the space $\ell^2(I)$.
\begin{prop}\label{L2decomposition}
Consider a family of Hilbert spaces $H_i:=L^2(X)$, such that $X$ is a $\sigma$-finite measure space
and consider a Hilbert basis $(e_i)_{i \in I}$ of the space $\ell^2(I)$. 
Then, for all $\psi \in H$ there exist uniquely determined functions $c_i \in L^2(X), i \in I$ that satisfy following conditions.
\begin{itemize}
\item The decomposition
\begin{equation}\label{fourier}
\psi(x)=\sum_{i \in I} c_i(x) e_i
\end{equation}
holds in $\ell^2(I)$ for almost every $ x\in \Omega $.
\item The estimate $\|c_i\|_{L^2(X)}\leq \|\psi\|_H$ holds for all $i \in I$.
\end{itemize}
The uniquely determined coefficients are denoted by
\begin{equation}\label{hilbertcoeff}
c^\psi_i(x):= c_i(x).
\end{equation}
\end{prop}

\begin{proof}
Fix a $\psi \in H$. By the definition of the norm in $H$,
$$
\|\psi\|^2_H=\sum_{i \in I} \|\psi_i\|^2_{L^2}=\sum_{i \in I} \int_\Omega |\psi_i(x)|^2 d\mu=\int_\Omega \sum_{i \in I}|\psi_i(x)|^2 d\mu.
$$
It follows that $\psi(x) \in \ell^2(I)$ for almost every $x\in \Omega$. Define 
$$c_i(x):=(\psi(x) \mid e_i)_{\ell^2}.$$
For almost every $x \in \Omega$  decompose
$$
\psi(x)= \sum_{i \in I}c_i(x) e_i,
$$
and we have obtained (a).
We obtain (b) using the Cauchy--Schwarz inequality and estimating
$$
\|c_i\|_{L^2(X)} = \int_\Omega |c_i(x)|^2 d \mu \leq \int_\Omega \|\psi(x)\|^2_{\ell^2} d \mu= \|\psi\|.
$$
To see that $c_i$ is a measurable function for all $i \in I$, observe that if $\chi_n$ is a sequence of simple functions converging to $\psi$, 
then for all $i \in I$, $\chi^i_n$ defined by
$$
\chi^i_n(x):=(\chi_n(x) \mid e_i), \qquad i,n \in \mathbb N
$$
is a sequence of simple functions converging to $c_i$.

Assume now that a second family $(d_i)_{i \in I} \subset L^2(X)$ satisfies both conditions. Then it follows
$$
(c_i(x)-d_i(x)) e_i = 0 \quad \mbox{for all $i\in I$ and almost all $x \in X$.}
$$
Denote by $N_i$ the exception set of the above equality, i.e., 
$$
N_i:=\{x \in X: (c_i(x)-d_i(x)) e_i \neq 0 \}.
$$
Since $N_i$ has measure 0 for all $i \in I$ and $I$ is countable, also $N:=\cup_{i\in I} N_i$ has measure 0.
So, $\sum_{i \in I} (c_i(x)-d_i(x)) e_i = 0$ pointwise almost everywhere.

By the dominated convergence theorem for vector-valued functions, we obtain $c_i=d_i$.

Since $(e_i)_{i \in I}$ is a Hilbert basis, $c_i(x)=d_i(x)$ almost everywhere in $X$ and the claim follows.
\end{proof}
We want to investigate the invariance of subspaces of $H$ induced by subspaces of $\ell^2(Y)$.
\begin{prop}
Fix a closed linear subspace $Y \subset \ell^2(I)$. The subspace $\mathcal Y \subset H$ defined by
$$
\mathcal Y:=\{\psi \in H : \psi(x) \in Y \mbox{ almost everywhere in } X\}
$$
is closed.
\end{prop}
\begin{proof}
Consider a sequence $\psi_n$ converging to $\psi$ in $H$. 
So, there exists a subsequence such that
$$
\lim_{k\to\infty} \psi_{n_k}(x) =\psi(x)\qquad \mbox{almost everywhere in } X.
$$
In particular, $\psi(x) \in Y$ almost everywhere since $Y$ is closed.
\end{proof}
Since $\mathcal Y$ is a closed subspace it is possible to decompose $H=\mathcal Y \bigoplus \mathcal Y^\perp$.
One sees that
$$
\mathcal Y^\perp = \{\psi \in H: \psi(x) \in Y^\perp \mbox{almost everywhere in } X \}.
$$
We now want to characterise the invariance of such subspaces $\mathcal Y$.
\begin{thm}\label{symmetries}
Consider a family of Hilbert spaces $H_i:=H_0:=L^2(X), V_i= V_0 \hookrightarrow H_0$, such that $X$ is a $\sigma$-finite measure space.
Fix a closed linear subspace $Y \subset \ell^2(I)$.
Assume $a=(a_{ij})_{i,j \in I}: V \times V \to \mathbb C$ to be a continuous and $H$-elliptic form.
Then for all Hilbert bases $(v_i)_{i\in I}$ of $H$ such that $(v_\ell)_{\ell \in I'}$ is a Hilbert basis of $Y$ and
$(v_k)_{k \in I''}$ is a Hilbert basis of $Y^\perp$ the following assertions are equivalent.
\begin{enumerate}[a)]
\item The subspace $\mathcal Y$ is invariant under the action of $\etasg$, i.e.,  $e^{ta}Y \subset Y$ for all $t\geq 0$.
\item The identity
\begin{equation}\label{autovk1}
\sum_{i,j \in I} \sum_{\ell \in I'} \sum_{k \in I''} (v_\ell)_j \overline{(v_k)_i} a_{ij}(\psi_\ell,\psi'_k)=0
\end{equation}
holds for all $\psi_\ell,\psi'_k \in V_0,$.
\end{enumerate}
\end{thm}
\begin{proof}
Observe that there exists a Hilbert basis $(v_\ell)_{\ell \in I}$ of the space $\ell^2(I)$ such that $(v_\ell)_{\ell \in I'}$ is a Hilbert basis of $Y$
and $(v_k)_{k \in I''}$ is a Hilbert basis of $Y^\perp$. In particular, $I' \cup I''= I$ and $I' \cap I''= \emptyset$. 
Denote $P_\mathcal Y $ the orthogonal projection of $H$ onto $\mathcal Y$.
According to Corollary \ref{ortho} the invariance property is equivalent to $P_\mathcal Y V \subset V$ and
$a(\psi,\psi')=0$ for all $\psi \in Y \cap V, \psi' \in Y^\perp \cap V$.

\smallskip
Denote $P_Y$ the orthogonal projection of $\ell^2(I)$ onto $Y$, and observe that the orthogonal projection $\mathcal P$ onto $\mathcal Y$ satisfies
$$
(P_\mathcal Y\psi)(x)= (P_Y(\psi(x))), \qquad \psi \in H, x \in X.
$$
So, it is in each component a convergent series of vectors in $V$. 
Since $V$ is a Hilbert space, then $P_\mathcal YV \subset V$ is satisfied.

\smallskip
We derive the algebraic relation~\eqref{autovk1}. 
Fix $\psi \in \mathcal Y \cap V$, $\psi' \in \mathcal Y^\perp \cap V$ and compute
\begin{eqnarray*}
a(\psi,\psi')&=&\sum_{i,j\in I}a_{ij}((\sum_{\ell \in I'}c^\psi_\ell v_\ell)_j, (\sum_{k \in I''}c^{\psi'}_k v_k)_i)\\
&=& \sum_{i,j \in I} \sum_{\ell \in I'} \sum_{k \in I''} (v_\ell)_j \overline{(v_k)_i} a_{ij}(c^\psi_\ell,c^{\psi'}_k)
\end{eqnarray*}
Since all coefficients may occur, this is equivalent to the relation in b).
\end{proof}

As an illustrative example we prove an irreducibility result.

\begin{exa}
Consider a family of Hilbert spaces $H_i:=H_0:=L^2(X), V_i= V_0 \hookrightarrow H_0$, such that $X$ is a $\sigma$-finite measure space.
Assume the form $a=(a_{ij})_{i,j \in I} : V \times V \to \mathbb C$ to be continuous and $H-$elliptic.
Then the subspace 
$$ \mathcal Y'=\{\psi \in H: \psi_i(x)=0, i \in I', x \in X \} $$ 
corresponds to the subspace 
$Y':=\ell^2(I')$ of $\ell^2(I).$
The projection on the subspace is given, of course, by 
$$(Pv)_i=\begin{cases}v_i & i \in I',\\ 0 & \mbox{otherwise.}\end{cases}$$
So, the subset $(e_\ell)_{\ell \in I'}$ of the canonical Hilbert basis $(e_\ell)_{\ell \in I}$ of $\ell^2(I)$
can be used as the Hilbert basis used in Theorem~\ref{symmetries}

Observe that in this case $v_{ij} =\delta_{ij}$ and so the sum~\eqref{autovk1} in Theorem~\ref{symmetries} can be computed by
\begin{eqnarray*}
\sum_{i,j \in I} \sum_{\ell \in I'} \sum_{k \in I''} (v_\ell)_j \overline{(v_k)_i} a_{ij}(\psi_\ell,\psi'_k)
&=&\sum_{i,j \in I} \sum_{\ell \in I'} \sum_{k \in I''} \delta_{\ell j} \delta_{ki} a_{ij}(\psi_\ell,\psi'_k)\\
&=&\sum_{\ell \in I'}\sum_{k \in I''}a_{k\ell}(\psi_\ell,\psi'_k)\\
%\sum_{j,\ell \in I'} \sum_{i,k \in I''}  a_{ij}(\psi_\ell,\psi'_k)\\
\end{eqnarray*}
where $I''=I \setminus I'$. 
Since all functions $\psi_\ell,\psi'_k$ may occur, then the sum only vanishes if
$a_{ij}=0$ for all $i \in I'', j \in I'$. 
So, the space $Y'=\bigoplus_{i \in I'} H_i $ is invariant if and only if $a_{ij}=0$ for all $i \in I'', j \in I'$. 

Conversely, the invariance of the subspace $\mathcal Y''=\bigoplus_{i \in I''} H_i$ induced by $\ell^2(I'')$ 
is equivalent to $a_{ij}=0$ for all $i \in I', j \in I''$.

Summing up, $Y'$ and $Y''$ are both invariant if and only if $(a,V)$ has block diagonal form, i.e., if and only if
$$
a_{ij}=a_{ji}=0, \quad  (i,j) \in I' \times I''.
$$
\end{exa}

%% file: sec_wave.tex
\section{A strongly damped wave equation}\label{sec:wave}
As an application of the theory developed so far, we investigate a strongly damped wave equation.
In Section~\ref{sec:histremI} references for the problem are given.
Let $\Omega$ be a bounded domain with boundary of class $C^\infty$, $\alpha \in \mathbb C$,
and consider the following second order problem.
\begin{equation}\label{problemwave}
\left\{
\begin{array}{rcll}
\ddot{u}(t,x)&=&\Delta(\alpha u+\dot{u})(t,x), &t\geq 0,\; x\in\Omega,\\
\frac{\partial u}{\partial \nu}(t,z)&=&\frac{\partial \dot{u}}{\partial \nu}(t,z)=0, &t\geq 0,\; z\in\partial \Omega,\\
u(0,x)&=&u_{0}(x), &x\in \Omega,\\
\dot{u}(0,x)&=&v_{0}(x), &x\in \Omega,
\end{array}
\right.
\end{equation}
The problem~\eqref{problemwave} can be reformulated as an abstract Cauchy problem of the form that we have presented in~\eqref{ACP}.
To see this, denote by $D(\Delta_N)$ the domain of the Laplacian with Neumann boundary conditions as in Example~\ref{multper}.
The operator governing the abstract Cauchy problem is given by
$$
 A=
\begin{pmatrix}
0 & I\\
\alpha \Delta & \Delta
\end{pmatrix},
\qquad 
D({ A})
=\{
\begin{pmatrix} \psi_1 \\ \psi_2 \end{pmatrix} \in H^1(\Omega)\times H^1(\Omega): 
\alpha \psi_1 + \psi_2 \in D(\Delta_N)\}.
$$
We want to investigate the problem by means of our theory. 
To this end we introduce a form matrix $a:V\times V\to\mathbb C$, where $V_1=V_2=H_1=H^1(\Omega)$, $H_2=L^2(\Omega).$ 
Set now
\begin{eqnarray*}
a_{11}(\psi,\psi')&:=&0, \\
a_{12}(\psi,\psi')&:=& -\int_\Omega \psi(x)\overline{\psi'(x)} dx - \int_\Omega \nabla \psi(x)\cdot \overline{\nabla \psi'(x)}dx,\\
a_{21}(\psi,\psi')&:=&\alpha \int_\Omega \nabla \psi(x)\cdot \overline{\nabla \psi'(x)}dx,\\
a_{22}(\psi,\psi')&:=&\int_\Omega \nabla \psi(x)\cdot \overline{\nabla \psi'(x)}dx.
\end{eqnarray*}
We prove that $(a,V)$ is the form associated with the operator $(A,D(A)).$
\begin{prop}\label{assocwave}
The operator $(A,D(A))$is associated with the form $(a,V)$.
\end{prop}
\begin{proof}
Denote $B$ the operator associated with $(a,V)$ and consider $\psi \in D(A), \varphi  \in V.$ 
By definition
\begin{eqnarray*}
a(\psi,\varphi)&=& a_{12}(\psi_2,\varphi_1)+a_{21}(\psi_1,\varphi_2)+a_{22}(\psi_2,\varphi_2) \\
&=&-(\psi_2,\varphi_1)_{H_1}+
\alpha \int_\Omega \nabla \psi_1(x)\cdot \overline{\nabla \varphi_2(x)}dx\\
&& +\int_\Omega \nabla \psi_2(x)\cdot \overline{\nabla \varphi_2(x)}dx.
\end{eqnarray*}
Applying on both second terms in right hand-side Green's first identity we obtain
\begin{eqnarray*}
a(\psi,\varphi)&=& -(\psi_2,\varphi_1)_{H_1} - \alpha \int_\Omega \Delta \psi_1(x) \overline{\varphi_2(x)} dx 
+ \alpha \int_{\partial \Omega} \frac{\partial \psi_1(x)}{\partial \nu} \overline{\varphi_2(x)} dx\\
&&-\int_{\Omega} \Delta\psi_2(x) \overline{\varphi_2(x)} dx + \int_{\partial \Omega} \frac{\partial \psi_2(x)}{\partial \nu} \overline{\varphi_2(x)} dx\\
&=&-(\psi_2,\varphi_1)_{H_1} 
+  \int_{\partial \Omega} (\alpha\frac{\partial \psi_1(x)}{\partial \nu} + \frac{\partial \psi_2(x)}{\partial \nu}) \overline{\varphi_2(x)} dx\\
&-&\int_{\Omega} \Delta\psi_2(x) \overline{\varphi_2(x)} dx - \alpha \int_\Omega \Delta \psi_1(x) \overline{\varphi_2(x)} dx \\
&=& - (A \psi, \varphi)_H,
\end{eqnarray*}
since $\alpha \psi_1 + \psi_2 \in D(\Delta_N)$. We have proved $A \subset B$. 

\smallskip
To see that the converse inclusion also holds, fix $\psi \in D(B)$.
Since $\psi \in D(B)$, there exists $f=(f_1,f_2) \in H$ such that
$
a(\psi,\varphi) = (f, \varphi)
$
for all $\varphi=(0,\varphi_2) \in V.$
This yelds
$$
a(\psi,\varphi) = \int_\Omega \nabla(\alpha  \psi_1(x)+ \psi_2(x))\overline{\varphi_2(x)} dx = \int_\Omega f_2(x) \overline{\varphi_2(x)} dx.
$$
By definition of the weak Laplacian, this means $\alpha  \psi_1 + \psi_2 \in D(\Delta_N)$ and $f_2=\Delta(\alpha  \psi_1 + \psi_2)$.
Choosing $\varphi=(\varphi_1,0) \in V$ yields $f_1 = \psi_2$ and so $B \subset A$. This completes the proof.
\end{proof}

The following properties hold.
\begin{itemize}
\item The form $a_{11}=0$, and so it is continuous on $V_1 \times V_1$ and is $H_1$-elliptic with constants $(\omega,\omega)$ for any $\omega \geq 0.$ 
\item The form $a_{22}$ is associated with the Laplace operator on $H_2$ with Neumann boundary conditions, thus it is elliptic and continuous. 
\item The sesquilinear mappings $a_{21}$ and $a_{12}$ are continuous on $H_1\times V_2$ and on $V_2\times H_1$, respectively. 
\end{itemize}
Combining these facts, we see that the problem (ACP) is well-posed by Proposition~\ref{perturbinter}.
\begin{prop}\label{wpwave}
The following assertions hold.
\begin{enumerate}[a)]
\item The sesquilinear form $(a,V)$ is continuous and $H$-elliptic for each $\alpha\in \mathbb C$.
\item The semigroup generated by $(A,D(A))$ is analytic of angle $\frac{\pi}{2}$.
\end{enumerate}
\end{prop}
\begin{proof}
Observe now that $\im a_{ii}(\psi)=0$ for all $\psi\in V_i$, $i=1,2$. 
Furthermore, there holds
 \begin{eqnarray*}
|\im \left(a_{12}(\psi,\psi')+a_{21}(\psi',\psi)\right)|&=&|\im(-\alpha \int_\Omega \nabla \psi(x)\cdot \overline{\nabla \psi'(x)}dx\\
&&+\int_\Omega \nabla \psi'(x)\cdot \overline{\nabla \psi(x)}dx)|\\
&&\leq (1+|\alpha|) \|\psi\|_{V_2}\|\psi\|_{H_1}.
\end{eqnarray*}
Thus, Proposition~\ref{cosine} applies and $A$ generates a cosine operator function, hence also an analytic semigroup of angle $\frac{\pi}{2}$. 
\end{proof}

We discuss now a particular invariance property of the problem~\eqref{problemwave}.

\begin{prop}\label{productwave}
Consider a closed product subspace $ Y= Y_1 \bigoplus Y_2$ of $H$ and denote by $P_1,P_2$ the projections of $H_1,H_2$ onto $Y_1,Y_2$ respectively. 
Then, $Y_1\bigoplus Y_2$ is invariant under $e^{ta}$ if and only if
\begin{itemize}
\item The subspace $Y_2$ is invariant under the semigroup $(e^{ta_{22}})_{t \geq 0}$ on $H_2$.
\item For all $\psi \in Y, \psi' \in Y^\perp$ 
$$
(\psi_2 \mid \psi'_1)_{H_1}=0, \quad \alpha (\nabla \psi_1 \mid \nabla \psi'_2 )_{H_2}=0.
$$
\end{itemize}
\end{prop}
\begin{proof}
The projection of $H$ onto $Y$ is given by
$$
P\psi = (P_1\bigoplus P_2)\begin{pmatrix} \psi_1 \\ \psi_2 \end{pmatrix} := \begin{pmatrix}P_1\psi_1\\P_2 \psi_2\end{pmatrix}
$$
By Proposition~\ref{ortho} we deduce that  $Y$ is invariant under the action of $(e^{ta})_{t\geq 0}$ if and only if
\begin{enumerate}[i)]
\item $P_{1} V_1 \subset V_1$ and $P_{2} V_2 \subset V_2$; 
\item for all $\psi \in (Y_1 \times Y_2) \cap V,\psi' \in (Y_1^\perp \times Y_2^\perp)\cap V$ 
$$
a(\psi,\psi')=( \psi_2,\psi'_1)_{H_1} + \alpha (\nabla \psi_1 \mid \nabla \psi'_2 )_{H_2} + (\nabla \psi_2 \mid \nabla \psi'_2)_{H_2}=0.
$$
\end{enumerate}
We prove the necessity of the conditions.

The first condition in i) is empty, since $V_1=H_1$. 
Fix $\psi_1=\psi'_1=0$. The second condition in i) and ii) imply that $e^{ta_{22}}Y_2 \subset Y_2$.

Choosing now $\psi_1=\psi_2'=0$, we obtain
$$
(\psi_2, \psi'_1)_{H_1}=0, \qquad \psi_2 \in Y_2, \psi'_1 \in Y_1^\perp,
$$
and choosing $\psi_2=\psi_1'=0$, we obtain
$$
\alpha (\nabla \psi_1 \mid \nabla \psi'_2 )_{H_2}=0, \qquad \psi_1 \in Y_1, \psi'_2 \in Y_2^\perp.
$$
So, we have shown that the conditions are necessary. The sufficiency is analogous.
\end{proof}
Using Proposition~\ref{productwave} we finally show that the space of radial functions is invariant under the action of the semigroup $\etasg$.
\begin{cor}
Assume $\Omega=B_R(0)$. 
Then the solution $u$ of~\eqref{problemwave} is radial provided that the initial data $(u_0,{v_0})$ are radial.
\end{cor}
\begin{proof}
Define $Y_1:=\{\psi \in L^2(\Omega): \psi \mbox{ is radial}\}$.  Then the claim is equivalent to the invariance of the subspace $Y=Y_1\cap H^1(\Omega) \times Y_1$.
The semigroup $(e^{ta_{22}})_{t \geq 0}$ is the heat semigroup with Neumann boundary conditions, and therefore it leaves invariant radial functions.
Observe that 
$$Y_1^\perp= \{\psi \in L^2 : \int_{\partial B_r(0)} \psi(x) d\sigma(x)=0 \}.$$

Let $\psi \in Y, \psi' \in Y^\perp$. 
In spherical coordinates $\psi_2(\theta, r)=f(r) $. We may thus compute
\begin{eqnarray*}
(\psi_2 \mid \psi'_1)_{H_1} &=& (\psi_2 \mid \psi'_1)_{L_2(\Omega)} + (\nabla \psi_2 \mid \nabla \psi'_1)_{L_2(\Omega)}\\
&=& \int_{\Omega} \nabla \psi_2(x) \cdot \overline{\nabla \psi'_1(x)} dx\\
&=& \int_{0}^R f'(r) \int_{\partial B_r(0)}  \nu(x) \cdot \overline{\nabla \psi'_1(x)} d\sigma(x)dr\\
&=& \int_{0}^R f'(r) \int_{\partial B_r(0)}  \frac{\partial \psi'_1(x)}{ \partial \nu} d\sigma(x)dr.
\end{eqnarray*}
Since $\int_{\partial B_r(0)} \psi'_1(x) d\sigma(x) =0$ for all $r$, in particular it is constant, 
and so 
$$
\frac{\partial }{\partial r}\int_{\partial B_r(0)} \psi'_1(x) d\sigma(x)=   \int_{B_r(0)} \frac{\partial \psi'_1(x)}{\partial \nu}=0.
$$
Plugging in this in the above equation, we obtain $(\psi_2 \mid \psi'_1)_{H_1}=0$. 
Exactly in the same way it is possible to see the second condition in Proposition~\ref{productwave}.
So, the claimed invariance holds.
\end{proof}

%% file: sec_dynamic.tex
\section{A heat equation with dynamic boundary conditions}\label{sec:dynamic}
As a second application of our theory, we investigate a heat equation with dynamic boundary conditions.
For references, see Section~\ref{sec:histremI}.

Consider a bounded open domain $\Omega \subset {\mathbb R}^n$ with $C^\infty$ boundary $\partial\Omega$ and set
$$V_1:=H^1(\Omega),\quad H_1:=L^2(\Omega),\quad V_2:=H^1(\partial\Omega),\quad H_2:=L^2(\partial \Omega).$$
Consider now the initial-boundary value problem
\begin{equation}\label{cenn2}
\left\{
\begin{array}{rcll}
 \dot{u}(t,x)&=&\Delta u(t,x), &t\geq 0, x\in\Omega,\\
 \dot{w}(t,z)&=& u(t,z)+\Delta_{\partial \Omega} w(t,z), &t\geq 0,z\in\partial\Omega,\\
 w(t,z)&=&\frac{\partial u}{\partial \nu}(t,z), &t\geq 0, z\in\partial \Omega,\\
 u(0,x)&=&f(x), &x\in\Omega,\\
 w(0,z)&=&h(z), &z\in\partial\Omega.
\end{array}
\right.
\end{equation}

Here, $\Delta_{\partial\Omega}$ denotes the Laplace--Beltrami operator, which is defined weakly as the operator associated with the form
$$a_{22}(\psi,\psi'):=\int_{\partial\Omega} \nabla \psi(x) \cdot \overline{\nabla \psi'(x)}d\sigma.$$
Moreover, define the forms
$$
\begin{array}{rcl}
a_{11}(\psi,\psi')&:=&\int_\Omega \nabla \psi(x)\cdot \overline{\nabla \psi'(x)}dx, \\
a_{12}(\psi,\psi')&:=&-\int_{\partial \Omega} \psi(x)\overline{\psi'(x)}_{|\partial \Omega} d\sigma,\\
a_{21}(\psi,\psi')&:=&-\int_{\partial \Omega} \psi_{|\partial \Omega}(x)\overline{\psi'(x)} d\sigma.
\end{array}$$
Define now
$$
{A}:=\begin{pmatrix}
\Delta & 0\\
\cdot_{|\partial \Omega} & \Delta_{\partial \Omega}
\end{pmatrix},\qquad 
D({A})=\left\{\begin{pmatrix} \psi_1 \\ \psi_2\end{pmatrix}\in H^2(\Omega)\times H^2(\partial \Omega): \frac{\partial \psi_1}{\partial \nu}=\psi_2\right\}.$$
Observe that the system~\eqref{cenn2} is governed by the operator $(A,D(A))$. 
Our aim is to prove well--posedness of the problem~\eqref{cenn2} on all $L^p$-spaces.

\begin{prop}\label{assocdyn}
The form $a=(a_{ij})_{i,j \in I}:V \times V \to \mathbb C$ is associated with the operator $(A,D(A))$.
\end{prop}
\begin{proof}
As in Proposition~\ref{assocwave}, we denote $(B,D(B))$ the operator associated with $(a,V)$.
We first prove $D(A) \subset D(B)$. Fix arbitrary $\psi \in D(A),\psi' \in V$.
Compute
\begin{eqnarray*}
-(A \psi' \mid \psi) &=& -\int_{\Omega} \Delta \psi_1(x) \overline{ \psi_1'(x)} dx\\
&&- \int_{\partial \Omega} \psi_{1_{|\partial \Omega}} (x) \overline{\psi'_2(x)} d\sigma(x)
- \int_{\partial \Omega} \Delta_{\partial \Omega} \psi_2(x) \overline{\psi'_2(x)} d\sigma(x)\\
&=& \int_{\Omega} \nabla\psi_1(x) \overline{\nabla \psi'_1(x)} dx - \int_{\partial \Omega} \frac{\partial \psi_1(x)}{\partial \nu} \overline{\psi'(x)} d\sigma(x)\\
&&- \int_{\partial \Omega} \psi_{1_{|\partial \Omega}} (x) \overline{\psi'_2(x)} d\sigma(x)
+ \int_{\partial \Omega} \nabla \psi_2(x) \overline{\nabla \psi_2'(x)} d\sigma(x)\\
&=& \int_{\Omega} \nabla\psi_1(x) \overline{\nabla \psi'_1(x)} dx - \int_{\partial \Omega} { \psi_2(x)} \overline{\psi'_{1|\partial \Omega}(x)} d\sigma(x)\\
&&- \int_{\partial \Omega} \psi_{1_{|\partial \Omega}} (x) \overline{\psi'_2(x)} d\sigma(x)
+ \int_{\partial \Omega} \nabla \psi_2(x) \overline{\nabla \psi_2'(x)} d\sigma(x)\\
&=& a(\psi,\psi').
\end{eqnarray*}

\smallskip
To see that $D(B) \subset D(A)$ also holds, fix $\psi \in D(B)$. 
Then there exists $f \in H$ such that for all $\psi' \in V$ holds
$$
(f, \psi' )_H = a(\psi, \psi').
$$
Choose $\psi'=(\psi'_1,0)^\top$ or $\psi'=(0, \psi'_2)^\top$. 

In the first case, choose $\psi'_1 \in C^\infty_c(\Omega)$. So,
$$
\int_{\Omega} f_1(x) \overline{\psi'_1(x)} dx = 
\int_{\Omega} \nabla\psi_1(x) \overline{\nabla \psi'_1(x)} dx 
$$
and, since $\psi'_1$ is arbitrary, this means $\psi \in D(\Delta), $ and $f_1=-\Delta \psi_1$ weakly.
Choose now an arbitrary $\psi'_1$. Now
$$
\int_{\Omega} f_1(x) \overline{\psi'_1(x)} dx = 
\int_{\Omega} \nabla\psi_1(x) \overline{\nabla \psi'_1(x)} dx 
- \int_{\partial \Omega} \psi_2(x)\overline{\psi'_1(x)} d\sigma (x).
$$
Since $f_1= \Delta \psi_1$ and the Green's first identity applies to the right hand-side, one obtains
\begin{eqnarray*}
- \int_{\Omega} \Delta \psi_1(x) \overline{\psi'_1(x)} dx &=& 
- \int_{\Omega} \Delta\psi_1(x) \overline{\psi'_1(x)} \\
&&+ \int_{\partial \Omega} \frac{ \partial \psi_1(x)}{\partial \nu} \overline{\psi'_1(x)} d\sigma(x)
- \int_{\partial \Omega} \psi_2(x)\overline{\psi'_1(x)} d\sigma (x).
\end{eqnarray*}
So, $\frac{\partial \psi_1(x)}{\partial \nu}(x) = \psi_2(x)$ weakly.

For the second choice compute
$$
\int_{\partial \Omega} f_2(x) \overline{\psi'_2(x)} dx = 
- \int_{\partial \Omega} \psi_{1_{|\partial \Omega}} (x) \overline{\psi'_2(x)} d\sigma(x)
+ \int_{\partial \Omega} \nabla \psi_2(x) \overline{\nabla \psi_2'(x)} d\sigma(x).
$$
This means $- f_2= \psi_{1_{|\partial \Omega}} + \Delta \psi_2$ weakly.
\end{proof}
We now address well--posedness properties of the system~\eqref{cenn2}.
\begin{prop}
 Following assertions hold.
\begin{enumerate}[a)]
\item The form $(a,V)$ is continuous and elliptic. Thus, it generates an analytic semigroup $\etasg$.
\item The semigroup $\etasg$ is analytic of angle $\frac{\pi}{2}$.
\item The semigroup $\etasg$ is real, positive.
\item Denote $(e^{ta_0})_{t \geq 0}$, where $a_0$ is associated with the uncoupled system, i.e., $a_{0_{12}}=a_{0_{21}}=0$. 
Then $\etasg$ dominates $(e^{ta_0})_{t \geq 0}$.
\item The semigroup $\etasg$ does not leave invariant the unitary ball of $L^\infty$.
\end{enumerate} 
\end{prop}
\begin{proof}
a) Observe that $a_{11}$ (resp. $a_{22}$) are continuous and $H_1$- (resp. $H_2$-) elliptic. 
Moreover, due to boundedness from $H^1(\Omega)$ to $L^2(\partial\Omega)$ of the trace operator, 
the forms $a_{12}$ and $a_{21}$ are bounded on $H_2\times V_1$ and on $H_1\times V_2$, respectively. 
Accordingly, $(a,V)$ is continuous and by Proposition~\ref{perturbinter} also $H_1\times H_2$-elliptic.

\smallskip
b) In order to estimate the angle of analyticity, we want to apply Proposition~\ref{cosine}. 
To this end, observe that $\im a_{ii}(\psi)=0$ for all $\psi\in V_i$, $i=1,2$. 
Moreover
\begin{eqnarray*}
|\im \left(a_{12}(\psi,\psi')+a_{21}(\psi',\psi)\right)|&=&|\im(\int_{\partial \Omega} \psi(x)\overline{\psi'_{|\partial \Omega}(x)} d\sigma(x) \\
&&+\int_{\partial \Omega} \psi'_{|\partial \Omega}(x)\overline{\psi(x)} d\sigma(x))|\\
&=&|\im (\int_{\partial \Omega} \psi(x)\overline{\psi'(x)}_{|\partial \Omega} d\sigma(x)\\
&&+\overline{\int_{\partial \Omega} \psi(x)\overline{\psi'_{|\partial \Omega}(x)} d\sigma(x)})|\\
&=&0.
\end{eqnarray*}
and so Proposition~\ref{cosine} applies. So, $A$ generates a cosine family and, thus, the analyticity angle is the maximal one. 

\smallskip
c) Theorem~\ref{positivity} promptly yields that the semigroup is real. 
To see that it is positive, 
observe that $a_{11}$ is associated with the Laplace operator with Neumann boundary conditions and $a_{22}$ with the Laplace--Beltrami operator on $\partial\Omega$. 
Therefore they generate positive semigroups and the first condition of Theorem~\ref{positivity} is satisfied. 
The second condition is also clear, since $\psi_{|\partial \Omega}$ is positive whenever $\psi$ is positive. 

\smallskip
d) It is a direct consequence of Proposition~\ref{domination}.

\smallskip
e) It also follows from Proposition~\ref{positivity}
that $(e^{ta})_{t\geq 0}$ is not $L^\infty(\Omega)\times L^\infty(\partial\Omega)$-contractive, 
since for non-constant $\psi\in H^1(\partial\Omega)$ such that $|\psi|\leq 1$ and for $\psi'\in H^1(\Omega)$ with $\psi'_{|\partial\Omega}=1+\psi$ the estimate
$a_{12}(\psi,\psi')=-\int_{\partial \Omega} |\nabla \psi|^2 d\sigma<0$ holds.
\end{proof}

The usual method to prove well--posedness in $L^p$-spaces is to prove $L^\infty$-contractivity of a semigroup and its adjoint. 
In this case, an additional perturbation argument is needed.

\begin{thm}
The semigroup $\etasg$ extrapolates to a consistent family of semigroups on all $L^p$-spaces.
\end{thm}
\begin{proof}
The strategy is to write the generator as a relatively compact perturbation of a the generator of an ultracontractive semigroup and
then obtain consistent semigroups on all $L^p$-spaces by a perturbation theorem.

If $p >2$, we define operators $\tilde{A},B$ by
$$
\tilde{A}:=\begin{pmatrix}
\Delta-C^* & 0\\
\cdot_{|\partial \Omega} & \Delta_{\partial \Omega}-Id
\end{pmatrix},
\qquad 
{ B}:=
\begin{pmatrix}
C^* & 0\\
0 &Id 
\end{pmatrix}.$$
Here, $C^*$ is the adjoint of the linear operator from $H^1(\Omega)$ to $L^2(\Omega)$ defined by
$$C\psi(x):=\nabla \psi (x) \cdot \overline{\nabla D_N1(x)},\qquad \psi\in H^1(\Omega), x\in\Omega,$$
where $D_N1$ denotes the unique (modulo constants) solution $u$ of the inhomogeneous Neumann problem
$$\left\{
\begin{array}{rcll}
\Delta u(x)&=& 0,\qquad &x\in\Omega,\\
\frac{\partial u}{\partial \nu}(z)&=&1,&z\in\partial \Omega.
\end{array}
\right.$$
The operator $\tilde{ A}$ is associated with the matrix form $\tilde{a}$ whose diagonal entries are given by
\begin{eqnarray*}
\tilde{a}_{11}(\psi,\psi')&:=&\int_\Omega \nabla \psi(x)\cdot \overline{\nabla \psi'(x)}dx
+\int_\Omega \psi \overline{(\nabla \psi'\cdot \nabla D_N 1)} dx,\\ 
\tilde{a}_{22}(\psi,\psi')&:=&\int_{\partial\Omega} \nabla \psi(x) \cdot \overline{\nabla \psi'(x)}d\sigma 
+\int_{\partial\Omega}\psi(x) \overline{\psi'(x)}d\sigma.\\
\end{eqnarray*}
The off-diagonal entries are given by
$$
\tilde{a}_{12}:=a_{12},\qquad \tilde{a}_{21}:= a_{21}.
$$
One sees that the perturbation $\tilde{a}_{11}-a_{11}$ is bounded on $H_1\times V_1$, 
thus by Proposition~\ref{perturbinter} 
$\tilde{a}$ is associated with a semigroup $(e^{t\tilde{a}})_{t\geq 0}$ on $H_1\times H_2$. 

Compute now for any $\psi'\in H^1(\partial\Omega)$ such that $|\psi'|\leq 1$ and all $\psi\in H^1(\Omega)$
\begin{eqnarray*}
|\tilde{a}_{12}(\psi',(|\psi|-1)^+ {\rm sign}{\psi})|
&\leq &\int_{\partial\Omega} |\psi'(x)| (|\psi(x)|-1)^+ d\sigma(x)\\
&\leq & \int_{\partial\Omega} (|\psi(x)|-1)^+ d\sigma(x)\\
&=&\int_{\partial\Omega} \frac{\partial D_N 1}{\partial \nu}(|\psi(x)|-1)^+ d\sigma(x)\\
&= & \int_{\Omega} \nabla(|\psi(x)|-1)^+\cdot \nabla D_N 1 dx.
\end{eqnarray*}
Splitting now the integral in two parts, the last expression can be written as
\begin{eqnarray*}
&&\int_\Omega \Big[ (1\wedge |\psi(x)|) \left(\nabla(|\psi(x)|-1)^+ \cdot \nabla D_N 1\right) {\mathbb 1}_{\{|\psi(x)|\geq 1\}} \Big] dx\\
&+&\int_\Omega \Big[ (1\wedge |\psi(x)|) \left(\nabla(|\psi(x)|-1)^+ \cdot \nabla D_N 1\right) {\mathbb 1}_{\{|\psi(x)|\leq 1\}} \Big]dx.
\end{eqnarray*}
The latter is
$$
\re \tilde{a}_{11} ((1\wedge |\psi|){\rm sign}\psi,(|\psi|-1)^+{\rm sign}\psi).
$$
since $\nabla(|\psi|-1)^+=0$ a.e. on $\{x\in \Omega : |\psi(x)|\leq 1\}$.

Likewise, for all $\psi\in H^1(\Omega)$ such that $|\psi|\leq 1$ and all $\psi'\in H^1(\partial\Omega)$ compute
\begin{eqnarray*}
|\tilde{a}_{21}(\psi,(|\psi'|-1)^+ {\rm sign}{\psi'})|
&\leq &\int_{\partial\Omega} |\psi(x)| (|\psi'(x)|-1)^+ d\sigma(x)\\
&\leq & \int_{\partial\Omega} (|\psi'(x)|-1)^+ d\sigma(x)\\
&=& \int_{\partial\Omega} (1\wedge |\psi'(x)|) (|\psi'(x)|-1)^+ {\mathbb 1}_{\{|\psi'(x)|\geq 1\}}dx\\
&& +\int_{\partial\Omega}(1\wedge |\psi'(x)|) (|\psi'(x)|-1)^+ {\mathbb 1}_{\{|\psi'(x)|\leq 1\}} dx
\end{eqnarray*}
since $|\psi_{|\partial\Omega}|\leq 1$.

The last term corresponds to
$$
\re \tilde{a}_{22} ((1\wedge |\psi'(x)|){\rm sign}\psi'(x),(|\psi'(x)|-1)^+{\rm sign}\psi'(x)).
$$
Thus, Theorem~\ref{positivity} applies, the semigroup $(e^{t\tilde{a}})_{t\geq 0}$ is $L^\infty$-contractive,
and we conclude that  extrapolates to a consistent family of semigroups on $L^p(\Omega)\times L^p(\partial\Omega)$, $p\geq  2$.
The generator of the semigroup in $L^p(\Omega)\times L^p(\partial \Omega)$ is the part of $\tilde{ A}$ in $L^p(\Omega)\times L^p(\partial \Omega)$.

The part of $ B$ is a compact operator from $W^{2,p}(\Omega)\times W^{2,p}(\partial \Omega)$ to $L^p(\Omega)\times L^p(\partial\Omega)$ 
for all $p=[1,\infty)$, and so we conclude that also (the part of) ${ A}=\tilde{ A}+{ B}$ generates a semigroup on $L^p(\Omega)\times L^p(\partial\Omega)$, $p\geq 2$.

\smallskip
The case of $p \in (1,2)$ can be carried out in the following manner. 
Introduce an operator $\tilde{ A}$ by replacing $C^*$ by $C$ and the corresponding form $\tilde{a}$. 
By the same arguments as in the first case, the semigroup associated with the adjoint $\tilde{a}^*$ is $L^\infty$-contractive.
By duality and perturbation arguments we conclude as above that $A$ generates a semigroups on $L^p(\Omega)\times L^p(\partial\Omega)$ also for $p\in [1,2]$.
\end{proof}

\begin{rem}
The standard way to deduce $L^p$-well--posedness of the semigroup $(e^{t\tilde{a}})_{t\geq 0}$ would be via ultracontractivity.
Observe now that 
$|\tilde{a}_{12}((|\psi'|-1)^+ {\rm sign}{\psi'},\psi)|\leq \re \tilde{a}_{22} ((|\psi'|-1)^+{\rm sign}\psi',(1\wedge |\psi'|){\rm sign}\psi')$ 
does not hold for all $(\psi,\psi')\in H^1(\Omega)\times H^1(\partial\Omega)$ such that $|\psi|\leq 1$. 
Thus condition (iii) in Theorem~\ref{ultra}.(1) is not satisfied and it is not possible to deduce ultracontractivity of $(e^{t\tilde{a}})_{t\geq 0}$ 
from  Theorem~\ref{ultra} 
and the Sobolev embeddings $H^1(\Omega)\hookrightarrow L^\frac{2n}{n-2}(\Omega)$, $H^1(\partial\Omega)\hookrightarrow L^\frac{2n-2}{n-3}(\partial\Omega)$. 
\end{rem}

%% file: sec_nondiag.tex
\section{Some remarks on non-diagonal domains}\label{sec:nondiag}

\begin{ass}
During the rest of this section we always assume the following.
\begin{itemize}
\item The subspace $Y$ is closed in the Hilbert space $V$.
\item The Hilbert spaces $Y$ and $V$ are densely continuously embedded into the Hilbert space $H$.
\item The linear operator $P_Y$ denotes the orthogonal projection of $V$ onto $Y$.
\item We denote by $a_Y$ the restriction of the form $a$ to the subspace $Y$, $A_Y$ the operator associated with the form $a_Y$ and $A_{|Y}$ the restriction of the operator $A$ to the subspace $Y \cap D(A)$.
\end{itemize}
\end{ass}
Under these assumptions the following easy result holds.
\begin{lemma}
If the sesquilinear form $a$ is continuous, then $a_Y$ is continuous.
The same holds for accretivity, coercivity and ellipticity.
\end{lemma}
In the general case, it is difficult to identify the domain of the operator associated with the form $a_Y$. 
Although the operator $A_Y$ is contained in the operator $A_{|Y}$, in the sense of operator theoretic inclusion, 
it is difficult to say whether the operator $A_Y$ coincides with the part of the operator $A$ in $Y$. 
In Chapter \ref{network} we will investigate the case in which the subspace $Y$ reflects the topological structure of a graph: for all those forms the domain of $A_Y$ does not, in fact, coincide with the part of the operator $A$ in the subspace $Y$. 

We turn our attention to characterisations of coercivity properties of the restricted form. 
This ensures that the operator $A_Y$ generates an analytic semigroups, even if it is not possible to explicitly identify the domain $D(A)$.
The following result is less useful than it might appear at a first glance.
\begin{lemma}\label{expstab}
Assume that the subspace $Y$ is invariant under the action of $\etasg$, where $a$ is a continuous, elliptic form. Then the form $a_Y$ is coercive with constant $\alpha$ if and only if 
$$
\re a(\psi,P_Y \psi) \geq \alpha \|P_Y\psi\|^2_V, \qquad \mbox{for all } \psi \in V. 
$$
\end{lemma}
\begin{proof}
By definition, the form is coercive if and only if $\re a(\psi) \geq \alpha \|\psi\|^2_V$ for all $\psi \in Y$. Since the projection $P_Y$ is, by definition, surjective onto $Y$, the last is equivalent to 
$$
\re a(P_Y \psi) \geq \alpha \|P_Y \psi\|^2_V, \qquad \mbox{for all } \psi \in V.
$$
Using Corollary \ref{ortho}, compute
\begin{eqnarray*}
\re a(P_Y \psi)&=& \re a(P_Y \psi, P_Y \psi)\\
&=& \re a(\psi, P_Y \psi) - \re a(\psi - P_Y \psi , P_Y \psi)\\
&=& \re a(\psi, P_Y \psi),
\end{eqnarray*}
and the claim holds.
\end{proof}
In fact, we assumed in the Lemma that the form $(a,V)$ is $H$-elliptic, and we derived a condition for the coercivity on $Y$.
In other words, Lemma~\ref{expstab} only yields a sufficient condition for the exponential stability along an invariant subspace. 
The very difficult challenge, however, is to find non-trivial conditions on a non-elliptic form $(a,V)$
which guarantee the ellipticity of $a_Y$. 

\begin{rem}
Assume the form $(a,V)$ to be continuous, accretive and symmetric. 
It is tempting to use a curve integral in the Hilbert space $H$ to characterise the coercivity of $a_Y$. 
Recall that since $(a,V)$ is continuous and densely defined, there exists a densely defined operator $A: D(A) \to H$ associated with $(a,V)$.  
For an arbitrary $\psi \in D(A)$, denote $\gamma: [0,1] \to H$ the curve defined by $\gamma : t \mapsto (1-t)\psi + tP_Y\psi$. 
Formally, we compute
\begin{eqnarray*}
 \frac{1}{2} a(P_Y \psi) 
&=& \int_0^1 (\partial a(\gamma(t)), \frac{d \gamma(t)}{dt}) dt \\
&=&  \int_0^1 (A \gamma(t), P_Y \psi - \psi) dt \\
&=&  \int_0^1 a((1-t) \psi + tP_Y \psi, P_Y \psi - \psi) dt \\
&=&  a(P_Y \psi, P_Y \psi - \psi),
\end{eqnarray*}
where $\partial a(\psi)$ has to be understood in the sense of monotone operators, see~\cite{Bre73}.
Thus, for all $\psi \in D(A)$ the above computation, if rigorously justified, yields an estimate for $a(P_Y \psi)$, 
for all $\psi \in D(A)$. 
\end{rem}

%% file: sec_histremI.tex
\section{Discussion and remarks}\label{sec:histremI}
In contrast to matrices of operators, matrices of sesquilinear forms have not received until now a great attention in the literature.
However, it seems to us that even in the simple case of diagonal domains, form matrices present some fundamental advantages, 
if compared to operator matrices. We now want to comment two papers that already contained the main ideas of this work.

In 1989, Rainer Nagel pointed out (see~\cite{Nag89}) that \emph{``many linear evolution equations $[...]$ can be written formally as a first order Cauchy problem $[...]$ with values in a product of different Banach spaces.''} 
There the generator properties of an operator matrix
$$
\mathcal A: = \begin{pmatrix} A& B \\ C& D\end{pmatrix}
$$ 
with domain
$$
D(\mathcal A):=D(A) \times D(D)
$$
are discussed, under the conditions that $C$ is relatively $A$-bounded and $B$ is relatively $D$-bounded. 

Observe that the in the case of operator matrices, diagonal domains lead to uncomfortable phenomena. 
Consider an arbitrary Banach space $X$ and an arbitrary closed operator $(A,D(A))$ on $X$, the operator matrix
$$
\mathcal A:=\begin{pmatrix} A & A\\A&A\end{pmatrix}
$$
is closed for the domain
$$
D(\mathcal A):= \left\{\begin{pmatrix} \psi_1 \\ \psi_2 \end{pmatrix} \in X \times X: \psi_1 + \psi_2 \in D(A) \right\}.
$$
However, the operator $\mathcal A$ on the diagonal domain $D(A) \times D(A)$ is not closed.
This observation led Rainer Nagel to claim that the diagonality of the domain \emph{``seem to be artificial''}. 
In fact, in subsequent works of Klaus-Jochen Engel (see e.g. \cite{Eng97}), the case of non-diagonal domain was investigated. 

One can ask the question whether diagonal domains are artificial also for forms.
We do not intend to answer this somehow philosophical question but do observe that the same example does not work in the case of forms.
To see this, assume that $X$ is a Hilbert space and that the operator $(A,D(A))$ on $X$ is associated with a continuous form $(a,Y)$ on a dense $Y$. 
Set $H_1=H_2=X$, $V_1=V_2=Y$ and $a_{ij}=a$ for $i,j \in \{1,2\}$. Denote $( B, D( B))$ the operator associated with $(a,V)$ and
fix a vector $\psi=(\psi_i)_{i=1,2} \in D( B)$. By definition, there exists $f=(f_i)_{i=1,2} \in H$ such that 
$$
a(\psi,\psi')=(f \mid \psi')_H, \qquad \psi' \in H.
$$
In particular, for $\psi'_2=0$ we obtain $a(\psi_1 +\psi_2,\psi'_1)= (f_1 \mid \psi'_1)$, i.e. $\psi_1 + \psi_2 \in D(A)$. 
The proof of the converse inclusion is analogous.
So the closed operator $(\mathcal A,D(\mathcal A))$ has non-diagonal domain, but it is associated with a form that has diagonal domain. 

A further important work to Felix Ali Mehmeti and Serge Nicaise, see~\cite{AliNic93}. 
The work is concerned with the analysis of \emph{``evolution phenomena on a (possibly infinite) family of domains''}. 
As they are interested in parabolic problems, they work in the framework of bilinear forms (and of maximal monotone operators for the nonlinear case). 
In fact, they introduce the form domain $V$ as in Section \ref{sec:findimarg}, i.e., they consider the case of diagonal domains. 
They also develop a theory for diagonal forms, i.e.,
$a_{ij}=0$ for all $i \neq j$ and they consider in the Example 1.2.3 a particular strongly coupled form, observing that it is elliptic. 
However, the major goal of that work was to prove well-posedness and regularity properties of non-linear systems; 
further, the invariance criterion of El-Maati Ouhabaz was not yet known.

In fact, our work can be seen as an extension of~\cite{AliNic93} to the case of non-diagonal sesquilinear forms, 
in which we are also able to investigate qualitative properties by Theorem~\ref{ouh}.

\subsubsection{Section \ref{sec:findimarg}}
The results in this section are a generalisation of results contained in~\cite{CarMug07}.
In Example~\ref{regularoperators} we exhibited an elementary example in which the boundedness constant is not optimal.
This phenomenon has been deeply investigated in the literature.
Assume that the form is of the type $a(\psi,\psi')=(C \nabla \psi \mid \nabla \psi')$, where $C$ is an operator on $\ell^2(I)$.
Now, the estimate in Lemma~\ref{boundedform} in this case means $|a(\psi,\psi')| \leq \||C|\|_{\mathcal L(\ell^2(I))}.$
So, this argument only makes sense if $C$ is a regular operator.
Consider a matrix $C$ which is a bounded, non regular operator on $\ell^2(\mathbb N)$, then the criterion in Lemma~\ref{boundedform} fails, 
but the form $(a,V)$ is continuous as a consequence of Proposition~\ref{contsystems}.
Observe that a matrix which is an operator on $\ell^2(I)$ but whose modulus is unbounded has been constructed by Wolgang Arendt and J\"urgen Voigt in~\cite{AreVoi91}. 
There, they also show that the space of regular operators is not even dense in $\mathcal L(\ell^2)$.

\subsubsection{Section \ref{sec:systemsI}}
The results in this section are unpublished.
Systems of PDE's are a major topic in the literature 
and we refer to~\cite{Ama84} for a systematic presentation of the theory in the finite case.
There, also the case of non-linear, non-autonomous equations is discussed.

It is worth mentioning that \emph{weighted positivity} has been discussed in the literature, 
see e.g.~\cite{LuoMaz07}, and references therein.

In~\cite{LuoMaz07} is considered the case of a strongly coupled operator $(A,D(A))$ in non-divergence form.
The operator $(A,D(A))$ is said to be weighted positive with weight $M$, if
$$
\int_{\Omega} (A \psi \mid M \psi) dx \geq 0,
$$
with a positive definite matrix $M$. 
This, in fact, equivalent to the accretivity of the form in the rescaled Hilbert space 
$L^2(\Omega, \ell^2(I); M dx)$, 
and so the theory developed for systems can be applied in this case. See also Example~\ref{multper}.

We also mention that it is possible to generalise the theory developed in Section~\ref{sec:histremI} substituting the operator $M^d_C$ by a more general operator of the form $\id \otimes (M_{C_1},\ldots,M_{C_d})$.

\subsubsection{Section \ref{sec:gencoercivity}}
The results in this section are unpublished. The idea of the proof of Proposition~\ref{propcoerc2x2} is due to Wolfgang Arendt.
In Section~\ref{sec:systemsI} we have investigated coercivity properties for matrices in the case that the single forms have a particular expression.
Proposition~\ref{propcoerc2x2} is an example of a characterisation exploiting a particular structure of the form matrix.
An open question is whether it is possible to translate other arguments holding for scalar-valued matrices in the context of forms, 
e.g. symmetric matrices or Toeplitz matrices.

\subsubsection{Section \ref{sec:operator}}
Proposition~\ref{contembed}.b) is due to Felix Ali-Mehmeti and Serge Nicaise in \cite{AliNic93}. We slightly simplified the proof.
Proposition~\ref{operdomain} is already contained in~\cite{CarMug07}.
Example~\ref{multper} uses Ulm's recipe for the investigation of multiplicative perturbations of the Laplacian by means of sesquilinear forms, see~\cite{Are06}. 
See also the comments to Section~\ref{sec:systemsI}.

\subsubsection{Section \ref{sec:evoleq}}
The results in this section mostly come from~\cite{CarMug07}, with some minor generalisations.
We chose not to give the proofs of any general results concerning the theory of sesquilinear forms, and we refer to the Appendix. 

However,  the reader interested in the general theory of evolution equations governed by operators associated with an energy form should consult the introductory manuscript \cite{Are06} and the fundamental monograph \cite{Ouh04}. 

Invariance of sets and subspaces have already been investigated in the literature, see e.g.~\cite{Ama78}.
Most works focus on systems of PDE's, that are a special case of our theory.

Theorem~\ref{positivity} deals with positivity of form matrices, when the single mappings are in general defined on the product of different $H^1(X_i)$-spaces.
If $X_i=X$ for all $i$ and $a_{ij}$ are bounded mappings from $H_j \times H_i$, then we obtain the case of \emph{weakly coupled} parabolic systems.
Such systems have been extensively investigated in the literature.
For the case $C(\overline{\Omega})$ see~\cite{Swe92} and \cite{MitSwe95}, and \cite{BerGazMit06} for the case of fourth order uncoupled operators in $L^2$.

Observe that in Theorem~\ref{positivity} positivity of the off-diagonal forms is required in order to achieve positivity.
In the case of systems this already implies that the coupling is of order 0, 
although it can also works with unbounded coefficients, see~\cite{AreTho05}. In fact, weak coupling is a necessary condition.

\subsubsection{Section \ref{sec:symmetriesmatrix}}
The results in this section are a generalisation of results contained in~\cite{CarMug07}.
We have always assumed that the set $I$ is a countable set. 
Indeed, this condition can in principle be weakened by requiring that for all $i \in I$, $a_{ij}=0$ vanishes for all but countable $j$ 
and that for all $j \in I$, $a_{ij}=0$ vanishes for all but countable $i$. 
Then, due to a graph theoretic argument (see Section~\ref{sec:histremII}), the uncountable matrix $a=(a_{ij})_{i,j\in I}$
can be decomposed in the uncountable direct sum countable matrices $a_s=(a^s_{ij})_{i,j=1}$, each of them being defined on a separable ideal $I_s$ of $\ell^2(I)$. 
Here $s$ is an index in an uncountable set $S$.
Then, the parabolic equations leaves each of these ideals invariants. 
As a consequence, if the non-countable matrix $a=(a_{ij})_{i,j \in I}$ has only countable many non vanishing terms in each of the rows and columns, 
then all properties of $\etasg$ can be studied on the level of the single separable ideals.

\subsubsection{Section \ref{sec:wave}}
The ideas contained in this section were introduced in~\cite{CarMug07}.
The problem is discussed in less generality and with the same techniques in~\cite[\S~XVIII.6]{DauLio88}. 
The case of a strictly positive damping term $\alpha$ was investigated in~\cite[Thm.~1.1]{CheTri89}.
See also~\cite{Mug08} for a thorough treatment of the problem with techniques of sesquilinear forms. 
There also well-posedness issues of related non-linear problems are discussed.

\subsubsection{Section \ref{sec:dynamic}}
The results in this section are taken from~\cite{CarMug07}. In particular, the perturbation argument was used there.
The definition of the Laplace-Beltrami operator is, in general, quite difficult and involves concepts coming from the differential geometry.
However, in the case of hypersurfaces it can be elementary defined, see~\cite{GilTru77}.
The same system has been investigated in~\cite{CasEngNag03} and~\cite{Mug06b}.
There, well-posedness and exponential stability of~\eqref{cenn2} has been proved by different methods.

\subsubsection{Section \ref{sec:nondiag}}
The case of non-diagonal domains is more difficult and we did not success in developing a satisfactory theory for the general case.
In fact, also for operators on product spaces it is quite difficult to develop an unified theory for operator matrices with non-diagonal domain. 
To be more precise, only in the case that the operator allows a matrix representation, it is indeed possible to ask formulate the question 
whether the domain is diagonal or not. An important class of operators in this context are the \emph{one-sided coupled} operator matrices, see \cite{Eng98}. Another approach usual in the context of operator theory is that of perturbation of domains, see \cite{Gre87}, and also \cite{GreKuh89}.
However, no general results for sesquilinear forms are known to us.

%% file: sec_infintro.tex
\section{Wolfgang's first network tutorial}\label{sec:infintro}
In this section, we want to introduce the basic tools necessary to define and study diffusion equations in network-shaped structures.
We use the word \emph{graph} and the word \emph{network} with two different meanings. A graph is a combinatorial object, whereas
a network is the measure theoretic object on which we define functions.

We start discussing the relations between these two objects in Figure \ref{fig:example}.

\begin{center}
\begin{pgfpicture}{0cm}{0cm}{5cm}{5cm}\label{fig:example}
\pgfsetendarrow{\pgfarrowtriangle{3pt}}
\pgfxycurve(1,1)(0,2)(0,3)(1,4)
\pgfxycurve(1,4)(2,3)(2,2)(1,1)
\pgfline{\pgfxy(1,4)}{\pgfxy(4,4)}
\pgfline{\pgfxy(4,4)}{\pgfxy(4,1)}
\pgfline{\pgfxy(1,1)}{\pgfxy(4,1)}
\pgfputat{\pgfxy(1,4.5)}{\pgfbox[center,center]{$\mv_1$}}
\pgfputat{\pgfxy(1,0.5)}{\pgfbox[center,center]{$\mv_2$}}
\pgfputat{\pgfxy(4,4.5)}{\pgfbox[center,center]{$\mv_3$}}
\pgfputat{\pgfxy(4,0.5)}{\pgfbox[center,center]{$\mv_4$}}
\pgfputat{\pgfxy(2.5,-0.5)}{\pgfbox[center,center]{Figure \ref{fig:example}: A simple example}}
\pgfputat{\pgfxy(2.5,4.5)}{\pgfbox[center,center]{$\me_3$}}
\pgfputat{\pgfxy(4.5,2.5)}{\pgfbox[center,center]{$\me_4$}}
\pgfputat{\pgfxy(2.5,0.5)}{\pgfbox[center,center]{$\me_5$}}
\pgfputat{\pgfxy(-0.2,2.5)}{\pgfbox[center,center]{$\me_1$}}
\pgfputat{\pgfxy(2.2,2.5)}{\pgfbox[center,center]{$\me_2$}}
\end{pgfpicture}
\end{center}
\vspace*{0.5cm}
One sees that we have fixed some labels on the graph, $\me_i$, $\mv_k$. 
This is always the case: we only consider labelled graphs, even if it were possible to develop a theory for diffusion equations on networks avoiding these combinatorial issues, see Section~\ref{sec:histremII} for these considerations. 
We assume that labels of the type $\mv_k$ always refer to \emph{vertexes} of the graph, whereas labels of the type $\me_i$ always refer to \emph{edges} of the graph. 
In fact in order to avoid confusion while juggling with a large number of indexes, we always use the indexes $k,\ell,n$ when we are handling with vertexes and the indexes $i,j,m$ if handling with edges. 

Since we are interested in second order problems, the orientation of the network could seem redundant or artificial. 
However, since we want to consider functions taking values at points of edges, we need to fix a system of coordinates on the network.
In this formalism each edge is assumed to have length one (possibly rescaling the coefficients of the equation), 
i.e.,  each edge is identified with a copy of the interval $[0,1]$. 
When we represent the network by a picture like in Figure~\ref{fig:example}, 
the arrows are oriented in such a way that the terminal endpoint of an arrow represents the origin of the coordinate system, 
i.e.,  the point 0, and the initial endpoint of an arrow represents the end of the coordinate system, i.e.,  the point 1. 

So far, we gave an intuitive definition of a graph. 
We now want to identify a \emph{graph} with a combinatorial (or algebraic object) and a \emph{network} with a measure-theoretic one.
Again, we start discussing the graph in Figure \ref{fig:example}. 
In this graph, there is an edge from $\mv_1$ coming into the vertex $\mv_3$. 
So, identifying the vertexes with their order numbers, we say that the couple $(3,1)$ belongs to the graph. 
Following this idea, we define a function $a: \V \times \V \to \mathbb N$ encoding the number of edges connecting two vertexes. 
In this case $a(3,1)=1$. This is the idea underlying the adjacency representation. 

\subsubsection{Adjacency representation} 
Each edge can be identified by specifying the vertex where it is starting, say $\mv_\ell$ and the vertex where it is ending, say $\mv_k$. 
So, each edge can be identified with a couple $(\mv_k,\mv_\ell)$. 
This identification can be made more formal by specifying a natural valued function $a:\V \times \V \to \mathbb N$ by
\begin{equation}\label{adjboolean}
a(\mv_\ell, \mv_k):=a_{k \ell}:=\#\{\me_j : \mbox{the edge $\me_j$ connects $\mv_k$ to $\mv_\ell$}\}.
\end{equation}
\begin{defin}\label{adjacency}
The \emph{adjacency matrix} of the graph $\G$ is the matrix $$A:=(a_{ij})_{i,j=1,\ldots,n}$$ defined as in \eqref{adjboolean}.
\end{defin}
The adjacency representation for the graph in Figure \ref{fig:example} is given by
$$
A=\left(\begin{array}{cccc}
0   &   1  &   1  &  0  \\
1 &   0    &   1  &  0  \\
0 &   0    &   0  &  1  \\
0   &   0   & 0  &  0
        \end{array}\right).
$$
\subsubsection{Incidence representation}
The adjacency representation has the disadvantage that it is not possible to determine the labels of the edges of the corresponding graph. 
We want to fix indexes for the edges in order to uniquely determine functions on the networks, and for that the adjacency representation is not suitable. 
In the incidence representation the \emph{incoming incidence matrix} $\mathcal I^+$  is defined by
\begin{equation}\label{incinc}
\iota^+_{kj}:=\begin{cases}
1 & \mbox{if the edge $j$ ends in the node $k$,}\\
0 & \mbox{otherwise.}
\end{cases}
\end{equation}
and its companion, the \emph{outgoing incidence matrix} $\mathcal I^-$, defined by
\begin{equation}\label{outinc}
\iota^-_{kj}:=\begin{cases}
1 & \mbox{if the edge $j$ starts in the node $k$,}\\
0 & \mbox{otherwise.}
\end{cases}
\end{equation}
\begin{defin}\label{incidence}
The \emph{incoming} (respectively, \emph{outgoing}) \emph{incidence matrix} of the graph $\G$ is the matrix $\mathcal I^+$, respectively $\mathcal I^-$, defined as in \eqref{incinc}, respectively \eqref{outinc}.
The \emph{incidence matrix} of the graph $\G$ is the matrix $\mathcal I=\mathcal I^+-\mathcal I^-.$
\end{defin}
As an example, we compute the incidence matrix of the graph in Figure~\ref{fig:example}
\begin{equation*}
\mathcal I=\left(\begin{array}{rrrrr}
1&-1&-1&0&0  \\
-1&1&0&-1&0  \\
0&0&1&-1&0  \\
0&1&1&0&0 
        \end{array}\right).
\end{equation*}
The incidence matrix has also a disadvantage: in contrast to the adjacency matrix, it is not possible to represent loops. 
The incoming and outgoing term would cancel each other, making impossible the detection of a loop in the matrix $\mathcal I$. 
However, it is possible to allow loops in the incidence representation, if only the ingoing and the outgoing incidence matrices are considered.

\subsubsection{Integration by parts on networks}
We now want to make a \emph{network} out of the graph in Figure~\ref{fig:example}. 
We define different spaces.

The Hilbert space $H:=(L^2(0,1))^5$, i.e.,  the space of square integrable functions defined on the edges of the graph will later be our state space.

Denote $\G:=\bigoplus_{i=1}^5[0,1]$. 
Via Lemma~\ref{Hrepresent}, $H$ is isomorphic to $L^2(\G)$, and so $\G$ can be identified with the measure-theoretic network.

Moreover, we want to define continuous functions on the network. To this aim define
$$
C(\G):=
\left\{\psi \in (C(0,1))^5: \exists d^\psi \in \mathbb C^4, \begin{array}{l}(\mathcal I^+)^\top d^\psi=\psi(0)\\(\mathcal I^-)^\top d^\psi=\psi(1)\end{array}\right\}.
$$
Later in this paragraph, we will show how the conditions involving the incidence matrices lead to functions that are continuous on the whole graph.

Now we can define weakly differentiable function on $\G$ by
$$
H^1(\G):=\left\{\psi \in (H^1(0,1))^5: \psi \in C(\G)\right\}.
$$

Finally, we define a space
$$
D(A):=\left\{\psi \in H^1(\G) \cap (H^2(0,1))^5: \mathcal I^+\frac{d}{dx}\psi(0)-\mathcal I^-\frac{d}{dx}\psi(1)=0\right\}.
$$
The two conditions involving the incidence matrices now lead to functions that are continuous in the nodes and satisfy the classical Kirchhoff node condition.

We now show that functions in $C(\G)$ are actually continuous in the nodes.
We fix a node, say $k=2$ and a function $\psi \in C(\G)$. 
Since it is in $C(\G)$ there exists $d^\psi$ as in the definition of $C(\G)$. 
With an abuse of notation, we can write
$$\me_1(1)=\mv_2, \quad \me_5(1)=\mv_2, \quad \me_2(0)=\mv_2.$$
Thus, we have to show that $\psi_1(1)=\psi_5(1)=\psi_2(0)$ for all $\psi \in C(\G)$. Since $\psi \in C(\G)$, $$\psi(0)=(\mathcal I^+)^\top d^\psi.$$ 
From this, compute
$$
d^\psi_2 =\psi_2(0)= \psi_5(1)=\psi_1(1),
$$
and so we are done. It is clear that it works also for the others vertexes.

We now define a form on the space $H^1(\G)$ and observe that the corresponding operator has domain $D(A)$. 
We are mainly interested in the Laplacian, it is thus natural to define $a: H^1(\G) \times H^1(\G) \to \mathbb C$ by
$$
a(\psi,\psi'):=\sum_{i=1}^5 \int_0^1 \frac{d}{dx}\psi_i(x)\overline{\frac{d}{dx}
\psi'_i(x)} dx=: 
\int_{\G} \frac{d}{dx}\psi(x)\overline{\frac{d}{dx}\psi'(x)} dx, 
$$
for all $\psi,\psi' \in H^1(\G).$
Our claim is the following.
\begin{prop}\label{illustrativeintegration}
The form $(a,V)$ is associated with the Laplacian on the edges of $\G$ with domain $D(A)$.
\end{prop}
\begin{proof}
Denote $(B,D(B))$ the operator associated with $(a,V)$. We show $D(A) \subset D(B)$. Fix an arbitrary $\psi \in D(A)$. 
Recall that whenever $\me_i$ end in $\mv_k$, the relation $f_i(1)=\iota_{ki}^+f(\mv_k)$ holds. 
For all $\psi' \in H^1(\G)$
\begin{eqnarray*}\label{parint1}
a(\psi,\psi')&=&\sum\limits_{i=1}^5 \int\limits_0^1 \frac{d}{dx}\psi_i(x)\overline{\frac{d}{dx}\psi'_i(x)}dx\\
&=&\sum\limits_{i=1}^5[\frac{d}{dx}\psi_i(x) \overline{\psi'_i(x)}]^1_0 - \sum\limits_{i=1}^5 \int\limits_0^1\frac{d^2}{dx^2}\psi_i(x)\overline{\psi'_i(x)}dx\\
&=&\sum\limits_{i,k} (\iota^+_{ki}-\iota^-_{ki})\frac{d}{dx}\psi_i(\mv_k)\overline{\psi'_i(\mv_k)} 
- \sum\limits_{i=1}^5 \int\limits_0^1\frac{d^2}{dx^2}\psi_i(x)\overline{\psi'_i(x)}dx.
\end{eqnarray*}
Since $\psi'\in C(\G)$, there exists $d^{\psi'} \in \mathbb C^5$ encoding the values of $\psi'$ in the nodes. Thus, we can compute
\begin{equation*}
\sum\limits_{i,k} (\iota^+_{ki}-\iota^-_{ki})\frac{d}{dx}\psi_i(\mv_k)\overline{\psi'_i(\mv_k)} 
=\sum\limits_{k=1}^4 d^{\psi'}_k \sum\limits_{i=1}^5(\iota^+_{ki}-\iota^-_{ki})\frac{d}{dx}\psi_i(\mv_k)=0,
\end{equation*}
since $\psi \in D(A)$. Choose $f= -(\frac{d^2}{dx^2}\psi_i)_{i=1,\ldots,5} \in L^2(\G)$. 
Then $a(\psi,\psi')=(f \mid \psi')_H$, i.e.,  $\psi \in D(B)$. Moreover, the operator $B$ acts as the Laplacian. 
The proof of the converse inclusion is also elementary.
\end{proof}

%% file: sec_graphoper.tex
\section{Graph and operator theoretic settings}\label{sec:graphoper}
We state in this section the general assumptions that we will use during the rest of this chapter.

Intuitively, a graph is a set of \emph{vertices} belonging to a set $\V$ connected by \emph{edges}, which we assume to belong to a set $\E$. 
We only consider edges of finite length. The extremal points of an edge are given by a map $\partial: E \to \V \times \V$.
Our final goal is to study global symmetry properties. For this, it is necessary to fix an orientation of the graph, which is implicitly given by the map $\partial$.
\begin{defin}
Consider sets $\V, \E$ and a map $\partial:\E \to \V \times \V$. Then, $(\V,\E,\partial)$ is said to be an oriented graph with vertex set $\V$ and edges set $\E$.
\end{defin}
This definition is a standard one in combinatorial graph theory, and it only requires $\V$ and $\E$ to be disjoint, without any cardinality assumption on the sets. 
However, we will only consider finite or countable graphs.
\begin{defin}\label{graphdefin}
Consider a graph $\G$. We fix the graph theoretic notation we are going to use. 
\begin{itemize}
\item Consider $\mv,\mv' \in \V$. 
If $(\mv,\mv')\in \range(\partial)$ or $(\mv',\mv)\in \range(\partial)$, we say that $\mv$ and $\mv'$ are \emph{adjacent} and we write $\mv \sim \mv'$.
The edge $\me$ such that tha above property holds is said to be \emph{incident} on both $\mv,\mv'$. 
We write indifferently $\me \sim \mv,\mv'$ or $\mv,\mv' \sim \me$.
\item The edge $\me$ is said to be \emph{outgoing} from $\mv$ if $\pi_1(\partial \me)=\mv$ and \emph{incoming} in $\mv$ if $\pi_2(\partial \me)=\mv$.
\item The \emph{incoming} (respectively, \emph{outgoing}) \emph{degree} of a node is the number of incoming (respectively, outgoing) edges. The \emph{degree} of a node is the sum of the incoming and outgoing degree. Notice that loops are counted twice.
They are denoted by $\deg{\mv}$, $\deg^+(\mv), \deg^-(\mv)$, respectively.
\item Consider subsets $\V'\subset \V,\E'\subset \E$. The \emph{subgraph $\G'$ induced by} $(\V',\E')$ is defined as following
$$
\mv \in \G' \Leftrightarrow \mv \in \V' \mbox{ or } \mv \sim \E',
$$
and
$$
\me \in \G' \Leftrightarrow 
\begin{cases} 
\me \in \E'
 \mbox{ or }\\
\exists \mv,\mv' \in \V' \mbox{ such that } \partial \me \in \{(\mv,\mv'), (\mv,\mv')\}.
\end{cases}
$$
We denote $\G'=g(\V',\E')$.
\item Consider a subgraph $\G'=(\V',\E')$ and the subgraph $\G''=(\V'',\E'')$ induced by $\E\setminus\E'$. We denote 
$$
\partial \G'=\partial \G'' = \V' \cap \V''.
$$
\item The sets $\Vfin^+, \V^+_\infty \subset \V$ are defined as
$$
\Vfin^+:=\{\mv \in \V: \deg^+(\mv)<\infty\}, \quad \V^+_\infty:=\{\mv \in \V: \deg^+(\mv)=\infty\}.
$$
The sets $\Vfin^-,\Vfin, \V^-_\infty, \V_\infty$ are defined analogously. 
\item The graph $\G$ is \emph{finite} if both $\V$ and $\E$ are finite.
\item The finite part $\Gfin$ of a graph is the subgraph induced by $(\Vfin,\emptyset)$, and $\Efin:=\{\me \in \E: \partial \me \subset \Vfin\}$.
\item The graph $\G$ is \emph{uniformly locally finite with degree $k$} if $\deg(\mv)\leq k$ for all $\mv \in \V$.
\item The \emph{incoming incidence matrix} is the map $\mathcal I^+: \E \times \V \to \{0,1\}$ defined by
$$
\iota(\me,\mv):=\iota_{\me\mv}:= \#\{\me \in \E, \pi_2(\partial(\me))=\mv\}
$$
The \emph{outgoing incidence matrices} $\mathcal I^-$ is defined analogously, 
and the \emph{incidence marices} $\mathcal I$ is defined by $\mathcal I= \mathcal I^+-\mathcal I^-.$
\item The \emph{incoming incidence list} is the map $\Gamma^+: \V \to \mathcal P(\E)$ defined by
$$\Gamma^+(\mv):=\{\me \in \E: \me(0)=\mv\}.$$
The \emph{outgoing incidence list} $\Gamma^-$ and the \emph{incidence list} of the graph is defined analogously.
%\item For a graph $\G'$, the dual graph $\G'$is defined by taking $\G':=(\E,\V,\partial')$, where $\partial'$ is choosen in such a way that
%$$\mathcal I'^+(\mv,\me):=\mathcal I^+(\me,\mv), \quad \mathcal I'^-(\mv,\me):=\mathcal I^-(\me,\mv), \qquad \mv \in \V, \me\in \E.$$
\item Fix two nodes $\mv,\mv' \in \V$. A \emph{path} $P$ is defined as a sequence of nodes $P:=(\mv=\mv_0,\ldots,\mv_n=\mv')$ such that $\mv_k$ and $\mv_{k+1}$ are adjacent. The length $l(P)$ of the path is $n$.
\item Analogously, fix two edges $\me,\me' \in \E$. 
A \emph{path} $P$ is defined as a sequence of nodes $P:=(\me \sim \mv_1,\ldots,\mv_n\sim\me')$ such that $\mv_k$ and $\mv_{k+1}$ are adjacent. 
The length $l(P)$ of the path is $n$.
\item We denote $P[\mv,\mv']$ the set of all paths from $\mv$ to $\mv'$ and $P[\me,\me']$ the set of all paths from $\me$ to $\me'$.
\item The distance $d(\mv,\mv')$ of two nodes is defined as
$$
d(\mv,\mv'):=\min_{P \in P[\mv,\mv']} l(P).
$$
If no path connects $\mv$ to $\mv'$, then $d(\mv,\mv'):=\infty$.
\item The distance between two edges is defined analogously.
\item The distance of two subgraphs $d(\G',\G'')$ is defined by
$$
d(\G',\G''):=\inf_{\mv' \in \G',\mv'' \in \G''} d(\mv',\mv'').
$$
\item A graph is \emph{connected} if $d(\mv,\mv')<\infty$ for all $\mv,\mv' \in \V$.
\end{itemize}
\end{defin}
This purely combinatorial setting is not suitable to define a network diffusion equation, due to lack of a metric structure. To this aim, we identify each edge with a bounded interval, which we parametrize as $[0,1]$. 
The endpoints of the interval are consequently identified with the vertexes. Fix now a graph $\G$, with finite or countably many edges $\me \in \E$. 
Following the approach of the first chapter, we can define $L^2(\G)$ as the product space of the $L^2(0,1)$ spaces corresponding to the single intervals.
\begin{defin}\label{l2grafo}
Consider a finite or countable graph $\G$. We define
$$
L^2(\G):=\bigoplus_{\me \in \E} L^2(0,1).
$$
For functions $\psi \in L^2(\G)$ we consequently write $\psi = (\psi_{\me})_{\me \in \E}$.
\end{defin}
We already observed in Section~\ref{sec:nondiag} that sesquilinear forms defined on graphs represent an important class of form matrices with non-diagonal domain. However, Definition \ref{l2grafo} does not yet contain any information about the topology of the graph, that is instead involved in the definition of the form domain. 

Consider now the space
$$
V_0:= \bigoplus_{\me \in \E} H^1(0,1).
$$
As a consequence of the boundedness of the trace operator on $H^1(0,1)$ both $\psi(0)$ and $\psi(1)$ are in $\ell^2(\E)$, and so the Definitions~\ref{adjacency} and~\ref{incidence} can be interpreted also in the case of a countable graph, if the corresponding matrices are allowed to be infinite. 
In fact, they have to be understood as (possibly unbounded) operators on the space from $\ell^2(\E)$ to $\ell^2(\V)$.

We now define a form domain by 
\begin{equation}\label{formdomain}
V:=\left\{ \psi \in V_0 :\exists d^\psi \in \ell^2(\V), \begin{array}{l}(\mathcal I^+)^\top d^\psi=\psi(0)\\(\mathcal I^-)^\top d^\psi=\psi(1)\end{array} \right\}.
\end{equation}
In the case of an infinite graph, $\mathcal I^{+^\top}, \mathcal I^{-^\top}$ denote the transpose of the incidence operators.
Consequently, the existence of a $d^\psi$ with the required properties has to be understood as the existence of $d^\psi$ in the domain of the incidence operators. 
We discuss in the next paragraph some issues concerning the operator theoretical properties of the incidence operators.

\subsubsection{Incidence matrices as operators}\label{incidenceop}
We start by specifying in which space the vectors $\psi(0), \psi(1)$ must be thought.

\begin{lemma}\label{inclusionboundaryspace}
The spaces
$$
\partial_0 V :=\{\psi(0): \psi \in V \}, \qquad \partial_1 V :=\{\psi(1): \psi \in V \}
$$
are isomorphic to subspaces of $\ell^2(\E)$.
\end{lemma} 
\begin{proof}
We show that $V \hookrightarrow C([0,1],\ell^2)$. The inclusion is clear by the definition of $V$. 
Further, the identity $i:V \to C([0,1],\ell^2)$ is a continuous mapping since it is continuous in each component with the same continuity constant,

So, we obtain $\psi(x) \in \ell^2$ for all $x \in [0,1]$ and the claim follows.
\end{proof}
It is a natural question whether the incidence operators have or not dense domain. This is the case, independently of the graph.
\begin{prop}\label{denseinc}
Consider a graph $\G$. Both $\mathcal I^+$ and $\mathcal I^-$ have dense domain as operators from $\ell^2(\E)$ to $\ell^2(\V)$.
\end{prop}

\begin{proof}
We split the proof in several cases. First, assume that the graph $G$ is locally finite. Then,
$$
\G= \bigcup_{n \in \mathbb N} \G_n:=\bigcup_{n\in \mathbb N}g(\V_n,\emptyset) ,\qquad \V_n:=\{\mv \in \V: \deg(\mv) \leq n\}.
$$
Further, for all $y \in \ell^2(\E)$
$$
\lim_{n \to \infty} \| y- \pi_{\G_n}(y)\|_{\ell^2(\E)}=0.
$$
The estimate
\begin{eqnarray*}
\|\mathcal I^+ \pi_{\G_n}(y)\|^2_{\ell^2(\V)}& = &\sum_{\mv \in \V_n }| \Big( \sum_{\me \in \Gamma^+(\mv)} (\pi_{\G_n}(y))_{\me}\Big) |^2 \\
&\leq & M \sum_{\mv \in \V_n}  \sum_{\me \in \Gamma^+(\mv)} |(\pi_{\G_n}(y))_{\me} |^2  \\
& = & M\|\pi_{\G_n}(y)\|^2_{\ell^2(\E)}
\end{eqnarray*}
yields that $\pi_{\G_n}(y)\in \dom(\mathcal I^+)$ for all $n$, and so the claim is proved for a locally finite graph.

\smallskip
If the graph consist of a single, incoming infinite star $S$ (see Definition~\ref{graphdef}), then $\ell^1(\E) \subset \dom(\mathcal I^+)$, 
and so $\dom(\mathcal I^+)$ is dense in $\ell^2(\E)$. So, the claim is true for such graphs.

\smallskip
Finally, for a general graph, fix a numbering $\V^+_\infty:=\{\mv_k : k \in \mathbb N\}$ decompose the graph as
\begin{equation}\label{hilfdenseinc1}
\G=(\bigcup_{k \in \mathbb N} S_k) \cup g(\Vfin, \emptyset), \quad S_k:=\{\me \in \E: \me \in \Gamma^+(\mv_k)\},
\end{equation}
and define for all $n \in \mathbb N$ the approximations 
$$
\G_n = \bigcup_{k \leq n-1} S_k \cup g(\Vfin,\emptyset).
$$
Fix an arbitrary $x \in \ell^2(\E)$. We define 
$$
v_0:=( \pi_{g(\Vfin,\emptyset)}(x)) , \quad v_k:=(\pi_{S_k}(x)).
$$
Since $g(\Vfin,\emptyset)$ is locally finite, there exists a sequence 
$$
(v^n_0)_{n \in \mathbb N} \in \dom(\mathcal I^+_{|g(\Vfin,\emptyset)}) \cap \ell^2(\Efin)
$$
such that
$$\|v^n_0 - v_0\| \leq \frac{1}{2^{n}}, \qquad n \in \mathbb N.$$
In particular, $\pi^{-1,r}_{\Efin,0}(v_0) \in \dom(\mathcal I^+)$, since $\mathcal I^+ \ell^2(\Efin) \subset \ell^2(\Vfin$).

Since the domain on infinite stars is also dense, for all $k \geq 1$ 
there is a sequence $(v^n_k)_{n \in \mathbb N} \in \dom(\mathcal I^+_{|S_k})$ such that 
$$\|v^n_k - v_k\| \leq \frac{1}{2^{n+k}}, \qquad k\geq 1, n \in \mathbb N.$$
Again by the same arguments as for finite part, $\pi^{-1,r}_{\Efin,0}(v_k) \in \dom(\mathcal I^+)$

\smallskip
Summing up, for all $k \in \mathbb N$ there is a sequence $(v^n_k)_{n \in \mathbb N} \in \dom(\mathcal I^+)$ such that 
$$\|v^n_k - v_k\| \leq \frac{1}{2^{n+k}}, \qquad k,n \in \mathbb N.$$
Define $x^n:=\sum_{k \leq n} v_k$ and fix $\varepsilon <0$.
Since~\eqref{hilfdenseinc1} holds, there exists $n_0 \in \mathbb N$, $\|x-\pi_{\G_n}(x)\| < {\varepsilon}$ for all $n \geq n_0$. 
For such an $n$ estimate
\begin{eqnarray*}
\|x-x^n\| & = & \|x-\pi_{\G_n} (x) + \pi_{\G_n}(x) - x^n\| \\
& \leq&  \|x-\pi_{\G_n} (x)\|+\|\pi_{\G_n}(x) - x^n\| \\
& <& {\varepsilon} + \frac{1}{2^n},
\end{eqnarray*}
and so $\lim_{n \to \infty }x^n = x$. 

\smallskip
We still have to prove $x^n \in \dom(\mathcal I^+)$, but this is clear since $x^n$ is a finite linear combination of elements in the domain.
\end{proof}
In the next result we investigate the boundedness of the incidence matrices.
\begin{prop}\label{incidencebound}
The following assertions hold.
\begin{enumerate}
\item The incidence matrices $\mathcal I^+, \mathcal I^-$ are bounded operators from $\ell^2(\E)$ to $\ell^2(\V)$ if and only if the graph $\G$ is uniformly locally finite.
\item The incidence matrices $\mathcal I^+, \mathcal I^-$ are bounded operators in $\ell^\infty(\E)$ to $\ell^\infty(\V)$ if and only if the graph $\G$ is uniformly locally finite.
\item The operators $\mathcal I^+, \mathcal I^- $ are contractive from $\ell^1(\E)$ to $\ell^\infty(\V)$ for every countable graph $\G.$
\end{enumerate}
\end{prop}

\begin{proof}
We stat proving a). Assume the graph $\G$ to be uniformly locally finite with a degree bound $D$, fix $x \in \ell^2(\E)$ and compute
\begin{eqnarray*}
\|\mathcal I^+x\|^2_{\ell^2(\V)}&=& \sum_{\mv \in \V}  |\sum_{\me \in \Gamma^+(\mv)}  x_{\me}|^2\\
&\leq& \sum_{\mv \in \V}  \|(x_{\me})_{\me \in \Gamma^+(\mv)}\|_{\ell^1(\Gamma^+(\mv))}^2\\
&\leq& \sum_{\mv \in \V} \deg^+(\mv) \|(x_{\me})_{\me \in \Gamma^+(\mv)}\|_{\ell^2(\Gamma^+(\mv))}^2\\
&\leq& D\|x\|^2_{\ell^2(\E)}.
\end{eqnarray*}
Assume now that the graph is not uniformly locally finite. 
We distinguish two cases. First, assume that there exists a node $\mv$ such that $\deg^+(\mv)=\infty$. 
It suffices to prove the claim for  an inward star $\G$ with center $\mv_1$ and infinitely many leaves. 
In this case, all vectors $0 \leq x \in \ell^2\setminus \ell^1$ are not in the domain of $\mathcal I^+.$ 
Alternatively, assume that there exists a sequence of nodes $(\mv_\ell)_{\ell \in \mathbb N}$ such that $\lim_{\ell \to \infty}\deg(\mv_\ell)=\infty$. 
Consider the vectors $x_{\ell}:= \frac{1}{\sqrt{\deg^+(\mv_\ell)}}\mathbb 1_{I(\mv_\ell)}.$ Then $\|x_{\ell}\|_{\ell^2(\E)}=1$, but $\|\mathcal I^+x_{\ell}\|_{\ell^2(\V)}=\deg^+(\mv_\ell).$ 
This shows that $\mathcal I^+$ is not bounded.

To see that b) holds observe that the operator $\mathcal I^+$ is a positive operator. Thus, it is sufficient to compute $\mathcal I^+ \mathbb 1_{\me}= (\deg^+ \mv)_{\mv \in \V}$.

To see that c) holds observe that $\mathcal I^+$ is a positive operator from $\ell^1(\E)$ to $\ell^1(\V)$, which is isometric on the positive cone. 
Thus, for arbitrary $x \in \ell^1$,
$$
\|\mathcal I^+x\|_{\ell^\infty}\leq \|\mathcal I^+x\|_{\ell^1} = \|x\|_{\ell^1}.
$$
\end{proof}

We turn now our attention to the transpose operators $(\mathcal I^+)^\star,(\mathcal I^-)^\star$.
%Observe that the transpose incidence operators acts on a function defined on the nodes by copying the value of the function into all edges that go out from this node.

\begin{prop}
Consider a graph $\G$. The transpose incidence operators $(\mathcal I^+)^\star,(\mathcal I^-)^\star$ are densely defined if and only if $\G$ is locally finite.
\end{prop}

\begin{proof}
Assume that $\G$ is locally finite and fix a vector $x \in c_{00}(\V)$. So, also $(\mathcal I^+)^\star \in c_{00}(\E) \subset \ell^2(\E)$.
This implies $x \in D((\mathcal I^+)^\star)$ and so the operator is densely defined. Conversely, assume that $(\mathcal I^+)^\star $ is densely defined. In particular, $\mathbb 1_{\mv}$ has to be in the domain for all $\mv \in \V$. 
Observe that $(\mathcal I^+)^\star \left( \mathbb 1_{\mv} \right)=\mathbb 1_{\Gamma^+(\mv)}$ which is in $\ell^2$ if any only if $\Gamma^+(\mv)$ is finite.
\end{proof}

\subsubsection{Continuity properties}
We can now investigate continuity properties of the sesquilinear form $(a,V)$.
\begin{ass}
During the rest of the section we always assume the following.
\begin{itemize}
\item $C$ is an operator valued function $C: x \mapsto (c_{ij}(x))_{i,j \in \E} $ of class $C^1$.
\item $M=(m_{k \ell})_{k,\ell \in V}$ is an operator in $\ell^2(V)$.
\item The sesquilinear form $a:V \times V \to \mathbb C$ is defined by
\begin{equation}\label{networkform}
a(\psi,\psi')=\left(C \frac{d}{dx} \psi \mid \frac{d}{dx} \psi'\right)_{L^2(\G)} - (Md^\psi \mid d^{\psi'})_{\ell^2(V)}.
\end{equation}
\end{itemize}
\end{ass}

We observe that on nodes with infinite degree Dirichlet node conditions are automatically imposed for all functions in the form domain.

\begin{prop}\label{infinitedirichlet}
Consider a graph $\G$ such that there exists $\mv \in \V$ with $\deg (\mv) = + \infty$. Then $\psi(\mv)=0$ for all $\psi \in V$.
\end{prop}
\begin{proof}
Recall that $H^1(0,1) \hookrightarrow C[0,1]$ and so $\|\psi\|_{H^1(0,1)} \geq M \|\psi\|_\infty$.

Fix an arbitrary node $\mv \in \V$ and assume that there exists $\psi \in V$, $\psi(\mv)\neq0$. Then
$$
\|\psi\|^2_V = \sum_{\me \in \E} \|\psi_{\me}\|^2_{V_{\me}} \geq  M \sum_{\me \in \Gamma(\mv)} \|\psi_{\me}\|^2_\infty \geq M \sum_{\me \in \Gamma(\mv)} |\psi(\mv)|^2,
$$
and so $\Gamma(\mv)$ has to be finite.
\end{proof}
Consequences of this results are discussed in Section~\ref{sec:infinitesymmetries}
\begin{prop}
The sesquilinear form $(a,V)$ is continuous if the function $C$ is essentially bounded in $ \mathcal L(\ell^2(\E))$ 
and $M \in \mathcal L(\ell^2(\Vfin)).$
\end{prop}
\begin{proof}
First observe that, since $d^\psi_{\mv} =0$ whenever $\deg(\mv)= \infty$, then $M$ has only to be investigate on the orthogonal complement of $\ell^2(\V_\infty)$.
So, the sufficiency of the two conditions for the continuity of the form is clear.
\end{proof}

The next step is to identify the operator associated with the sesquilinear form $(a,V)$. 
Observe that the matrix $C$ causes difficulties in the integration by parts: on the one hand, the node conditions are requiring some boundary coupling of the different components $\psi_i$ of a function $\psi \in V$. 
On the other hand, for non-diagonal $C$ some internal coupling of the functions is prescribed. 
These issues are the object of the next section.

%% file: sec_operatordomain.tex
\section{Identification of the network operator domain}\label{sec:operatordomain}
Aim of this section is to identify the domain of the operator associated with the form defined in \eqref{networkform}, 
thus deriving the boundary conditions that can be used in the formulation of a parabolic equation on the network. 

In this section we fix a numbering both of the vertices and of the edges, i.e. $\V:=\{\mv_1,\mv_2,\ldots \}$ and $\E:=\{\me_1,\me_2,\ldots\}$, and
we want to derive the expression of the \emph{incidence tensor} for a strongly coupled network equation.

Integrate by part the expression of the sesquilinear form
\begin{eqnarray*}
a(\psi,\psi') &=&
\sum\limits_{i,j \in \E} \int\limits_0^1c_{ji}(x)\frac{d}{dx}\psi_i(x)\overline{\frac{d}{dx}\psi'_j(x)}dx - \sum\limits_{k,\ell \in \V} m_{k\ell}d^\psi_k\overline{d^{\psi'}_\ell}\\
&=&\sum\limits_{i,j \in \E} [c_{ji}(x)\frac{d}{dx}\psi_i(x)\overline{\psi'_j(x)}]^1_0\\
&-&\sum\limits_{i,j \in \E} \int\limits_0^1\frac{d}{dx}(c_{ji}\frac{d}{dx}\psi_i)(x)\overline{\psi'_j(x)}dx
-\sum\limits_{k,\ell \in V} m_{k\ell}d^\psi_\ell\overline{d^{\psi'}_k}.
\end{eqnarray*}
It remains to derive the dependence on the incidence matrices of the first term in the sum. To this aim we define tensors 
$$\mathfrak I^-=(\hat{\iota}^{kj}_{\ell i})_{i,j \in \E, k,\ell\in \V}, \qquad\mathfrak I^-=(\check{\iota}^{kj}_{\ell i})_{i,j \in \E, k,\ell\in \V}$$ by setting
$$
\hat{\iota}^{kj}_{\ell i}:=
\begin{cases}
1, & \mbox{ if } i_i(1)=\ell_\ell \mbox{ and } i_j(1)=\ell_k,\\
0, & \mbox{otherwise.}
\end{cases}
$$
The tensor $\mathfrak I^-$ is defined analogously. 
In a more compact way we also write $\mathfrak I^+ = \mathcal I^+ \otimes \mathcal I^+$ and $\mathfrak I^- = \mathcal I^- \otimes \mathcal I^-$. 
For this reason, we call $\mathfrak I^+$ the \emph{incoming incidence tensor} and $\mathfrak I^+$ the \emph{outgoing incidence tensors}. 
In fact, the incoming incidence tensor contains all possible pairs of couples $(i_i, \ell_\ell)$ such that $i_i(1)=\ell_\ell$, and the outgoing can be interpreted analogously. 

Further, we introduce the \emph{weighted incidence tensor} 
$$
\mathfrak W=(\omega^{k j}_{\ell i})_{i,j \in \E, k,\ell\in \V}, \quad 
\omega^{k j}_{\ell i}:=c_{ij}(\ell_\ell) \iota^{k j}_{\ell i}, \quad  k,\ell \in \V, i,j \in \E.
$$
Using these conventions, we obtain
\begin{eqnarray*}
\sum\limits_{i,j \in \E} [c_{ji}(x)\frac{d}{dx}\psi_i(x)\overline{\psi'_j(x)}]^1_0&=&
\sum\limits_{i,j \in \E}\sum\limits_{\ell, k \in \V} c_{ji}(\ell_\ell) (\hat{\iota}^{kj}_{\ell i}-\check{\iota}^{kj}_{\ell i}) \frac{d}{dx}\psi_i(\ell_\ell)\overline{\psi'_j(\ell_k)}\\
&=&\sum\limits_{i,j \in \E}\sum\limits_{k,\ell \in \V} {\omega}^{kj}_{\ell i} \frac{d}{dx}\psi_i(\ell_\ell)\overline{\psi'_j(\ell_k)}\\
&=&\sum_{k \in \V} \overline{d^{\psi'}_k}\sum_{i,j \in \E}\sum_{\ell \in \V} \omega^{kj}_{\ell i} \frac{d}{dx}\psi_i(\ell_\ell)\\
&=&\sum_{k\in \V} \overline{d^{\psi'}_k} \sum_{\ell\in \V}\sum_{i,j \in \E}\omega^{kj}_{\ell i} \frac{d}{dx}\psi_i(\ell_\ell).\\
\end{eqnarray*}
Now we have derived the correct condition. Assume that
\begin{equation}\label{nonlockirch}
\sum_{\ell \in V} \sum_{i,j \in \E}\omega^{kj}_{\ell i} \frac{d}{dx}\psi_i(\ell_\ell)=\sum_{\ell \in V} m_{k\ell} d^{\psi}_\ell, \qquad \psi \in D(A), k \in \V.
\end{equation}
Then,
$$
\sum_{k\in \V} \overline{d^{\psi'}_k} \sum_{\ell\in \V}\sum_{i,j \in \E}\omega^{kj}_{\ell i} \frac{d}{dx}\psi_i(\ell_\ell)=
\sum_{ k, \ell\in \V}m_{k\ell} d^\psi_{\ell} d^{\psi'}_k.
$$
Define now
\begin{equation}\label{networkoperdomain}
D(A):=\left\{ \psi \in V \cap H^2(\G): \psi \mbox{ satisfies \eqref{nonlockirch}} \right\}.
\end{equation}
We have shown that $D(A) \subset D(B)$, where $B$ is the operator associated with the sesquilinear form $a$
and $A$ is the operator 
\begin{equation}\label{operaction}
A\psi:=\left( \sum_{j \in \E} \frac{d}{dx}(\frac{d}{dx} c_{ji}(x) \psi_i(x)) \right)_{i \in \E}^\top, \qquad \psi \in D(A) .
\end{equation}
Since $d^{\psi'}$ is arbitrary, also the converse inclusion can be proved.
\begin{thm}
Assume $C$ to be an uniformly bounded operator-valued function and $M$ to be a bounded operator.
Then, the operator $(A,D(A))$ defined as in \eqref{networkoperdomain}, \eqref{operaction} is the operator associated with the form $a$.
\end{thm}
\begin{proof}
Denote $(B,D(B))$ the operator associated with the form $(a,V)$.We have already proved that $D(A) \subset D(B).$ Consider now $\psi \in D(B)$.
By definition, there exists $f \in H$ such that $a(\psi,\psi')=-(f \mid \psi')$ for all $\psi' \in V$.
Thus,
$$
\sum_{i,j \in \E} \int_0^1 c_{j i}(x) \frac{d}{dx} \psi_{i}(x) \overline{\frac{d}{dx} \psi'_{j}(x)}dx
- \sum_{\ell,k \in \V} m_{\ell k} d^\psi_{k} \overline{d^{\psi'}_{\ell}} =
- \sum_{i \in \E} \int_0^1 f_{i}(x) \overline{\psi'_{i}(x)} dx.
$$
Integrating by part the left hand side, we obtain that
\begin{eqnarray*}
- \sum_{i \in \E} \int_0^1 f_{i}(x)\overline{\psi'_{i}(x)}dx &=&
\sum_{\ell \in \V}\overline{d^{\psi'}_{\ell}}\sum_{i,j \in \E} \sum_{k \in \V} \omega^{\ell j}_{k \ell} \frac{d}{dx} \psi_{i}(\mv_k) \\
&&- \sum_{i,j \in \E} \int_0^1\frac{d}{dx}(c_{ji}\frac{d}{dx}\psi_{i})(x)\overline{\psi'_{j}(x)}dx \\
&&- \sum_{ \ell,k \in \V} m_{\ell k} d^\psi_{k} \overline{d^{\psi'}_{\ell}}.
\end{eqnarray*}
This holds for all $\psi' \in V$.
In particular, considering $\psi' \in V$  vanishing on all but one edge of the network, we conclude that 
\begin{equation*}\label{coincideriv}
f_{i}(x)=\sum_{j \in \E}\frac{d}{dx}(c_{j i}\frac{d}{dx}\psi_{j})(x) \qquad \mbox{for all } x\in (0,1) \mbox{ and all }i \in E.
\end{equation*}
Similarly, considering $\psi'$ with arbitrary nodal values and arbitrary small $H$-norm, we obtain
\begin{equation*}
\sum_{i,j \in \E} \sum_{k \in \V} \omega^{\ell j}_{k i} \frac{d}{dx}\psi_{i}(\mv_k) 
- \sum_{k \in \V} m_{\ell k}d^\psi_{k } = 0\qquad\mbox{for all } \ell \in \V.
\end{equation*}
This shows that $\psi \in D(A)$ and completes the proof.
\end{proof}

%% file: sec_systemsII.tex
\section{Systems of network equations}\label{sec:systemsII}
In this section we want to extend the arguments presented in Section \ref{sec:systemsI} to the case diffusion taking place on a network. 
We start characterizing the well-posedness in $H:= \bigoplus_{\me \in \E} L^2(0,1)$ of the abstract Cauchy problem
\begin{equation}\label{acpgraph}
\left\{\begin{array}{rclr}\dot{u}(t)&=& Au(t), & t \geq 0, \\
u(0)&=& u_0, & u_0 \in H.\end{array}\right.
\end{equation}
Here $A$ is the operator associated with the form $(a,V)$ defined in \eqref{networkform} on the domain defined in \eqref{formdomain}. 
According to the theory of sesquilinear forms, the problem~\eqref{acpgraph} is well--posed if the form $(a,V)$ is continuous, elliptic and densely defined. 
First observe that the form domain contains $\ell^2(\E, H^1_0(0,1))$. Since the latter is dense in $H$, the form $(a,V)$ is densely defined.
The continuity of the form has already been discussed in Section~\ref{sec:graphoper}, 
whereas the characterization of ellipticity given in Section~\ref{sec:systemsI} also works in the case of infinite networks. 
We thus assume that the form $(a,V)$ is densely defined, continuous and elliptic and we turn our attention to the analysis of the qualitative properties of \eqref{acpgraph}. 
In particular, we want to discuss symmetry properties of the system. 
To this end, we introduce a concept of symmetry which is suitable for network equations.

Recall that in the elementary calculus a function $f : \mathbb R \to \mathbb R$ is called even 
if it is symmetric with respect to the $y$-axis, i.e., if $f(x)=f(-x)$ for all $x\in \mathbb R$. 
Assume now the function $f$ to be in  $ H^1(-1,1)$. 
Then $f$ is also continuous, and thus the definition makes sense also for functions in $H^1(-1,1)$. 
More generally, assume that the graph $\G$ is an outward star, see Definition~\ref{graphdef}.
We denote $\mv_0$ the centre of the star and we number the other edges such that $\mv_\ell \sim \me_\ell$ for all $\ell=1,\ldots, |\E|$.
We say the a function $f: \G \to \mathbb C$ is \emph{even} if $f_i(x)=f_j(x)$ for all $x\in [0,1], i,j=1,\ldots,|\E|$.

\begin{lemma}
A function $f \in V$ on a star $\G$ is even if and only if $f(x)\in \langle \mathbb 1 \rangle$ for all $x\in [0,1]$.
\end{lemma}
\begin{proof}
Let $f \in V$ even. Then $f_i(x)=f_j(x)$ for all $i,j \in \E, x\in (0,1)$. 
So, there exists $\lambda \in H^1(0,1)$ such that $f_i(x)=\lambda(x)$ for all $i \in \E$. This means $f(x) \in \langle \mathbb 1 \rangle$ for all $x$.

Conversely, if $f(x) \in \langle \mathbb 1 \rangle$ for all $x$, then for all $x$ there exists $\lambda(x)$, $f(x)=\lambda (x) \mathbb 1$ and so
$f_j(x)=f_i(x)$ for all $i,j\in \E$, $x \in (0,1)$.
\end{proof}

Observe now that this concept also makes sense for $L^2$-functions if it is understood almost everywhere

We now generalise this concept. For an arbitrary linear subspace $Y \subset {\ell^2(\E)}$, 
we interpret $Y$ also as a subspace $\mathcal Y$ of $L^2$ by considering those functions that belong to $Y$ almost everywhere.
\begin{defin}
Consider a function $\psi \in H$ and a closed linear subspace $Y \subset \ell^2(\E)$. We say that $\psi$ is \emph{$Y$-symmetric} if 
$$
\psi(x) \in Y, \qquad \mbox{for almost all } x\in (0,1). 
$$
Moreover, we define
$$
\mathcal Y:=\{\psi \in H: \psi \mbox{ is $Y$-symmetric}\}.
$$
Further, consider the semigroup of the solution operators $(e^{ta})_{t \geq 0}$ of the problem~\eqref{acpgraph}.
 We say that the semigroup $(e^{ta})_{t \geq 0}$ is \emph{$Y$-symmetric} if the subspace $\mathcal Y$ is invariant under the action of $e^{ta}$ for all $t \geq 0$.
\end{defin}

In the following, $H$ will always denote the product $\bigoplus_{\me \in \E} L^2(0,1)$ and $V$ the form domain as in \eqref{formdomain}. 
It is useful to decompose the space $H=\mathcal Y\bigoplus \mathcal Y^\perp$. 
We formulate a result analogous to Proposition~\ref{L2decomposition}.
\begin{lemma}
Let $Y \subset \mathbb \ell^2(\E)$ a closed linear subspace. 
Fix orthonormal bases $(x'_{\me})_{\me \in \E'}$ of $Y$ and $(x''_{\me})_{\me \in E''}$ of $Y^\perp$. 
Then
$$
\mathcal Y= \{ \sum_{\me' \in \E'} \psi_{\me'}x'_{\me'}: (\psi_{\me'})_{\me' \in \E'} \in \ell^2(\E', L^2(0,1)) \},
$$
and
$$
\mathcal Y^\perp= \{ \sum_{\me \in \E''} \psi_{\me''}x''_{\me''}: (\psi_{\me''})_{\me'' \in \E''} \in \ell^2(\E'', L^2(0,1)) \}.
$$
\end{lemma}
\begin{proof}
Since the subspace $Y$ is closed, it is an Hilbert space with orthonormal basis $\{(x_{\me'})_{\me' \in \E'}\}$. 
So, the claim follows by Proposition \ref{L2decomposition}.
\end{proof}
Observe that having computed $\mathcal Y, \mathcal Y^\perp$ allows us to apply Corollary~\ref{ortho} for the investigation of the invariance properties of $\mathcal Y$.
In particular, for a system of network equations, the second condition in Corollary~\ref{ortho} is equivalent to the statement b) in Theorem~\ref{symmetries}.

In contrast to the case of systems associated with forms with diagonal domains, however, the condition $P_\mathcal Y V \subset V$ is not automatically fulfilled. 

In the case of network equations, one has in particular to prove
that if $\psi$ is a function which is continuous in the nodes, then also $P_\mathcal Y\psi$ is continuous in the nodes.
If this is the case, we call the projection $P_\mathcal Y$ \emph{admissible}.

We show now how it is possible to reduce the issue of admissibility to linear algebraic statements.
We start by proving a linear algebraic auxiliary result.

\begin{lemma}\label{decomp}
Consider an orthogonal projection $K\in \mathcal L(\ell^2(\E))$ and a closed linear subspace $Y \subset \ell^2(\E)$. 

Then the following assertions are equivalent.
\begin{enumerate}
\item $KY \subset Y$.
\item $Y=(\Ker K \cap Y) \bigoplus (\range  K \cap Y)$.
\end{enumerate}
\end{lemma}
\begin{proof}
To see that b) implies a) fix an arbitrary $y \in Y$. By hypothesis $y=y_1 +y_2$, where  $y_1 \in \Ker K \cap Y$ and $y_2 \in 
\range  K \cap Y$. Then $Ky=Ky_1 + Ky_2=y_2 \in Y$, which proves the claim.

\smallskip
We now prove that a) implies b). Since $KY\subset Y$, the restriction $K_{|Y}$ is a linear operator on $Y$. 
In particular, since $Y$ is a closed subspace, then $Y$ is an Hilbert space and $K_{|Y}$ is the orthogonal projection onto its own range.
As a consequence $Y=\Ker K_{|Y} \bigoplus \range K_{|Y}$.  
Obviously $\Ker K_{|Y} = \Ker K \cap Y$ and, by a), $\range K_{|Y} = \range K \cap Y$ because $K$ is a projection. The claim follows.
\end{proof}

To characterise admissibility of projections we introduce some additional notation. 
We define the operators 
$$
\It:\ell^2(\V) \to \ell^2(\E) \times \ell^2(\E)
\quad \mbox{and} \quad
\widehat{P_Y} : \ell^2(\E) \times \ell^2(\E) \to  \ell^2(\E) \times \ell^2(\E)
$$
by
\begin{equation}\label{defwidehat}
\It:=(\mathcal I^+,\mathcal I^-)^\top=\begin{pmatrix} (\mathcal I^+)^\top \\(\mathcal I^-)^\top\end{pmatrix} \qquad \text{and} 
\qquad \widehat{P_Y}:=\K.
\end{equation}
Observe that $\widehat{P_Y}$ is an orthogonal projection.

\begin{thm}\label{charadmiss}
If the graph $\G$ is connected, then the following assertions are equivalent.
\begin{enumerate}[a)]
\item\label{charadmiss:a} The projection ${P}_{\mathcal Y}$ is admissible.
\item\label{charadmiss:b} The range of $\It$ is invariant under $\widehat{P_Y}$, i.\,e., 
$$\K \range \It \subset 
\range \It.$$
\item\label{charadmiss:c} There exists a basis of $\range \It$ consisting of eigenvectors of $\widehat{P_Y}$.
\end{enumerate}
\end{thm}

\begin{proof}
We start by proving the equivalence of (\ref{charadmiss:a}) and (\ref{charadmiss:b}). 
Recall that for every $\psi \in V$ there exists a vector $d^\psi \in \ell^2(\E)$ such that
$$(\mathcal I^+)^\top d^\psi=\psi(0), \qquad (\mathcal I^-)^\top d^\psi=\psi(1).$$
The admissibility of the projection is equivalent to the fact that for every $\psi \in V$ there exists a vector $d^{P_{\mathcal Y}\psi} \in \ell^2(\V)$ such that
$$(\mathcal I^+)^\top d^{ {P}_{\mathcal Y}\psi}= {P}_{\mathcal Y}\psi(0)=P_Y\psi(0), 
\qquad (\mathcal I^-)^\top d^{{P}_{\mathcal Y}\psi}={P}_{\mathcal Y} \psi(1)=P_Yf(1).$$
Inserting the first equation into the second and observing that for all $x \in \ell^2(\V)$ 
there exists a function $\psi \in \V$ which is continuous in the nodes and such that $d^\psi=x$,
we obtain that (\ref{charadmiss:a}) is equivalent to the fact that for all $x \in \ell^2(\V)$ 
there exists $y\in \ell^2(\V)$ such that
$$(\mathcal I^+)^\top y = P_Y (\mathcal I^+)^\top x, \qquad  (\mathcal I^-)^\top y = P_Y (\mathcal I^-)^\top x,$$
which can equivalently be stated as
$$
\widehat{P_Y} \range \It \subset \range \It.
$$

\smallskip
To see the second equivalence, observe first that the existence of the claimed basis is equivalent to 
$\range \It$ being decomposable into 
$\range \It=(\Ker \widehat{P_Y} \cap \range \It ) \bigoplus (\range \widehat{P_Y} \cap \range \It).$ 
Now we can apply Lemma \ref{decomp}.
\end{proof}

The conditions in Theorem are \ref{charadmiss} are difficult to check. In practice one maybe wants to use criteria which are easier.
The next two Lemmas are criteria of this type.

\begin{lemma}\label{graphnecessary}
Consider a connected finite graph $\G$. If $P_{\mathcal Y}$ is admissible, then $\mathbb 1$ is an eigenvector of $P_Y$.
\end{lemma}
\begin{proof}
First observe that $\mathbb 1 \in V$, since it is differentiable and continuous on the graph. 
Since the continuity conditions on the nodes must be satisfied, 
there exists $\lambda \in \mathbb C$ such that $P_Y \mathbb 1 = \lambda \mathbb 1$, i.e., $P_Y\mathbb 1 = \lambda \mathbb 1$.
So, $\mathbb 1$ is an eigenvector.
\end{proof}

\begin{lemma}\label{disjointdec}
Consider a subset $\E' \subset \E$ of the edge set, 
the induced graph $G':=(\emptyset, E')$,
and a non-admissible orthogonal projection $\mathcal P_{Y'}$ on $\ell^2(E')$.
Then the projection $\mathcal P_Y$ of $\ell^2(\E)$ onto $Y'\bigoplus \ell^2(\E \setminus E')$ is given by
\[
P_Y:=\begin{pmatrix}
P_{Y'}&0\\
0&\id
\end{pmatrix}
\]
and is not admissible.
\end{lemma}
\begin{proof}
Define $V':=\bigoplus_{\me \in \E'} H^1(0,1)$ with continuity node conditions in the nodes of $\G'$.
Since $\mathcal P_{Y'}$ is not admissible, there exists a function $\psi \in V'$ such that $\mathcal P_{Y'}\psi \not\in V'$, 
i.e., such that the continuity condition is violated in a node $\mv_{k_0}$. 
It is possible to extend the function $\psi$ to a function $\tilde{\psi}$ on the whole graph, such that $d^{\tilde{\psi}}=0$ in all nodes of $\G \setminus \G'$. 
Then the function $\mathcal P_Y\tilde{\psi}$ does not satisfy the continuity condition in $\mv_{k_0}$, either.
\end{proof}

%% file: sec_infinitesymmetries.tex
\section{Irreducibility in infinite networks}\label{sec:infinitesymmetries}
It is well known that the heat semigroup on a finite, connected graph (as well as on smooth domains) is positive and irreducible.
A consequence of Proposition~\ref{infinitedirichlet} is that the same does not need to hold for infinite graphs that are not locally finite.

Roughly speaking, since on nodes with infinite degree Dirichlet boundary conditions are automatically imposed, 
initial data localised on a side of such nodes cannot propagate to the other side. 
Thus, irreducibility is reduced to the purely combinatorial problem whether the graph is still connected if such nodes are ``cut out" of the graph.

In order to characterise this property precisely, we need some auxiliary results, which we shortly discuss now, for the sake of readability.
Lemma~\ref{hilflemma3} states that Dirichlet boundary conditions are only imposed on nodes with infinite degree; 
Proposition~\ref{infiniteinvariance} deals with invariance of dual decompositions; 
in Proposition~\ref{connected} we show that a graph is path-connected if and only if it is topologically connected;
Definition~\ref{zusatzgraph} sets up a language for the main theorem, which is stated in Theorem~\ref{infiniteirreducibility}.
The Lemmas~\ref{hilflemma2}--\ref{hilflemma1} and Proposition~\ref{decompositiongraph} are the last steps to the proof.

\begin{lemma}\label{hilflemma3}
The following implication holds for each graph $\G$:
$$[\forall \mv \in \V, \psi \in V: \pi_{\mv}(d^\psi)=0] \Rightarrow [\deg(\mv)=\infty].$$
\end{lemma}

\begin{proof}
We prove the equivalent statement
$$
[\deg(\mv)< \infty] \Rightarrow [\exists \mv \in \V, \psi \in V: \pi_{\mv}(d^\psi)\neq 0].$$
To see this, choose $\mathbb C \ni \lambda \neq 0$ and set
$$
\psi(\mv')=
\begin{cases}
\lambda, & \mv' =\mv,\\
0, & \mbox{otherwise}.
\end{cases}
$$
For all $x \in \G \setminus \V$ interpolate $\psi$ by affine functions. Then $\psi_{\me} \in H^1(0,1)$ for all $\me \in \E$ and moreover
$$
\|\psi\|^2_{L^2}= \deg(\mv) \frac{|\lambda^2|}{3},\quad \|\psi\|^2_{H^1}= \deg(\mv) \frac{4 |\lambda|^2}{3}.
$$
Finally, $\psi$ is continuous in the nodes. So, $\psi \in V$ has the claimed properties.
\end{proof}

\begin{prop}\label{infiniteinvariance}
Consider the form $(a,V)$ associated with the Laplacian on an graph $\G$ with Kirchhoff node conditions, i.e., 
$V$ is defined in~\eqref{formdomain} and in the Definition~\eqref{networkform} we choose $c_{ij}=\mathbb 1 \otimes \delta_{ij}$ and $M=0$. 
Fix a nodal decomposition of the graph $\G=\G_1 \cup \G_2$. 
Then the following assertions are equivalent.
\begin{enumerate}
\item The ideal $I_1:=L^2(\G_1)$ is invariant.
\item If $\mv_1 \in \G_1,\mv_2 \in \G_2$ are adjacent, then $ \deg(\mv_2) = \infty$.
\end{enumerate}
\end{prop}
\begin{proof}
To see that b) implies a), we have to prove that the conditions in Theorem \ref{ortho} hold. 
We only prove that $P_{I_1} V \subset V$. The second condition is clear, since $P_{I_1} \psi$ and $(\id-P_{I_1}) \psi$ have disjoint support.

First observe that $P_{I_1}$ acts on a function $\psi$ by
$$
P_{I_1}\psi(x)=
\begin{cases}
\psi(x) & x \in \G_1, \\
0 & x \in \G_2.
\end{cases}
$$
So, we only have to prove that $P_{I_1}\psi$ is continuous in each node of $\G_2$. Fix an arbitrary node $\mv_2 \in \G_2$. 
If $\mv_2$ is not adjacent to a node in $\G_1$, then $P_{I_1}\psi(x)=0$ for all $x$ in a neighbourhood of $\mv_2$ and so it is continuous.
If $ \mv_2$ is adjacent to a node on $\G_1$, then $\deg(\mv_2)=\infty$ and thus $\psi(\mv_2)=0$ by Proposition~\ref{infinitedirichlet}. 
As a consequence, also $P_{I_1}\psi$ is continuous in $\mv_2$.

To prove the converse implication, observe that the boundary space $\partial V \subset \ell^2(V)$ satisfies
$$
\partial V:=\{ d^\psi: \psi \in V\} \subset \{(x_{\mv})_{\mv \in \V} \in \ell^2(\V): [\deg(\mv')=  \infty \Rightarrow x_{\mv'}=0]\}.
$$
Assume that $L^2(\G_1)$ is invariant. By Theorem \ref{ouh} $P_{I_1}\psi$ is continuous in the nodes if $\psi$ is continuous in the nodes.
In particular, $P_{I_1} \psi$ has to be continuous in all nodes $\mv_2 \in \G_2$ which are adjacent to $\G_1$.
For those nodes $P_{I_1}\psi(\mv_2)=0$. 
Fix now an arbitrary $\mv_2$ of this kind, a neighbourhood $N$ of $\mv_2$ and an arbitrary $\psi \in V$.
On each point $x \in (N \cap \G_1) \setminus \{\mv_2\}=:N_1$, the projection $P_{I_1}$ acts as the identity. 
Further, the ideal $I_1$ is invariant and so $P_{I_1}$ is a continuous function. We compute
$$
0=\lim_{N_1 \ni x \to \mv_2}P_{I_1} \psi(x)= \lim_{N_1 \ni x \to \mv_2} \psi(x)= \psi(\mv_2).
$$
Since $\psi$ is arbitrary, $\deg(\mv_2)=\infty$ follows from Lemma~\ref{hilflemma3}.
\end{proof}

Proposition \ref{infiniteinvariance} allows us to characterize irreducible semigroups. 
For finite graphs, irreducibility is equivalent to the graph $\G$ being connected in the sense of Definition~\ref{graphdefin}.

Before proving similar results for infinite graphs, we show that the definition of connectedness in Definition~\ref{graphdefin} is equivalent to the topological one.

\begin{prop}\label{connected}
The following assertions are equivalent.
\begin{enumerate}
\item The graph $\G$ is connected in the sense of Definition \ref{graphdefin}, i.e.,
$$d(\mv,\mv')<\infty, \qquad \mv,\mv' \in \V.$$
\item The graph $\G$ is topologically connected, i.e., if $\emptyset \neq \V_1, \V_2 \subset \V$ are sets such that
$$
\V_1 \cap \V_2 = \emptyset, \quad \mbox{and} \quad  \V_1 \cup \V_2 =\V,
$$
then there exists $ \me \in \E$ such that $\me \sim \V_1, \me \sim \V_2$.
\end{enumerate}
\end{prop}

\begin{proof}
Assume that a) holds and fix a decomposition $\V=\V_1 \cup \V_2$. 
Since the graph is connected, for every $\mv_1 \in \V_1, \mv_2 \in \V_2$ there exists a path $P$ of finite length connecting $\mv_1$ to $\mv_2$.
In other words, $P$ does not lie entirely in either $\G_1$ or $\G_2$. In particular, there exists an edge $\me \in P$ connecting $\G_1$ to $\G_2$.

Assume that a) does not hold. Fix an arbitrary $\mv\in \V$ and define
$$
\V_1:= \bigcup_{n=1}^\infty \{\mv' \in \V: d(\mv,\mv')=n\}.
$$
Obviously, $\V_1$ is connected subgraph of $\V$. Since a) does not hold, $\V \neq \V_1$. Denote $\V_2:=\V \setminus \V_1$.
We observe that the subgraphs induced by $\V_1$ and $\V_2$ are disjoint.

In fact, assume that there exists $\me \in \E$ such that $\me(0) \in \V_1,\me(1) \in \V_2$. 
Then $d(\mv ,\me(1)) < \infty$ for all $\mv \in \V$ and so $\me(1) \in \V_1$, which is a contradiction.
Hence the subgraphs induced by $\V_1,\V_2$ are disjoint and so the graph $\G$ is not topologically connected.
\end{proof}

In order to to characterise irreducibility, we need some additional definitions.
\begin{defin}\label{zusatzgraph}
During the rest of this section, we use the following definitions.
\begin{itemize}
\item The generalized distance function $\dv: \V \times \V \to \mathbb N \cup \{ \infty\}$ is defined by
$$
\dv(\mv,\mv')=\min_{P \in P[\mv,\mv']} \sum_{\mv'' \in P} \deg (\mv'').
$$
\item The generalised distance function $\de$ between two edges is defined analogously.
\item The finite span $\Sfin(\me) \subset \E$ of $\me$  is the defined by
$$
\Sfin(\me):=\{\me' \in \E: \de(\me,\me') < \infty\}.
$$
We denote, if there is no danger of confusion by $\Sfin(\me)$ the finite span on the edge set as well as the induced subgraph.
\end{itemize}
\end{defin}
We are now in the position of stating the main theorem.
\begin{thm}\label{infiniteirreducibility}
The following assertions are equivalent.
\begin{enumerate}
\item The semigroup $e^{ta}$ is irreducible. 
\item $\Sfin(\me)= \G$ for one $\me \in \E$.
\item $\Sfin(\me)= \G$ for all $\me \in \E$.
\end{enumerate}
\end{thm}
\begin{cor}
If $\G$ is a connected, locally finite network, then $\etasg$ is irreducible.
\end{cor}
\begin{rem}
The functions $d,\dv,\de$ are not distances in the usual sense. 

In fact, $d$ is not defined everywhere if the graph is not connected, and $d_{\V}, d_{\E}$ if are not defined everywhere if the graph is not locally finite.
\end{rem}
For the sake of the simplicity, we split the proof in several lemmas. 
The idea is to prove that the invariant ideals of the semigroup $\etasg$ are all of the form 
$\Sfin(\me)$ for some $\me \in \E$. 

As a preliminary remark, we observe that the only possible invariant ideals are of the form $L^2(\G')$, 
where $\G'=(V',\emptyset)$ is some subgraph of $\G$ induced by a subset of the node set.
In fact, all ideals of $L^2(\G)$ have the form $L^2(\omega)$, where $\omega \subset \bigoplus_{\me \in \E}[0,1]$. 
As we already observed, $\omega=\bigoplus_{\me \in \E} \omega_{\me}$. 
Thus, we are claiming that if $L^2(\omega)$ is invariant, then $\omega_{\me}\in \{\emptyset, [0,1]\}$, 
but this is a consequence of the irreducibility of hte heat semigroup on $L^2[0,1]$.

Next, we show that ideals of the form $\Sfin(\me)$ are invariant.

\begin{lemma}\label{hilflemma2}
Consider a connected graph $\G$. Fix $\me \in \E$. Then $\Sfin(\me)$ is invariant under the action of $\etasg$.
\end{lemma}
\begin{proof}
We use Theorem \ref{ouh}. To prove that both condition hold, fix $\me \in \E$ and denote $P$ the projection onto $\Sfin(\me)$. 
Observe that $P\psi(x)=\mathbb 1_{\Sfin(\me)}(x) \psi(x)$ for all $x \in \G$.
So, the first condition of Theorem \ref{ouh} holds since $P\psi$ and $(I-P)\psi$ have disjoint support. 

\smallskip
We prove that $P V\subset V$. 
Recall that the boundary of $\Sfin(\me)$ consists of those nodes that are adjacent to $\Sfin(\me)$ and to its complement, according to Definition~\ref{graphdefin}.
As usual, one only has to prove continuity in the nodes $\mv'\in \partial \Sfin(\me)$.
To see this, arguing as in the proof of Proposition \ref{infiniteinvariance}, it suffices to prove that $\partial \Sfin(\me)$ only consists of nodes of infinite degree.
We split the proof in two steps. First, we assume that the boundary only consist of one node $\mv'$.
For that node we define
$$
d(\me,\mv'):= \min_{\mv'\sim \me' \in \Sfin(\me)} \de(\me,\me').
$$
This number is finite, since $\mv' \in \Sfin(\me)$ and so there is a path of finite generalized length going from $\me$ to an edge $\me'$ which is incident on $\mv'$.
Since the boundary of $\Sfin(\me)$ only consists of $\mv'$, 
all paths from an edge in $\Sfin(\me)$ to an edge in its complement have to pass through $\mv'.$
As a consequence, we compute for all edges $\me' \in \Gamma(\mv') \setminus \Sfin(\me)$
\begin{equation}\label{auxdegree}
\de(\me',\me) = d(\me,\mv') + \deg(\mv').
\end{equation}
Assume that the degree of $\mv'$ is finite. Then $\de(\me',\me)$ is also finite, and $\me'$ belongs to $\Sfin(\me)$ which is a contradiction.

\smallskip
If $\partial \Sfin(\me)$ contains more than a single node, then $\de(\me,\mv')$ equals the minimum of a degree sum,
where the minimum is taken over all possible paths. Since $d(\me,\mv')=\infty$, then each term in the minimum is infinite, and the claim follows. 
We omit the precise technical details of this construction.

So we obtained that exactly the nodes on the boundary of $\Sfin(\me)$ have infinite degree. As a consequence for all $\mv \in \partial \Sfin, \psi \in V$, $\psi(\mv)=0$.
As a consequence $P\psi$ is continuous in the nodes and the proof is complete.
\end{proof}

The next step is to identify the subgraphs of the form $\Sfin(\me)$. Recall that the boundary of a subgraph is introduced in Definition~\ref{graphdefin}.

\begin{lemma}\label{hilflemma1}
Consider a connected graph $\G$ and a connected subgraph $\G'$. The following assertions hold.

\begin{enumerate}
\item For each $\me \in \E$, $\deg(\mv)=\infty$ for all $\mv \in \partial \Sfin(\me)$. 

\item If $\de(\me,\me') <\infty$ for all $\me,\me' \in \G'$ and 
$\deg(\mv)=\infty$ for all $\mv \in \partial \G'$, then
$
\G' = \Sfin(\me)
$
for all $\me \in \G'.$
\item If $\deg(\mv)=\infty$ for all $\mv \in \partial \G'$, then 
$
\Sfin(\me ) \subset \G'
$
for all $\me \in \G'$.
\end{enumerate}
\end{lemma}

\begin{proof}
a). This claim is the first step of proof of the Lemma \ref{hilflemma2}. 

b) Fix a subgraph with the required properties.
First observe that since $\de(\me,\me')<\infty$ for all $\me,\me' \in \G'$, then $\G' \subset \Sfin(\me) $ for all
$\me \in \G'$.

Assume now that $\exists \me' \in \Sfin(\me)\setminus \G'$. 
Without loss of generality, let $\me' \sim \G'$, i.e., $\me'\sim \mv, \mv \in \G'$ and assume that the boundary $\partial \G'$ only consists of $\mv$. 

By hypothesis $\deg(\mv)= \infty$ and so $\de(\me,\me')=\infty$, which is a contradiction.

c) Fix an arbitrary $\me \in \G'$, $\me' \not\in \G'$, and a path $P \in P[\me,\me']$. 
By definition of $\partial \G'$, there exists $\mv \in P \cap \partial \G'$. As a consequence $\de(\me,\me')= \infty$ and the proof is complete.
\end{proof}

The following is a consequence of the above lemma.

\begin{prop}\label{decompositiongraph}
Consider a connected graph $\G$. Then there exists $\E' \subset \E$ such that
$$
\bigcup_{\me \in \E'} \Sfin(\me)= \G,
$$
and
$$ 
\partial \Sfin(\me) = \Sfin(\me)\cap \Sfin(\me')=\partial \Sfin(\me'), \qquad \me,\me' \in \E'.
$$
Moreover, for all $\me \in \E'$, $\deg(\mv)=\infty$ for all $\mv \in \partial \Sfin(\me)$.

The same decomposition holds for each subgraph of $\G$ whose boundary only consists of nodes with infinite degree.
\end{prop}
Now we can characterise irreducibility.
\begin{proof}[Proof of Theorem \ref{infiniteirreducibility}]
Observe that, as a consequence of Lemma \ref{hilflemma3}, $\psi(\mv)=0$ for all $\psi \in \V$ if and only if $\deg(\mv)=\infty$.

By definition, a semigroup is irreducible if there are no non-trivial invariant ideals. 
Thus, it suffices to prove that all invariant ideals of $\etasg$ have the form $L^2(\Sfin(\me))$, i.e.,
$$
\{I \subset L^2(\G): I \mbox{ is an invariant ideal}\}=\{L^2(\Sfin(\me)): \me \in \E\}.
$$
By Lemma \ref{hilflemma2}, $L^2(\Sfin(\me))$ are invariant and so ``$\supset$'' holds.

Consider $I$ invariant ideal of $\etasg$. As we have already observed, $I=L^2(\G')$ for some subgraph $\G' \subset \G$. 

The projection $P\psi$ of a function $\psi$ on $L^2(\G')$ vanishes in all points of $\G\setminus\G'$ and so, 
by continuity it vanishes in all points of the boundary of $\G'$. 
As a consequence, each function $\psi \in V$ also has to vanish on all points of the boundary of $\G'$.

By Lemma \ref{hilflemma3}, all those points $\mv'$ satisfy $\deg(\mv')=\infty.$ 
So, there exists a decomposition of $\G'$ as in Lemma \ref{hilflemma1}, i.e.,
$$
L^2(\G')=\bigoplus_{\me \in \E''} L^2(\Sfin(\me)).
$$
Now, $L^2(\G')$ has to be an ideal, so $\E''$ has cardinality one, i.e., $\G'=\Sfin(\me)$ for some $\me$.

Summing up, we have proved that $I$ is an invariant ideal of $\etasg$ if and only if $I=L^2(\Sfin(\me))$ for some $\me \in \E$.

Concluding, irreducibility is then equivalent to $\Sfin(\me)= \G$ for all $\me \in \G.$ 
In order to see that b) and c) also are equivalent observe that if $\me' \in \Sfin(\me)$, then $\Sfin(\me)=\Sfin(\me')$.
\end{proof}

%% file: sec_treelayer.tex
\section{Symmetries in special networks}\label{sec:treelayer}

In this section we will study symmetry issues for some special classes of graphs.
We start presenting some graph theoretical definitions.

\begin{defin}\label{graphdef}
Let $\G$ a graph with no isolated nodes, i.\,e., such that
$\deg(\V)\geq 1$ for all $\mv \in \V$.
\begin{itemize}
\item We call the graph $\G$ \emph{completely unconnected} if $\G$ is the union of disjoint compact intervals, i.\,e., if $\G$ is a regular graph of degree $1$.
\item We call the graph $\G$ an \emph{inbound} (respectively, \emph{outbound}) star, if there exists a node $\V$ such that
$\me(1)=\V$, (respectively, if $\me(0)=\V$), for all $\me \in \E$. 
We call the graph $\G$ a star if it is an inbound or outbound star and $\V$ the \emph{centre} of the star.
\item We call the graph $\G$ \emph{simple} if there are no parallel edges, i.e., if $\partial$ is an injective mapping.
\item We call the graph $\G$ \emph{bipartite} if each node has only either ingoing or outgoing edges.
\item We call the graph $\G$ \emph{Eulerian} if all nodes have the same number of ingoing and outgoing edges.
\item We call a graph $\G$ a \emph{layer graph} if there  exist disjoint sets $V_1,\ldots,V_L$ such that 
\begin{itemize}
\item $V=\cup_{p=1}^L V_p$,
\item $\me(0) \in V_p$ implies $\me(1) \in V_{p+1}$  for all $p=1,\ldots,L-1$, and 
\item $\me(0) \in V_L$ implies $\me(1) \in V_{1}$.
\end{itemize}
Nodes belonging to $V_p$ are said to lie in the $p$\textsuperscript{th} layer. 
Edges outgoing from nodes in the $p$\textsuperscript{th} layer are also said to lie in the $p$\textsuperscript{th} layer.
\item We call a layer graph \emph{symmetric} if the ingoing and outgoing degrees of the nodes only depends on the layer, 
i.e., if there exist numbers $I(p), O(p) \in \mathbb N_0$ such that $\deg^+(\V)=I(p), \deg^-(\V)=O(p)$ for all nodes $\V$ in the $p$\textsuperscript{th} layer.
\end{itemize}
\end{defin}
During this section, it is always assumed that the graph $\G$ is finite, and that numberings of the edges and of the vertices are fixed, i.e.,
$$
\E:= (\me_1,\ldots,\me_m) \quad \mbox{and} \quad \V:=(\mv_1,\ldots,\mv_n)
$$

\subsubsection{Bipartite and Euler Graphs}
It is possible to characterize some classes of graphs by the
admissibility of the projection $P_Y$ on $Y= \langle \mathbb 1 \rangle$.
\begin{thm}\label{mittelunggraph}
Consider a finite graph $\G$ and the orthogonal projection $P_Y$ of $\ell^2(\E)$ onto $Y\langle \mathbb 1 \rangle $.
Then $P_Y$ is admissible if and only if $\G$ is bipartite or Eulerian.
\end{thm}

\begin{proof}
As usual, we have to prove that for a continuous function $\psi \in \V$ also $P_\mathcal Y \psi$ is continuous in the nodes.
Observe that the projection $P_Y$ of $\ell^2(\E)$ onto $Y$ satisfies
$$
P_Y x = \sum_{i =1}^m \frac{x_i}{m}.
$$
Let $\V_1 \subset \V$ denote the set of all vertices
having outgoing edges, and let $\V_2 \subset \V$
denote the set of all vertices having ingoing edges. 
We distinguish
two cases. First, assume  $\V_1 \cap \V_2 = \emptyset$.
Then $\G$ is a bipartite graph.

On the other hand, if $\V_1 \cap \V_2  \not= \emptyset$, then
$d^{{P}_\mathcal Y \psi}$ exists if and only if
\begin{equation} \label{sameaverage}
\sum_{j=1}^m\frac{\psi_j(0)}{m} = \sum_{j=1}^m\frac{\psi_j(1)}{m}.
\end{equation}
We show now that the equality
\eqref{sameaverage} is equivalent to the graph being Eulerian.
First, assume that \eqref{sameaverage} holds for every $\psi \in V$.
Fix an arbitrary $\mv_k \in \V$ and choose $\psi \in V$ such
that $d^\psi = \mathbb 1_{\{i\}}$. Then
$$
\frac{1}{m}  \deg^+(\mv_k) = \sum_{j=1}^m\frac{ \psi_j(0)}{m} =
\sum_{j=1}^m\frac{\psi_j(1)}{m} = \frac{1}{m} \deg^-(\mv_k).
$$
Thus it is necessary that $\deg^-(\mv_k) = \deg^+(\mv_k)$ holds
for every $k=1,\ldots,n$. 

\smallskip
Conversely, assume that
$$
\deg^-(\mv_k) = \deg^+(\mv_k), \qquad k=1,\ldots,n.
$$
Then
\begin{eqnarray*}
\sum_{j=1}^m\frac{\psi_j(0)}{m} & = & \frac{1}{m}\sum_{k=1}^{n}\deg^+(\mv_k)d^\psi_k\\
& = &\frac{1}{m}\sum_{k=1}^{n}\deg^-(\mv_k)d^\psi_k\\
& = & \sum_{j=1}^m\frac{\psi_j(1)}{m}.
\end{eqnarray*}
Hence \eqref{sameaverage} is satisfied, so this condition is also sufficient.

\smallskip
It only remains to show that indeed for every bipartite graph $P_Y$ is admissible.
To see this, note that for an arbitrary $\psi \in V$ the vector $d^{P_ \mathcal Y \psi}$ can
be chosen to equal $\sum_{i=1}^m\frac{\psi_i(0)}{m}$ in all components belonging to nodes in
$\V_1$ and to equal $\sum_{i=1}^m\frac{\psi_i(1)}{m}$ in all
components belonging to $\V_2$. 
This shows continuity of $P_\mathcal Y \psi$ in the nodes,
thus implying $ P_\mathcal Y \psi \in V$.
\end{proof}

\subsubsection{Stars}
The main result of this subsection is a characterization of stars in the class of the simple graphs.
We first investigate the admissibility of projections.

\begin{prop}\label{unconnstar}
The following assertions hold.
\begin{enumerate}
\item The graph $\G$ is completely unconnected if and only if
$P_\mathcal Y$ is admissible for all linear subspaces $Y \subset \ell^2(E)$.
\item Let $\G$ be a simple, connected graph. 
Then $\G$ is a star if and only if $P_\mathcal Y$ is admissible for all linear subspaces $Y \subset \ell^2(E)$ such that
$\langle \mathbb 1 \rangle \subset Y.$
\end{enumerate}
\end{prop}

\begin{proof}
a) Since the graph $\G$ is completely unconnected, the
continuity condition in $V$ is empty, and therefore each $P_\mathcal Y$ is admissible. 
Conversely, if $\G$ is not completely unconnected, then it is possible to decompose $\G$
into the disjoint union of a connected graph $\G_1$ with $m_1$
edges, $m_1\geq 2$ and the remaining graph $\G_2$. 
Let $K_1$ be an orthogonal projection of $\mathbb C^{m_1}$, which does not have $\mathbb 1$ as
an eigenvector. 
Lemma~\ref{graphnecessary} asserts that the orthogonal projection
$$\begin{pmatrix} K_1&0\\0&\id \end{pmatrix}$$ 
is not admissible.

b) Let the graph $\G$ be a star and $Y$ a subspace containing $\langle \mathbb 1 \rangle$.
Without loss of generality, we prove the claim for an outgoing
star with centre $\V_1$ and with the natural numbering of the other
nodes. 

In fact, for this star
$$
\It=\begin{pmatrix}\mathbb 1&0\\0&\id_m\end{pmatrix}.
$$ 
Since now $P_Y$ has
$\mathbb 1$ as eigenvector to the eigenvalue 1, 
$$
\begin{pmatrix}  P_Y \mathbb 1 &0 \\0 & P_Y \id_m\end{pmatrix} =
\begin{pmatrix}  \mathbb 1 &0 \\0 & P_Y \end{pmatrix}.
$$
So, $\range\widehat{P_Y }\It \subset \range\It $, and this
implies the admissibility of $\mathcal {P}_Y$. 

\smallskip
Conversely, assume that the graph $\G$ is not a star. This implies
the existence of an undirected path of length
$3$. We will denote it by
$\me_1,\me_2,\me_3$, possibly relabelling the edges.
Our strategy is the following: for each graph that is a path consisting of $3$ edges we construct a non-admissible projection $P_L$ where $L\mathbb 1 = \mathbb 1$.
We then consider the projection $\mathcal P_Y$, where $P_Y$ is
\[
P_Y:=\begin{pmatrix}
P_L & 0\\0&\id
  \end{pmatrix}.
\]
Then, by Lemma~\ref{disjointdec}, we conclude that $\mathcal P_Y$ is not admissible,
although $\mathbb 1$ is an eigenvector of $P_Y$.

\smallskip
First, consider cycles of length $3$. Since each edge can be
directed arbitrarily, there are $8$ such graphs. Let us
start with the case of a not strongly connected graph.

Such graphs are neither Eulerian nor bipartite.
Thus, Theorem~\ref{mittelunggraph} provides an example of an $L$ as requested.
If the graph is a (directed) cycle such that
$\me_1(0)=\mv_1$, consider the projection

\[ L:=\begin{pmatrix}
\frac{1}{2} & \frac{1}{2} & 0\\
\frac{1}{2} & \frac{1}{2} & 0\\
  0&0&1
\end{pmatrix}\]
and the function $\psi$ defined by $\psi(x):=(x,1-x,0)^\top \in V$. Then, $\psi \in V$ but
$P_\mathcal Yf \not\in V$, since $ P_Y f(x) = (\frac{1}{2}, \frac{1}{2},0)^\top$ for a.\,e. $x\in(0,1)$.

Consider now the lines of length $3$. We split this into three possible
cases: $\G$ may be bipartite line, a (directed) line, or neither a (directed) line nor a bipartite graph.
In the last two cases the graphs is neither bipartite nor
Eulerian, and hence we can use Theorem~\ref{mittelunggraph} again.

In the case of a bipartite line, let us consider the projection
\[
L:=\begin{pmatrix}
\frac{1}{2} & \frac{1}{2} & 0\\
0&0&1\\
\frac{1}{2} & \frac{1}{2} & 0\\
\end{pmatrix}
\]
for the parametrization
$\me_1(0)=\mv_1$, $\me_1(1)=\me_2(1)=\mv_2$, $\me_2(0)=\me_3(0)=\mv_3$, and
$\me_3(1)=\mv_4$. 
Consider the function $\psi(x):=(x,x,0)^\top$. 
Again, $\psi \in V$ but  $ P_Y \psi \not\in V$, 
since $P_Y \psi(x)= (\frac{x}{2}, x,\frac{x}{2})^\top$ for a.\,e. $x\in(0,1)$.
\end{proof}

\begin{rem}
In Proposition \ref{unconnstar}, (2) we have assumed the graph $\G$ to have no multiple edges. 
In fact, it is not possible to relax this
condition, since all orthogonal projections with eigenvector $\mathbb 1$ are admissible on all connected graphs consisting of $2$ nodes and
$m$ edges for each $m \in \mathbb N$ and each orientation of the edges.
\end{rem}

\subsubsection{Layer Graphs}
In this section we prove an admissibility result for symmetric layer graphs.
We start fixing a canonical numbering of the edges of a layer graph.
First observe that the node decomposition induces an edge
decomposition $\E=\bigcup_{p=1}^L \E_p$ by setting
$$
\E_p:=\{\me \in\E : \me \mbox{ lies in the $p$\textsuperscript{th} layer}\}.
$$
After relabelling the edges we may assume that there exist $L_p$,
$p=1,\ldots,L+1$ satisfying
\begin{enumerate}
\item $L_1=0$;
\item $\me_i(0)=\me_j(0) \mbox{ or } \me_i(1)=\me_j(1)$ implies $L_{p-1} < i,j \leq L_{p}$ for some $p$;
\item $\me_i(0)=\me_j(1)$ implies $L_{p-1} < j \leq L_p < i \leq L_{p+1}$ for some $p$.
\end{enumerate}
The numbering obtained in such a way has the property that $\me_i$ is
in the $p$\textsuperscript{th} layer if and only if $L_p < i \leq L_{p+1}$. In fact,
all edges $\me_i$ such that $i \leq L_{p+1}$ are in any of the first
$p$ layers.

We are going to exhibit a class of admissible projections.

\begin{prop}\label{layer}
Consider a symmetric layer graph $\G$ and the orthogonal projection $P_Y$
\begin{equation}\label{layerproj}
P_Y=
\begin{pmatrix}
(\frac{1}{|E_1|})_{i,j=1,\ldots,|E_1|} &&0\\
&\ddots&\\
0&&(\frac{1}{|E_L|})_{i,j=1,\ldots,|E_L|}
\end{pmatrix},
\end{equation}
where $|E_p|$, $p=1,\ldots,L$ denotes the number of edges in the
$p^{\rm th}$ layer. 
Then $\mathcal P_K$ is admissible.
\end{prop}

\begin{proof}
One has to check the continuity condition for each $p=1,\ldots, L-1$ in every node of the $p$\textsuperscript{th} layer. 
Define the auxiliary function
$$
\lambda: k \mapsto \mbox{layer of the node $\mv_k$}.
$$
We thus have to check continuity in those nodes $\mv_k$ such that $\lambda(k) = p$, $p=1,\ldots, L-1$.

The set $\lambda^{-1}(p)$ can be represented in the form
$$
\lambda^{-1}(p)=\{k: \exists i \in \{L_p+1,\ldots,L_{p+1}\} \mbox{ s.\,t.
}\me_i(1)=\mv_k\},
$$
as well as in the form
$$
\lambda^{-1}(p)=\{k: \exists i \in \{L_{p+1}+1,\ldots,L_{p+2}\} \mbox{
s.\,t. }\me_i(0)=\mv_k\},
$$
whenever the expression is defined.
By the definition of $P_Y$, for all $p=1,\ldots,L-1$ and
all $i,j=L_p+1,\ldots,L_{p+1}$ the identities
\begin{equation}\label{indimitt1}
\mathcal P_Y\psi_i(1)=P_Y \psi_j(1), \qquad
\mathcal P_Y \psi_i(0)=P_Y\psi_j(0).
\end{equation}
hold.
As a consequence, for layers having ingoing or outgoing degree $0$,
the continuity is obvious. Assume now that $I(p)\neq 0$ and
$O(p)\neq 0$.

For an edge $\me_i$ in the $p$\textsuperscript{th} layer and for $\psi \in \mathcal V$,
\begin{eqnarray*}
\mathcal (P_Y \psi)_i(1)=\sum_{i=L_p+1}^{L_{p+1}} \frac{ \psi_i(1)}{|E_p|}
=\sum_{k\in \lambda^{-1}(p)} \deg^-(\mv_k)\frac{ \psi_i(1)}{|E_p|}.
\end{eqnarray*}
Recall that since our graph is symmetric, the incidence degree
$\deg^-(\mv_k)$ only depends on the layer, and therefore we can write
$$
\mathcal (P_Y \psi)_i(1)=\sum_{k\in \lambda^{-1}(p)} |I(p)|
\frac{\psi(\mv_k)}{|E_p|}
=\frac{|I(p)|}{|E_p|}\sum_{k\in \lambda^{-1}(p)} f(\mv_k).
$$
With analogous computations we obtain for edges in the $p+1$ layer
$$
\mathcal (P_Y \psi)_i(0)=\frac{|O(p)|}{|E_{p+1}|}\sum_{k\in
\lambda^{-1}(p)} \psi(\mv_k).
$$
Observe that the identities
$|E_{p+1}|=|\lambda^{-1}(p)||O(p)|$ and $|E_p|=|\lambda^{-1}(p)||I(p)|$
imply
${|O(p)|}{|E_{p+1}|}^{-1}={|I(p)|}{|E_p|}^{-1}$.
We have thus proved that $
(P_Y \psi)_i(1)=(P_Y \psi)_j(0)$ for all $\me_i,\me_j$ such that $i \in
\lambda^{-1}(p), j\in \lambda^{-1}(p+1)$. This completes the proof.
\end{proof}

\begin{rems}\mbox{}
\begin{enumerate}[(1)]
\item  	The above result can be restated as follows. 
	Let $\G$ be a symmetric layer graph.  Then, if $M=0$ and $C=\id$,
	$$
	\mathcal Y:=\{\psi \in L^2: \psi_i=\psi_j \mbox{ for all } i,j \in \ell^{-1}(p), \; p=1,\ldots,L\}
	$$
	is invariant under the action of $(e^{tA})_{t\geq 0}$. In fact, this means that ``radial'' function are invariant not only for trees but also for layer graphs.

\item 	The class of the \emph{layer graphs} is a common object in the graph theoretical literature. 
	In fact, layer graphs are nothing but (directed) $p$-partite
	graphs, for which collapsing the components of the graphs to
	a single vertex leads to a finite line or to a cycle.  In particular, homogeneous trees of finite depth are symmetric
	layer graphs. Such graphs play a role in the investigation of
	biological neural networks.
\item	The symmetry condition in Proposition~\ref{layer} cannot be
	relaxed. To see this, consider the following simple example. Let $\G$
	be an outbound star of order two and consider two copies of $\G$.
	Identifying two of the external nodes defines a layer graph. It is possible to 
	show that the orthogonal projection defined in~\eqref{layerproj} is not admissible, due to the two free nodes in the second layer.
\end{enumerate}
\end{rems}

%% file: sec_general.tex
\section{Generalised networks}\label{sec:general}
We want to discuss in this section an alternative setting which contains network equations as a special case.
Therefore, we only assume that $\E$ is a countable set. 
We stress that we are not in graph theoretic setting, and so $\E$ should not be considered as the edge set of a graph.
We start with a definition justified by Lemma~\ref{inclusionboundaryspace}.
\begin{defin}
The bounded operators $\partial_0, \partial_1 \in \mathcal L (H^1((0,1), \ell^2(\E)) $ are defined by
$$
\partial_0 \psi:= \psi(0), \quad \partial_1\psi:= \psi(1), \qquad \psi \in H^1((0,1), \ell^2(\E)).
$$
Analogously, the operator $\partial \in \mathcal L (H^1(0,1), \ell^2(\E) \bigoplus \ell^2(\E) )$ is defined by
$$
\partial \psi:= \begin{pmatrix} \psi(0)\\ \psi(1) \end{pmatrix}, \qquad \psi \in H^1((0,1), \ell^2(\E)).
$$
\end{defin}

Fix a subspace $Y$ of $\ell^2(\E)\bigoplus \ell^2(\E)$. We define a form domain by setting
\begin{equation}\label{gendomain}
V_Y:=\{\psi \in H^1((0,1),\ell^2(\E)): \partial \psi \in Y\}.
\end{equation}
In analogy with~\eqref{networkform}, we define the form $a:V_Y \times V_Y \to \mathbb C$ by
\begin{equation}\label{genform}
a(\psi,\psi'):=( \frac{d}{dx}\psi \mid \frac{d}{dx}\psi') - (M \partial \psi \mid \partial \psi'), \qquad \psi,\psi' \in V_Y.
\end{equation}
Here $M \in \mathcal L(\ell^2(\E \bigoplus \E))$. Exactly as in the sections before, the form $(a,V_Y)$ has domain dense in $L^2((0,1), \ell^2(\E))$. 
Further, it is possible to characterise continuity and ellipticity of the form.

We focus on the identification of \emph{nodes} of the form domains. 
Observe that, since we are not in a graph theoretic setting, there are no nodes in the sense we discussed in the previous sections.

The idea for doing this is the following: 
the condition $\partial \psi \in Y$ in the definition of the form domain will contain informations about couplings in the boundary values of $\psi$.
We identify several boundary points with a node, if the boundary condition of two different nodes are mutual independent, 
and if they are minimal in some sense.

This identification relies on the following lemma.
\begin{lemma}
Fix a countable index set $I$ and consider a linear bounded operator $A \in \mathcal L(\ell^2(I))$.
Then there exists an index set $\V$ and an uniquely determined family $(I_{\mv})_{\mv \in \V}$ with the following properties:
\begin{enumerate}[a)]
\item $(I_{\mv})_{\mv\in \V}$ is a partition of $I$;
\item For all $\mv \in \V$, $\ell^2(I_{\mv})$ is invariant under the action of $A$ and $A^*$;
\item For all $I' \subset I$, $I' \neq \emptyset$, if $\ell^2(I')$ is invariant under the action of $A$ and $A^*$, 
then there exists $\mv' \in \V$ such that $I_{\mv } \subset I'$.
\end{enumerate}
We call the $I_{\mv}$ the \emph{invariant blocks} of $A$.
\end{lemma}
\begin{proof}
We first observe that the statement can be reformulated in the language of infinite matrices. 
In fact, if $A$ is an infinite matrix (not necessarily countable), i.e., a doubly indexed family $(a_{ij})_{i,j \in I}$, then the statement is equivalent to the existence
of an index set $\V$ and an uniquely determined family $(I_{\mv})_{\mv \in \V}$ with the following properties:
\begin{enumerate}[1)]
\item $(I_{\mv})_{\mv\in \V}$ is a partition of $I$;
\item For all $\mv \in \V, i \in I_{\mv}, j \in I \setminus I_{\mv}$, $a_{ij}=a_{ij}=0$;
\item For all $ I' \subset I$, $I' \neq \emptyset$, if b) holds for $I'$, 
then there exists $\mv' \in \V$ such that $I_{\mv } \subset I'$.
\end{enumerate}

\smallskip
Without loss of generality, we assume that $a_{ij}\in \{0,1\}$. If this is not the case, we observe that the matrix $B:=(b_{ij})_{i,j \in I}$ obtained by setting
$$
b_{ij}:=
\begin{cases}
1, &  a_{ij}\neq 0,\\
0, & \mbox{otherwise},
\end{cases}
$$
has the same invariant blocks as $A$ and satisfies the assumption.

\smallskip
Now, we interpret $I$ as the vertex set of a graph $\G$ and $A$ as the adjacency matrix,
and we consider the partition $I=(I_{\mv})_{\mv \in \V}$ in  connected components.
We claim that $(I_{\mv})_{\mv \in \V}$ has the claimed properties.

First, $(I_{\mv})_{\mv \in \V}$ is a partition of $I$ by construction.

Second, fix $\mv \in \V$, $i \in I_{\mv}, j \in I \setminus I_{\mv}$. 
Since $i,j$ are not in the same connected component, there is no edge connecting the two nodes. Thus, $a_{ij}=a_{ji}=0$.

Finally, assume that $I' \subset I$ has the property 2). 
Consider the graph with index set $I'$ and adjacency matrix $A_{|I'}$, and decompose it in its connected components.
Since $I'$ is disconnected from the rest of the graph, these connected components are also connected components of the graph $\G$.

\smallskip
Conversely, assume that a family $(I_{\mv})_{\mv \in \V}$ has the properties 1)--3). We have to show that $I_{\mv}$ are the connected components.
First observe that by 2) each of the node sets $I_{\mv}$ is disconnected from the rest of the graph.
Thus, subsets of $I_{\mv}$ are connected in the graph $I_{\mv}$ with adjacency matrix $A_{|I_{\mv}}$ if and only if they are connected in $\G$.
Consider now the decomposition $I_{\mv}:=(I_{\mv_{\ell}})_{\ell \in L}$ in connected components. Observe that 2) holds for all $I_{\mv_{\ell}}$.
By assumption 3), for all $\ell$ there exists $I_{k} \subset I_{\mv_{\ell}} \subset I_{\mv}$. 
Since $(I_{\mv})_{\mv \in \V}$ is a partition by assumption 1), $I_{\mv}=I_{\mv_\ell}$ for all $\ell$. As a consequence $I_{\mv}$ is connected.
\end{proof}

Observe that $\ell^2(\E) \bigoplus \ell^2(\E) \simeq \ell^2(\E\bigoplus \E)$, thus the above Lemma can be applied in our situation.
\begin{defin}
Consider the form $(a,V_Y)$ defined in~\eqref{gendomain}--\eqref{genform} and denote by $P_Y$ the orthogonal projection of $\ell^2(E \bigoplus E)$ onto $Y$.
The \emph{nodes} of the form domain $V_Y$ are the invariant blocks $I_{\mv}$ of $P_Y$.
\end{defin}
\begin{rem}
If $Y:=\langle \{\mathbb 1_{I_{\mv}}: I_{\mv} \mbox{ finite}\}\rangle \bigoplus \{0\}$ and $M=0$, then $(a,V_Y)$ is the form associated with the 
Laplace operator with Kirchhoff boundary conditions on $L^2(\G)$. 

Here the graph $\G$ has vertex set $\V$ and it is obtained identifying the boundary points in the same connected component $I_{\mv}$.
\end{rem}
Next, we identify the domain of the operator associated with $(a,V_Y)$. To this aim we give the following definition.
\begin{defin}
For all $\psi \in H^2((0,1), \ell^2(\E))$ we define the bounded operator $\frac{\partial}{\partial \nu} \in \mathcal L(H^2((0,1), \ell^2(\E)), \ell^2(\E \bigoplus \E))$ by
\begin{equation}\label{normaldergen}
\frac{\partial}{ \partial \nu} \psi:=
\begin{pmatrix}
- \frac{d}{dx} \psi(0) \\
\frac{d}{dx} \psi (1)
\end{pmatrix}.
\end{equation}
\end{defin}
\begin{prop}
The operator $(A,D(A))$ associated with $(a,V_Y)$. Then, the operator associated with $(a,V_Y)$ is given by
\begin{equation}\label{generaliseddomain}
D(A):=\{\psi \in H^2((0,1), \ell^2(\E)): \big( \frac{\partial}{\partial \nu}  - M \big)\partial \psi \in Y^\perp\},
\end{equation}
and
\begin{equation}\label{genoperator}
A\psi:= \frac{d^2}{dx^2} \psi.
\end{equation}
\end{prop}
\begin{proof}
Denote $(B,D(B))$ the operator associated with $(a,V_Y)$. We first show $D(A) \subset D(B)$. Fix arbitrary $\psi \in D(A), \psi' \in V_Y$ and compute
\begin{eqnarray*}
a(\psi,\psi')&=&\sum_{\me\in \E} \int_0^1 \frac{d}{dx}\psi_{\me}(x)\overline{\frac{d}{dx}\psi_{\me}'(x)}dx 
-(M\partial \psi \mid \partial \psi' )\\
&=&-\sum_{\me \in \E} \frac{d^2}{dx^2} \psi_{\me} (x) \overline{\psi'_{\me}(x)} dx + \sum_{\me \in \E}\big[\frac{d}{dx}\psi_{\me} \overline{\psi'_{\me}} \big]^1_0
-(M\partial \psi \mid \partial \psi' )\\
&=&-\sum_{\me \in \E} \frac{d^2}{dx^2} \psi_{\me} (x) \overline{\psi'_{\me}(x)} dx + \big(\frac{\partial }{\partial \nu}\psi \mid \partial \psi'\big)
-(M\partial \psi \mid \partial \psi' )\\
&=&-\sum_{\me \in \E} \frac{d^2}{dx^2} \psi_{\me} (x) \overline{\psi'_{\me}(x)} dx + \big((\frac{\partial }{\partial \nu} - M )\psi \mid \partial \psi'\big)\\
\end{eqnarray*}
The last term vanishes, since $\psi \in D(A)$. Thus $a(\psi,\psi')=(A\psi \mid \psi')$ for arbitrary $\psi' \in V_Y$. 
Choosing $f= A\psi$ shows that $\psi \in D(B)$.

\smallskip
Conversely, fix $\psi \in D(B)$. Then, there exists $f \in L^2((0,1),\ell^2(\E))$ such that $a(\psi,\psi'):=(f \mid \psi')$ for all $\psi' \in V_Y$.
Integrating by part shows that 
$$
-\sum_{\me \in \E} \frac{d^2}{dx^2} \psi_{\me} (x) \overline{\psi'_{\me}(x)} dx + \big((\frac{\partial }{\partial \nu} - M )\psi \mid \partial \psi'\big)
 = (f \mid \psi')
$$
holds for all $\psi \in V_Y$ and thus in particular for $\psi \in H^1_0((0,1), \ell^2(\E))$. 
As a consequence $f:=\frac{d^2}{dx^2}\psi$ and 
$$
\big((\frac{\partial }{\partial \nu} - M )\psi \mid \partial \psi'\big)= 0, \qquad \psi' \in V_Y.
$$
The latter is equivalent to $(\frac{\partial}{\partial \nu}-M)\partial \psi \in Y^\perp$ since $\partial \psi'$ can assume all values in $ Y$.
This means $\psi \in D(A)$ and the proof is complete.
\end{proof}

It is in principle possible to systematically analyse qualitative properties of the heat semigroup associated with the form $a_Y$.
We refer to Section~\ref{sec:histremII} for some remarks about this topic.

%% file: sec_histremII.tex
\section{Discussion and remarks}\label{sec:histremII}
The main topic of our work is the investigation of \emph{symmetry properties} of network equations.
We shortly discuss how our use of the term symmetry is related to the physical meaning of the word.

In a physical context, it is often useful to work in a Lagrangian framework.
This means, a Lagrange function is given for our system, and the solutions are the orbits that minimise this Lagrangian function.

To be more precise, fix a smooth, bounded domain $\Omega$ and consider the space $V:=H^1(\Omega)$. Choose now a positive time $t >0$ and consider the space
$
X:=C^1([0,T], V).
$
Now, the solutions of the heat equation
$$
\left\{
\begin{array}{rclr}
\dot{\psi}(t,x)&=& \Delta \psi(t,x)& t\in [ 0,T], x\in \Omega, \\
\psi(0,x)& =& f(x) & f \in C^\infty(\Omega), x \in \Omega, 
\end{array}
\right.
$$
are the minima of the Lagrange function $\mathcal L: X \to \mathbb C$
$$
\mathcal L (\psi):= \int_{0}^T \int_\Omega \frac{1}{2}\psi(t,x) \dot{\psi}(t,x) - |\nabla \psi(t,x)|^2 dxdt. 
$$
One sees that the second term in the Lagrangian function is exactly the diagonal term of the sesquilinear form $a:H^1(\Omega) \times H^1(\Omega) \to \mathbb C$ defined by
$$
a(\psi,\psi'):= \int_{\Omega} \nabla\psi(x) \overline{\nabla \psi'(x)} dx, \quad \psi,\psi' \in H^1(\Omega).
$$
Now, assume that $(S_z)_{z \in Z}$ is a family of mappings $S_z : H^1(\Omega) \to H^1(\Omega)$ that is a group for the binary operation $\circ$ defined by 
$$
(S_{z_1} \circ S_{z_2}) \psi := S_{z_1}(S_{z_2} \psi), \qquad z_1,z_2 \in Z, \psi \in H^1(\Omega).
$$
Further, assume that each $S_z$ commutes with the time derivative
and that $a(\psi(x))=a(S_z \psi)$ for all $z \in Z, \psi \in H^1(\Omega)$. 
Thus, for all $z \in Z$, $S_z$ maps solutions of the problem above into other solutions:
$(S_z)_{z \in Z}$ is a \emph{symmetry group} of the physical system.

It is possible to operate a distinction between two main type of symmetries: \emph{space symmetries} and \emph{gauge symmetries}.
We now show that our type of symmetries are gauge symmetries.

Space symmetries are symmetries for which $S_z$ is given by a transformation of the domain $\Omega$. 
This means, we choose a family of mapping $(s_z)_{z \in Z}$ such that $s_z: \Omega \to \Omega$ and $(s_z)_{z \in Z}$ is a group for the composition.
Then, one defines
$$
S_z \psi(t, x):=  \psi(t, s_z x), \qquad z \in Z, \psi \in H^1(\Omega), x \in \Omega, t \in [0,T]. 
$$
Such a transformations obviously commutes with the time derivative. Thus, if $a(S_z \psi)=a(\psi)$ for all $z \in Z, \psi \in H^1(\Omega)$, then $(S_z)_{z \in Z}$ is a space symmetry of the system.

We now turn our attention to the case of gauge symmetries. 
Assume that $V:=H^1(\Omega, H)$, where $H$ is a complex Hilbert space (in the easiest case, $H=\mathbb C)$.
Fix a group of mappings $(s_z)_{z \in Z}$, where $s_z \in \mathcal L(H)$ and define
$$
S_z \psi (t, x) := s_z \psi(t, x), \qquad z \in Z, \psi \in H^1(\Omega, H), x \in \Omega, t \in [0,T].
$$
Then $(S_z)_{z \in Z}$ is also a group for the composition, and it also commutes with the time derivative. 
If $a(S_z \psi)=a(\psi)$ for all $z \in Z, \psi \in H^1(\Omega,H)$, then one has a global gauge symmetry of the system.

The symmetries that we discuss in this work are of this type. 
Observe that if $Y$ is invariant, then the unitary group $s_z:=(e^{izP_Y})_{z \in \mathbb R}$ 
defines a group such that $a(\psi)=a(S_z\psi)$ for all $z \in Z, \psi \in H^1(\Omega)$. That is, it is a gauge symmetry of the system.

To see what is the action of this group, observe that writing the $e^{izP_Y}$ as the exponential series $\sum_{n=0}^\infty \frac{(izP_Y)^n}{n!}$ yields
$$
e^{izP_Y}= e^{iz}P_Y + (\id -P_Y), \qquad z \in Z.
$$
Decompose now $\psi:=\psi_{Y} + \psi_{Y^\perp}$. Now, the gauge group acts as 
$$e^{izP_Y}: \psi \mapsto e^{iz}\psi_{Y} + \psi_{Y^\perp}, \qquad z \in Z, \psi \in H^1(\Omega,H).$$
In other words, it acts as the standard gauge group on the $Y$-component of the function and as the identity on the $Y^\perp$-component, see also~\cite{BolCarMug07}.

Finally, observe that also the symmetries considered in Chapter~\ref{hesse} were of this type.

\subsubsection{Section~\ref{sec:infintro}}
The Wolfgang's first network tutorial was written for my supervisor Wolfgang Arendt.
The approach we use to represent graphs and networks is mostly combinatorial. 
The reader interested in details and in general topics of graph theory should consult the monographs \cite{Die05}, or \cite{Vol96}, 
from which we have borrowed most of the notations.
Another approach on the theory of network equations it is based on the theory of $CW$-complexes, see~\cite{Whi49}-\cite{Whi49a}. 
It avoids combinatorial issues and it is, for instance, the approach used in~\cite{Cat97}.

\subsubsection{Section~\ref{sec:graphoper}}
The results in this section are unpublished.
The issues discussed in Section~\ref{sec:graphoper} are not new in the mathematical literature. 
An amusing introduction to the topic can be found in~\cite[Chap.~4]{Hal67}. 
There, the boundedness of infinite matrices in the space $\ell^2(\V)$ set is discussed. 
As Paul Halmos points out, \emph{``there are no elegant and usable necessary and sufficient conditions''} for characterizing the boundedness of such an operator. 
In Proposition~\ref{incidencebound} is an example of such a result for a restricted class of operators. It be compared with the well-known fact that the stochastic version of the adjacency matrix is a bounded operator in $\ell^2$ if and only if the graph is uniformly locally finite. 
We reformulate the result in Proposition \ref{incidencebound} for infinite matrices. To this end, assume that the set of edges and of vertexes have the same cardinality, i.e., both are infinite countable sets. 
Then fixing a numbering of edges and vertexes, i.e., identifying $\V$ and $\E$ with $\mathbb N$, shows that a stochastic, Boolean valued, infinite matrix is bounded in $\ell^2(\mathbb N)$ if and only if there is a uniform bound on the number of non-zero entries of each row. 
This result is similar to the classical result of Otto Toeplitz (see~\cite{Toe10}), that if $A$ is an operator in $\mathcal L(\ell^2(\V))$, then there exists an orthonormal basis such that the number of non-zero entries of each row is finite.

On the other side, one could ask whether the question of the boundedness of the incidence operator have been already discussed in the literature.
We refer to \cite{MohWoe89} for a complete survey of the standard result about discrete operators on graphs, and for boundedness issues to \cite{Var85}. 
There, the case of a transition matrix corresponding to a reversible Markov chain is discussed.
%\begin{prop}A Boolean valued (possibly infinite) matrix is reversible if and only if it is symmetric.\end{prop}
%\begin{proof}We recall that a Markov chain is reversible if $\exists (\lambda_i)_{i \in I}, \lambda_i > 0,$ such that for the transition matrix $\mathcal I$ one has 
%$$\lambda_i \iota_{ij}= \lambda_j \iota_{ji}, \qquad i,j\in I.$$ 
%So, if $\iota_{ij}=0$, then the above equation is solvable only for $\iota_{ji}=0$. If, conversely, $\iota_{ij} \neq 0$, then also $\iota_{ji} \neq 0$ and since $\mathcal I$ is Boolean valued one obtains $\iota_{ij}= \iota_{ji}=1$.\end{proof}
It should be mentioned that in \cite{Moh82} some results about compactness, boundedness and spectral properties of the adjacency operator are investigated

\subsubsection{Section~\ref{sec:operatordomain}}
The results in this section come from~\cite{CarMugNit07} for the case of a finite network.

In particular, the identification of the domain of $L^2$-realisations of the strongly coupled elliptic operators on networks did not attract the interest of researchers.
In~\cite{KosSch99}, however, some issues regarding the domain of the Laplacian are considered; a more general case is discussed there, since the Laplacian is not seen as the operator associated with a sesquilinear form.

A related question is the study of spectral properties of the Laplacian on a network.
This is mostly studied in a $C^2$-setting, see~\cite{Bel85}-\cite{Bel88} or in  $L^\infty$-setting for locally finite networks, see~\cite{BelLub04}-\cite{BelLub05}.

\subsubsection{Section~\ref{sec:systemsII}}
The results in this section are taken from~\cite{CarMugNit07}, with some generalisations and modifications.
In particular, in that paper symmetries were introduced by fixing an orthogonal projection  $P \in \mathcal L(\ell^2(I))$.
Starting from subspaces has the advantage of making immediately clear what is the connection with the classical meaning of symmetries.
In that work is also shown that for the Laplacian is equivalent to study symmetries for the heat or Schr\"odinger equation.
It should be mentioned that other forms of symmetries have been considered in connection with diffusion (or Schr\"odinger) equations on graphs.
In fact, if the graph possesses some symmetry group, then it is possible to exploit this feature in order to construct
counterexamples to the Kac's problem~\cite{Kac66} for graphs. 
This has been done in~\cite{Bel01}--\cite{GutSmi01}, and a systematic approach has been announced in a forthcoming paper~\cite{BanShaSmi06}.
For a discussion of the issue of symmetries in quantum graphs see also~\cite{BolCarMug07}.

\subsubsection{Section~\ref{sec:infinitesymmetries}}
The results in this section are unpublished.
As we already pointed out in the comments to Section \ref{sec:graphoper}, most of the authors have only considered spectral issues 
for the continuous Laplacian on infinite graphs, neglecting the investigation of the parabolic problem. 
The standard assumption is the uniform locally finiteness of the graph.
Observe that by Theorem \ref{infiniteirreducibility} the heat semigroups is irreducible on every locally finite, connected graph.

On the other hand, many authors considered the case of the the discrete Laplacian on infinite graphs.
Such considerations are interesting because of their connection with the theory of Markov chains and Riemannian manifolds.
For an introduction to the topic, see~\cite{Soa94}.

\subsubsection{Section~\ref{sec:treelayer}}
The results in this section come from~\cite{CarMugNit07}.
Direct graph theoretic characterisations of non-spectral properties of the parabolic equations are rare in the literature.
In~\cite{Cat98}-\cite{Cat99} heat kernels of the parabolic equation are derived; this derivation involves graph theoretic objects.
%In~\cite{Bel94} the spectrum of infinite periodic graphs is studied.

\subsubsection{Section~\ref{sec:general}}
The results in this section are unpublished.
Observe that it is possible to relate qualitative properties of the associated semigroup to functional analytic properties of the orthogonal projection $P_Y$. 
This is a work in progress with Delio Mugnolo and Olaf Post.
The most general boundary condition have been investigated in~\cite{KosSch99}. 
Although not all the boundary conditions considered in that work can be obtained as described in Section~\ref{sec:general},
we stress that our method allows us to characterise symmetry and qualitative properties of the semigroup coming from the form.

%% file: sec_sesquilinear.tex
\section{Sesquilinear forms on complex Hilbert spaces}\label{sec:sesquilinear}
We state here all results we used in the work. We also give the basic definitions. 
The reader interested in the theory of sesquilinear forms and analytic semigroups should consult the monographs~\cite{Are06},\cite{Dav89},\cite{Ouh04} 
and the survey~\cite{Are04}.
\begin{defin}[Sesquilinear mappings and forms]\label{franceform}
Consider Hilbert spaces $U,V,$ and a mapping $a: U \times V \to \mathbb C$.
\begin{enumerate}
\item The mapping $a$ is a \emph{sesquilinear mapping} if $a$ is linear in the first component and antilinear in the second one, i.e. if
$$
a(\lambda \psi + \mu \psi', \psi'')= \lambda a(\psi,\psi'') + \mu a(\psi', \psi''), \qquad \mbox{for all } \psi,\psi' \in U, \psi'' \in V
$$
and
$$
a( \psi'', \lambda \psi + \mu \psi')= \overline{\lambda} a(\psi,\psi'') + \overline{\mu} a(\psi', \psi''), \qquad \mbox{for all } \psi,\psi' \in V, \psi'' \in U.
$$
If $U=V$ we use the term \emph{sequilinear form} and we specify the Hilbert space writing $a=(a,V)$. 
\item With an abuse of notation, we denote the \emph{energy functional} $a: V \to \mathbb C$ associated with the form $(a,V)$ also $(a,V)$. 
The energy functional is defined by
$$
a(\psi):= a(\psi,\psi), \qquad \mbox{for all } \psi \in V.
$$
\item A sesquilinear mapping $a$ is \emph{continuous} if there exists $M \in \mathbb R$ such that
$$
a(\psi,\psi') \leq M \|\psi\|_U \|\psi'\|_V, \qquad \mbox{for all } \psi \in V.
$$
For a sesquilinear form this is equivalent to the existence of a $M \in \mathbb R$ such that
$$
a(\psi) \leq M\|\psi\|^2_V, \qquad \mbox{for all } \psi \in V.
$$
\item A sesquilinear form $(a,V)$ is \emph{coercive} if there exists $\alpha>0$ such that
$$
\re a(\psi) \geq \alpha \|\psi\|^2_V,\qquad \mbox{for all } \psi \in V.
$$
\item A sesquilinear form $(a,V)$ is \emph{accretive} if
$$
\re a(\psi) \geq 0, \qquad \mbox{for all } \psi \in V.
$$
\item Fix a second Hilbert space $H$, $V\hookrightarrow H$. 
A sesquilinear form $(a,V)$ is $H$\emph{-elliptic} if there exist $\alpha > 0$, $\omega \in \mathbb R$ 
such that the form $a_\omega:V\times V \to \mathbb C$ defined by
$$
a_\omega(\psi,\psi'):= a(\psi,\psi') + \omega(\psi,\psi')_H
$$
is coercive.
\end{enumerate}
\end{defin}
If we fix a second Hilbert space $H$ such that $V \hookrightarrow H$ densely, we can then canonically associate an operator $(A,D(A))$ on $H$ with the form $(a,V)$.
\begin{defin}
Let $(a,V)$ a continuous form. The operator $(A,D(A))$ associated with the form $(a,V)$ is defined on
$$
D(A):=\{\psi \in V : \exists f \in H,\mbox{ such that } a(\psi,\psi') = (f,\psi'), \mbox{ for all } \psi' \in V\}
$$
by
$$
A\psi := - f.
$$
\end{defin}
We are interested in studying evolution equations in Hilbert spaces which have the form of an abstract Cauchy problem
\begin{equation}\label{EE}
\left\{\begin{array}{rclr}
\dot{u}(t)&=& Au(t), & t \geq 0, \\
u(0)&=& u_0, & u_0 \in H.
\end{array}\right.
\end{equation}
Assume that the problem~\eqref{EE} is well-posed. 
Then, in analogy with vector-valued ordinary differential equations, 
we denote the solution $u(t)$ to~\eqref{EE} for the initial data $u_0$ by $e^{tA} u_0$. 
The family $(e^{tA})_{t\geq 0} \in \mathcal L(H)$ is called a \emph{semigroup of operators} (in the following: semigroup).
Reader interested in the general theory of semigroups should consult~\cite{EngNag00} and~\cite{AreBatHie01}.
There are many possible notions of continuity and differentiability for semigroups.
We focus on analytic semigroups and we refer to~\cite[Chapter~2]{Are06} for a short introduction in this topic.
\begin{thm}\label{WPform}
Assume that the operator $(A,D(A))$ is associated with a sesquilinear form $(a,V)$. 
Then, if the form $(a,V)$ is $H$-elliptic, then $(A,D(A))$ generates an analytic semigroup $(e^{tA})_{t \geq 0}$. 
In this case we write $\etasg:=(e^{tA})_{t\geq 0}$.
\end{thm}
One can characterise properties of $\etasg$ by properties of the form $a$.
As an example, we observe that accretivity and coercivity of the form already imply contractivity estimates for the semigroup.
\begin{thm}[Forms and analytic semigroups]
On an Hilbert space $H$ consider a continuous, $H-$elliptic sesquilinear form $(a,V)$ with constants $\alpha,\omega$. 
The following assertions hold.
\begin{enumerate} 
\item If the form $(a,V)$ is coercive, then there is $\epsilon >0 $ such that $\|e^{ta}\|_{\mathcal L(H)} \leq M e^{-\epsilon t}$, for all $t \geq 0$.
\item The estimate $\|e^{ta}\|_{\mathcal L(H)} \leq 1$ holds for all $t \geq 0$ if and only if $(a,V)$ is accretive.
\end{enumerate}
\end{thm}
An important tool for the analysis of qualitative properties of semigroups coming from forms is the following result due to El-Maati Ouhabaz, see \cite[~Thm. 2.1]{Ouh04}.
\begin{thm}[Ouhabaz's criterion]\label{ouh}
The closed convex subset $C \subset H$ is invariant under the action of the contractive semigroup $(e^{ta})_{t\geq 0}$ if and only if for the Hilbert space projection $P$ onto $C$ holds
$$
\psi \in V \Longrightarrow P\psi \in V,\quad \re a(P\psi, \psi-P\psi) \geq 0.
$$
\end{thm}
A consequence of Theorem \ref{ouh} is the following result.

\begin{cor}[Invariance criterion for closed subspaces]\label{ortho}
The closed linear subspace $Y$ is invariant under the action of the  semigroup $(e^{ta})_{t\geq 0}$ if and only if
\begin{enumerate}
\item The orthogonal projection $P_Y$ of $H$ onto $Y$ satisfies $\psi \in V \Rightarrow P_Y\psi \in V$.
\item $a(\psi,\psi')=0$ for all $\psi \in Y \cap V$, $\psi' \in Y^\perp \cap V$.
\end{enumerate}
\end{cor}

\begin{rem}
We will refer to the first condition in Corollary \ref{ortho} as to the \emph{admissibility condition} and to the second one as to the \emph{orthogonality condition}.
\end{rem}

For $\alpha \in [0,1)$ an interpolation space of order $\alpha$ between $V$ and $H$ is any linear space $H^\alpha$ such that
$$
V \hookrightarrow H^\alpha \hookrightarrow H, \quad \mbox{and} \quad \|\psi\|_{H^\alpha} \leq M \|\psi\|^\alpha_{V}\|\psi\|^{1-\alpha}_{H}, \quad \psi \in V.
$$
The following result is~\cite[Lemma~2.1]{Mug08}.

\begin{lemma}\label{delioperturb}
Let $a:V \times V \to \mathbb C$ a sesquilinear form. Consider two continuous sesquilinear mappings 
$a_1: H \times H_{\alpha} \to \mathbb C$, $a_2: H_{\alpha} \times H \to \mathbb C$, where $H_{\alpha}$ is an interpolation space between $H$ and $V$ of order $\alpha$.
Then $a$ is $H$-elliptic if and only if $a+a_1+a_2$ is $H$-elliptic.
\end{lemma}

\emph{Cosine operator functions} are defined analogously to semigroups of operators.
A cosine operator function is a mapping ${\mathrm Cos}: \mathbb R \to \mathcal L(H) $ such that
\begin{equation}\label{cosinedef}
{\mathrm Cos}(0)=Id, \quad
2{\mathrm Cos}(t){\mathrm Cos}(s)={\mathrm Cos}(t+s)+{\mathrm Cos}(t-s), \quad t,s\in \mathbb R.
\end{equation}
We do not enter the details of the theory of cosine operator functions, for which we refer to~\cite[Chapter~10]{AreBatHie01}.
We only mention that it is possible to define a \emph{generator} $(A,D(A))$ of a cosine operator function. If we define
\begin{equation*}
u(t):= {\mathrm Cos}(t) u_0 + \int_0^t {\mathrm Cos} (s) w_0 ds,
\end{equation*}
then $u$ is a weak solution of the second-order abstract Cauchy problem
$$
\left\{
\begin{array}{rclr}
 \ddot{u}(t)&=& Au(t), & t \geq 0, \\
u(0) & =& u_0, & u_0 \in H, \\
\dot{u}(0) &=& w_0,&  w_0 \in H.
\end{array}\right.
$$
Conversely, if $(A,D(A))$ is a closed operator such that the abstract Cauchy problem is weakly well-posed, then $(A,D(A))$ is the generator of a cosine function.
In the case of an operator $(A,D(A))$ associated with a continuous, elliptic, sesquilinear form $(a,V)$
the following result holds, see~\cite[p.~212]{Haa06}, or~\cite{Cro04}.
\begin{prop}\label{crouzeix}
Assume that 
$$
|\im a(u)| \leq M \|u\|_H\|u\|_V.
$$
Then $A$ generates a cosine operator function.
\end{prop}

%% file: sec_zusammen.tex
\chapter*{Zusammenfassung in deutscher Sprache}

Wir beginnen mit der Beobachtung, dass eine auf der direkten Summe von Hilbertr\"aumen
$$V:= \bigoplus_{i \in I} V_i$$
definierte Sesquilinearform $(a,V)$ eine Matrixdarstellung zul\"asst, die durch
$$
a_{ij}(\psi,\psi'):= a(\mathbb 1_j \otimes \psi , \mathbb 1_i \otimes \psi' )
$$
hergestellt wird.

Anhand dieser \"Uberlegung ist es m\"oglich, eine linear algebraische Theorie f\"ur auf direkten Summen von Hilbertr\"aumen definierte Sesquilinearformen entwickeln,
um Stetigkeit, Elliptizit\"at und \"ahnliche Eigenschaften zu charakterisieren, oder um wenigstens hinreichende Bedingungen beweisen zu k\"onnen.

Wenn $V_i=V_j$ f\"ur alle $i,j\in I$, kann man diese linear algebraische Theorie verwenden, 
um \emph{Symmetrieeigenschaften} der Halbgruppe $\etasg$ assoziert mit der Sesquilinearform $(a,V)$ in einem Hilbertraum $H$ zu untersuchen, 
die sich als Invarianz gewisser Teilr\"aume von $H$ formulieren lassen.

Der Schl\"ussel dazu ist die Tatsache, dass das Kriterion von Ouhabaz f\"ur die Invarianz abgeschlossener, konvexer Mengen f\"ur Teilr\"aume eine besonders einfache Form annimmt.
Um das zu sehen, bezeichne man $P_Y$ die Hilbertraum Projektion von $H$ auf den abgeschlossenen Unterraum $Y$.
Somit, gilt genau dann
$e^{ta}\mathcal Y \subset \mathcal Y$ f\"ur alle $t \geq 0$ 
wenn
\begin{equation}\label{intro3}
P_\mathcal Y V \subset V
\end{equation}
und
\begin{equation}\label{intro4}
a(\psi,\psi')=0, \qquad \psi \in \mathcal Y, \psi' \in \mathcal Y^\perp.
\end{equation}
Es zeigt sich, dass f\"ur Teilr\"aume, die mit Symmetrieeigenschaften assoziert sind, die erste Bedingung immer erf\"ullt ist, wenn der Formbereich $V$ 
eine direkte Summe ist.
Die Invarianz von $Y$ ist dann \"aquivalent zur zweiten Bedingung, und diese kann durch Eigenschaften der einzelnen Abbildungen $a_{ij}$ charakterisiert werden.

Wenn aber Kopplungsterme in der Definition von $V$ vorhanden sind, das heisst, wenn $V \subset \bigoplus_{i \in I} V_i$ und $V$ kein Ideal ist,
wird es aufw\"andiger eine linear algebraische Theorie zu entwickeln. Eine L\"osung hierf\"ur ist Spezialklassen von Formbereichen und Kopplungstermen zu untersuchen.

Ein Beispiel daf\"ur sind Netzwerkgleichungen. In diesem Fall, der Formbereich kann als
$$
V:=\{\psi \in \bigoplus_{i \in I} H^1(0,1): \psi(0)\oplus \psi(1)  \in Y\}
$$
geschrieben werden, wobei $Y$ ein geeigneter Unterraum von $\ell^2(I) \bigoplus \ell^2(I)$ ist.

Auch in diesem Fall ist es m\"oglich die Symmetrieeigenschaften der assoziierten Halbgruppe $\etasg$ systematisch zu studieren.

F\"ur die W\"armeleitungsgleichung, insbesondere, ist die Orthogonalit\"at~\eqref{intro4} automatisch erf\"ullt, 
so dass die Invarianz \"aquivalent zur Zul\"assigkeit~\eqref{intro4} ist.

Die Herausforderung ist Zul\"assigkeit graphentheoretisch zu charakteririsieren, und, m\"oglicherweise, einige Aussagen auf allgemeinere Situationen zu verallgemeinern.

\medskip
%Im Kapitel~\ref{hesse} entwickeln wir eine allgemeine Theorie f\"ur Sesquilinearformen, die auf der direkten Summe von Hilbertr\"aumen definiert werden und
%wir zeigen wie qualitative Eigenschaften der mit der Sesquilinearform $(a,V)$ in einem geeigneten Hilbertraum assozierte Diffusionshalbgruppe $\etasg$
%durch die Analyse der Matrixdarstellung der Form untersucht werden k\"onnen.

%Im Kapitel~\ref{network} wenden wir diese Theore auf den Fall von Netwzerkgleichungen an.
%Die Techniken, die wir im Kapitel~\ref{hesse} entwickelt haben, k\"onnen hier angewendet werden,
%um Symmetrieeigenschaften zu charakterisieren, aber zus\"atzliche Argumente werden ben\"otigt.

Wir beginnen mit der Beobachtung, dass eine Sesquilinearform $(a,V)$, 
die auf einem Hilbertraum $V:=\bigoplus_{i \in I} V_i$ definiert ist, als eine Matrix $(a_{ij})_{i,j \in I}$ von Sesquilinearabbildungen betrachtet werden kann.
Von dieser Darstellung kann man hinreichende Bedingungen f\"ur unterschiedliche Eigenschaften von $(a,V)$ angeben, die linear algebraischer oder elementar funktionalanalytischer Art sind. Dieses \"Uberlegung wirde im Abschnitt~\ref{sec:findimarg} durchgef\"uhrt.

Solche Argumente reichen oft nicht aus, um einfache Eigenschaften sowie Stetigkeit und Koerzivit\"at unendlicher Formmatrizen zu charakterisieren.
Im Abschnitt~\ref{sec:systemsI} betrachten wir ein unendliches, stark gekoppeltes System und charakterisieren Stetigkeit und Koerzivit\"at.

Im Abschnitt~\ref{sec:gencoercivity} wenden wir uns dem Thema der Koerzivit\"at im allgemeinen Fall zu und zeigen, dass es ausreichend ist, die endlichen Untermatrizen von $a$ zu untersuchen. Eine Charakterisierung der Koerzivit\"at zwei dimensionaler Formmatrizen duch Eigenschaften der einzelnen Abbildungen wird bewiesen.

Im Abschnitt~\ref{sec:operator} beginnen wir die Untersuchung von Evolutionsgleichungen und versuchen den Definitionsbereich des Operators zu identifizieren.
Wir setzen unsere Untersuchungen im Abschnitt~\ref{sec:evoleq} fort, in dem wir Wohlgestelltheit-\, und Kontraktivit\"atseigenschaften von $(\etasg)$ untersuchen.
Beispielweise, beweisen wir in Theorem~\ref{positivity}, dass die Halbgruppe $\etasg$ genau dann positiv ist,
wenn die Halbgruppen, die von der Sesquilinearformen auf der Hauptdiagonal erzeugt werden, positiv sind, und alle Abbildungen ausserhalb der Hauptdiagonale negativ.

Wir untersuchen im Abschnitt~\ref{sec:symmetriesmatrix} Symmetrieeigenschaften und Theorem~\ref{symmetries} charakterisiert vollst\"andig diese Eigenschaften.

Die Abschnitte~\ref{sec:wave} und~\ref{sec:dynamic} enthalten zwei Anwendungsbeispiele.
Im ersten betrachten wir eine stark ged\"ampfte Wellengleichung und im zweiten eine W\"armeleitungsgleichung mit dynamischen Randbedingungen.
In beiden F\"allen erhalten wir Wohlgestelltheit f\"ur eine gro{\ss}e Auswahl an Parametern dank Satz~\ref{perturbinter}. Darauf folgend untersuchen wir qualitative Eigenschaften der Gleichungen.

Schliesslich, f\"uhren wir im Abschnitt~\ref{sec:nondiag} das Thema der nicht diagonalen Formbereiche ein.
Netzwerkgleichen sind wahrscheinlich das wichtigste Beispiel von nicht diagonalen Formbereichen, und die werden im Kapitel~\ref{network} untersucht.

Die Ergebnisse des Kapitels~\ref{hesse} sind im Abschnitt~\ref{sec:histremI} diskutiert.
Insbesondere, wir zeigen die Vorteile eines Matrixformalismus f\"ur Sesquilinearformen gegen\"uber einem Formalismus f\"ur Operatoren.
Im gleichen Abschnitt beschreiben wir die Geschichte der verschiedenen Fragsestellungen des Kapitels.

Wir beginnen in Kapitel~\ref{network} mit einer informellen Erkl\"arung der Idee, die hinter der Benutzung von Sesquilinearformen f\"ur den Laplace Operator auf Netzwerken steht.
Der Abschnitt~\ref{sec:infintro} stellt eine komprimierte Erl\"auterung über partielle Integration in Netzwerken dar.

Da wir an unendliche Netzwerke im $L^2$-Fall interessiert sind, zeigen wir im Abschnitt~\ref{sec:graphoper}, dass alle Definitionen, die f\"ur endliche Netzwerke \"ublich sind, auch f\"ur unendliche Netzwerke formuliert werden k\"onnen.
Insbesondere untersuchen wir operatorentheoretische Eigenschaften der Inzidenzmatrizen.
Nach diesen Untersuchungen definieren wir in~\eqref{formdomain}-\eqref{networkform} den Formbereich und die Wirkung der Form $(a,V)$.

Im Abschnitt~\ref{sec:systemsII} charakterisieren wir Symmetrieeigenschaften f\"ur Diffusionsysteme in Netzwerken.
Dieser Abschnitt ist die Fortsetzung des Abschnittes~\ref{sec:systemsI}.
In Theorem~\ref{charadmiss} charakterisieren wir die Zul\"assigkeit von Projektionen assoziert mit Symmetrieeigenschaften vollst\"andig.

Die beiden Abschnitte~~\ref{sec:infinitesymmetries} und~\ref{sec:treelayer} sind eine Anwendung dieser Resultate.
Im Ersten wenden wir unsere Aufmerksamkeit der Irreduzibilit\"at der Halbgruppe $\etasg$ im Fall unendlichere Netzwerke zu, und charakterisieren die Netzwerke, so dass die Halbrguppe $\etasg$ irreduzibel ist. Im Zweiten untersuchen wir Symmetrieeigenschaften einiger bestimmter Klassen von Netzwerken.

Im Abschnitt~\ref{sec:general} verlassen wir die Netzwerke und zeigen, dass alle Systeme von Diffusionsgleichungen auf $H^1(0,1)$
als Diffusionsysteme in Netzwerken interpretiert werden k\"onnen, wenn die Randbedingungen geeignete Eigenschaften erf\"ullen.

Der Begriff der Symmetrie wird in seiner physikalischer Bedeutung im Abschnitt~\ref{sec:histremII} erl\"autert,
Wir unterscheiden zwischen \emph{Raum-} und \emph{Eichsymmetrien} und wir zeigen, dass die Symmetrien, die wir untersuchen, Eichsymmetrien in der urspr\"unglichen
physikalischen Bedeutung sind.
Abschlie{\ss}end diskutieren wir die Geschichte der unterschiedliche Fragestellungen des Kapitels.

%% file: sec_leben.tex
\pagestyle{empty}
\begin{tabular}{r|l}
& {\bf \large Lebenslauf}\\[2cm]
& Stefano Cardanobile \\[1cm]
21.09.1980 & geboren in Bari (Italien)\\[0.3cm]
1987 - 1991 & Scuola elementare in Bari\\[0.3cm]
1991 - 1994 & Scuola media in Bari\\[0.3cm]
1994 - 1998 & Liceo classico in Bari\\[0.3cm]
1998 - 1999 & Studium der Maschinenbau \\[0.3cm]
&\hspace{0.3cm} an dem Politecnico di Bari\\[0.3cm]
1999 - 2001 & Studium der Mathematik \\[0.3cm]
&\hspace{0.3cm}an der Universit\`a di Bari\\[0.3cm]
2001 - 2002 & Studienaufenthalt \\[0.3cm]
&\hspace{0.3cm}an der Universit\"at T\"ubingen\\[0.3cm]
2002 - 2005 & Studium der Mathematik\\[0.3cm]
&\hspace{0.3cm}an der Universit\"at T\"ubingen\\[0.3cm]
2005 & Diplom in Mathematik\\[0.3cm]
2005-2008 & Stipendiat im Promtionskolleg\\[0.3cm]
&\hspace{0.3cm}Mathematische Analyse\\[0.3cm] 
&\hspace{0.3cm}von Evolution, Information und Komplexit\"at\\[0.3cm]
&\hspace{0.3cm}an der Universit\"at Ulm.
\end{tabular}

%% file: sec_publikationen.tex
{\bf \large Publikationsliste}\\[1cm]
\begin{enumerate}[1.]
\item \emph{Analysis of a FitzHugh--Nagumo model of a neuronal network with excitatoric semilinear node conditions}, mit Delio Mugnolo. 
Erschienen in ``Mathematical Methods in the Applied Sciences'', Vol. 30, no. 18, pp. 2281-2308, 2007.
\item \emph{Well-posedness and symmetries of strongly coupled network equations}, mit Delio Mugnolo und Robin Nittka. 
Erschienen in ``Journal of Physics A: Mathematical and Theoretical'', 41(5):055102 (28pp), 2008.
\item \emph{Qualitative properties of coupled parabolic systems of evolution equations}, mit Delio Mugnolo. 
Im Druck in ``Annali della Scuola Normale Superiore - Classe di Scienze'', Pisa.
\item \emph{Symmetries in Quantum Graphs}, mit Jens Bolte, Delio Mugnolo und Robin Nittka. 
Im Druck in ``Mathematical Analyisis of Evolution, Information and Complexity'', editors Wolfgang Arendt und Woflgang Schleich, Wiley, Berlin.
\item \emph{Investigation of input-output gain in dynamical systems for neural information processing}, 
mit Michael Cohen, Silvia Corchs, Delio Mugnolo und Heiko Neumann. 
Im Druck in ``Mathematical Analyisis of Evolution, Information and Complexity'', editors Wolfgang Arendt und Woflgang Schleich, Wiley, Berlin.
\item \emph{Relating Simulation and Modelling of Neural Networks}, mit Heiner Markert, Delio Mugnolo, Günther Palm und Friedhelm Schwenker.
Im Druck in ``Mathematical Analyisis of Evolution, Information and Complexity'', editors Wolfgang Arendt und Woflgang Schleich, Wiley, Berlin.
\end{enumerate}